\documentclass[11pt]{article}
\usepackage{amssymb}
\usepackage{amsmath}
\usepackage{mathrsfs}
\usepackage{graphics}
\usepackage{graphicx}
\usepackage{subfigure}
\usepackage[T1]{fontenc}
\usepackage{latexsym,amssymb,amsmath,amsfonts,amsthm}\usepackage{txfonts}
\newcommand{\be}{\begin{eqnarray}}
\newcommand{\ben}{\begin{eqnarray*}}
\newcommand{\en}{\end{eqnarray}}
\newcommand{\enn}{\end{eqnarray*}}

\newtheorem{theorem}{Theorem}[section]
\newtheorem{lemma}{Lemma}[section]
\newtheorem{prp}[theorem]{Proposition}

\newtheorem{dfn}{Definition}[section]

\newtheorem{remark}{Remark}

\begin{document}
\renewcommand{\theequation}{\arabic{section}.\arabic{equation}}
\begin{titlepage}
\title{\bf Markov Selection and $\mathcal{W}$-strong Feller for 3D
Stochastic Primitive Equations
}
\author{Zhao Dong$^{1,}$,\ \ Rangrang Zhang$^{1,}$\\
{\small $^1$ Institute of Applied Mathematics,
Academy of Mathematics and Systems Science, Chinese Academy of Sciences,}\\
{\small No. 55 Zhongguancun East Road, Haidian District, Beijing, 100190, P. R. China}\\
({\sf dzhao@amt.ac.cn},\ {\sf rrzhang@amss.ac.cn} )}
\date{}
\end{titlepage}
\maketitle

%\vspace{.2in} \noindent
%Proposed running head: {\bf something about it}\\
%\vspace{.1in}

%\noindent
%Manuscript Correspondence Address:\\
%Professor Bo Zhang\\
%Institute of Applied Mathematics\\
%Academy of Mathematics and Systems Science\\
%Chinese Academy of Sciences\\
%Beijing 100080, China\\
%Email: b.zhang@amt.ac.cn\\
%Telephone: +86 10 6265 1358\\
%Fax:       +86 10 6254 1689
%
%\newpage
%\begin{center}
%\title{\Large\bf Newton iteration of inverse scattering problem}
%\end{center}

%\vspace{.2in}

\begin{abstract}
This paper studies some analytical properties of weak solutions of 3D stochastic primitive equations with periodic boundary conditions. The martingale problem associated to this model is shown to have a family of solutions satisfying the Markov property, which is achieved by means of an abstract selection principle. The Markov property is crucial to extend the regularity of the transition semigroup from small times to arbitrary times. Thus, under a regular additive noise, every Markov solution is shown to have a property of continuous dependence on initial conditions, which follows from employing the weak-strong uniqueness principle and the Bismut-Elworthy-Li formula.

\textbf{Key words}: Markov selection; $\mathcal{W}$-strong Feller;
%Moreover, uniqueness of invariant measures are stated clearly and t

%\vspace{.2in} {\bf Keywords:}.
\end{abstract}

\section{Introduction }
The primitive equations (PEs) derived by Boussinesq approximation are a basic model in the study of large oceanic and atmospheric dynamics. These systems form the analytical core of the most advantaged general circulation models. For this reason and due to their challenging nonlinear and anisotropic structure, the PEs have recently received considerable attention from the mathematical community.

The mathematical study of the PEs originated in a series of articles by J.L. Lions, R. Temam, and S. Wang in the early 1990s \cite{L-T-W-1,L-T-W-2,L-T-W-3,L-T-W-4}. They set up the mathematical framework and showed the global existence of weak solutions. For the existence and uniqueness of strong solution, many works are concerned on it. For example, C. Hu, R. Temam and M. Ziane proved the global existence and uniqueness of strong solutions to the viscous primitive equations in thin domains for a large set of initial data whose size depends on the thickness of the domain in \cite{H-T-Z}. In \cite{GG-MN}, F. $\rm Guill\acute{e}n-Gonz\acute{a}lez$, N. Masmoudi and M.A. Rodriguez-Bellido showed the local existence and uniqueness of strong solutions to the viscous primitive equations for any initial data. C. Cao and E.S. Titi developed a beautiful approach to dealing with the $L^6$-norm of the fluctuation $\tilde{v}$ of horizontal velocity and obtained the global well-posedness for the 3D viscous primitive equations in  \cite{C-T-1}. For the uniqueness of weak solutions, in \cite{L-T}, J. Li and E.S. Titi established some conditional uniqueness of weak solutions to the viscous primitive equations under periodic boundary conditions, and they proved the global existence and uniqueness of weak solutions with the initial data taken as small $L^{\infty}$ perturbations of functions in the space $X=\Big\{v\in (L^6(\mathcal{O}))^2| \partial_z v\in (L^2(\mathcal{O}))^2\Big\}$.

%Despite these great successful developments for the deterministic
%primitive equations, introducing uncertainty in `exact' model
%equations is reasonable and necessary in the study of the primitive
%equations of the large scale ocean. The introducing of stochastic
%processes is aimed at accounting for a number of uncertainties and
%errors.
For the primitive equations in random case, many authors paid attention to it. In \cite{Guo}, B. Guo and  D. Huang obtained
the existence of universal attractor of strong solution under the assumptions that the momentum
equation is driven by an additive stochastic forcing and the
thermodynamical equation is driven by a fixed heat source. A. Debussche, N. Glatt-Holtz, R. Temam and M. Ziane
established the global well-posedness of strong solution, when the primitive
equations are driven by multiplicative random noises in \cite{D-G-T-Z}. For the ergodicity, in \cite{RR}, the authors obtained the existence of global weak solutions, and also obtained the exponential mixing property for the weak solutions which are limits of spectral Galerkin approximations of 3D stochastic primitive equations driven by regular multiplicative noise. For a special case that the stochastic primitive equations are in two space dimensions with small linear multiplicative noise, H. Gao and C. Sun obtained a Wentzell-Freidlin type large deviation principle for by weak convergence method in \cite{G-S}, where they omit the spatial variable $y$ and only take $(x,z)$ into account. Furthermore, they established the Hausdorff dimension of the global attractor is finite in \cite{G-S-1}. When the primitive equations are driven by an infinite-dimensional additive fractional noise with Hilbert-space-valued, G. Zhou obtained the existence of random attractor in \cite{Z}.

As we know, both in deterministic and stochastic case, the uniqueness of weak solutions is an important open problem, which results in many properties of weak solutions disappear. Thus, in order to have a deeper understanding of weak solutions and have some development on their uniqueness, it's natural to explore more properties of them. This article presents a step in this direction. We establish that there exists an almost sure Markov family of the primitive equations forced by multiplicative noise. Furthermore, we obtain that every Markov solution has a property of continuous dependence on the initial conditions ($\mathcal{W}$-strong Feller) if the primitive equations are driven by a regular additive noise.
 In comparison with \cite{RR}, we stress that the main improvement of our paper is that the $\mathcal{W}$-strong Feller is valid for all Markov solutions and not restricted to solutions which are limits of Galerkin approximations. Moreover, the conditions on the noise here is much weaker than those in \cite{RR}.

When uniqueness of weak solutions is open, Markov property has no direct meaning but a natural question is the existence of a Markov selection. A sufficient condition for the existence of almost sure Markov selections was provided by B. Goldys, M. R\"{o}ckner and X. Zhang in \cite{B-M-X}, where they dealt with an abstract stochastic evolution equations. Here, we apply this sufficient condition to our equations (\ref{eq5-1})-(\ref{eqq20}) and obtain
\begin{theorem}\label{thm-2}
Under \textbf{Hypothesis H0}, there exists an almost sure Markov family $(P_x)_{x\in H}$ of  (\ref{eq5-1})-(\ref{eqq20}).
\end{theorem}
The definition of weak solution of (\ref{eq5-1})-(\ref{eqq20}) is in Sect 4.2.

The important part of this paper is to investigate the continuity with respect to the initial conditions (strong Feller property) for the Markov family $(P_x)_{x\in H}$.
To achieve this, $\mathcal{W}$-strong Feller is considered which is weaker than strong Feller in $H$ when $\mathcal{W}$ is a subspace of $H$.
In the past two decades, there are several works concerned on $\mathcal{W}$-strong Feller for stochastic evolution equations. In particular, F. Flandoli and M. Romito established an abstract criterion Theorem 5.4 to obtain $\mathcal{W}$-strong Feller property for Markov selections of 3D Stochastic Navier-Stokes equations in \cite{F-M}. It says that if a Markov process coincides on a small positive random time with a strong Feller process, then it is strong Feller itself. The idea behind this is to use an approximation by a regularised problem, which has itself strong Feller solutions. For the concrete proof, two key points are needed: weak-strong uniqueness principle and the Bismut-Elworthy-Li formula. It's worth mentioning that this technique is usually applied to handle locally Lipschitz nonlinearities in stochastic equations.
%As Bismut-Elworthy-Li formula is used, which restricts us to deal with a regular noise.
%In \cite{F-M}, the authors obtained that Markov solutions of stochastic 3D Navier-Stokes is strong Feller in $\mathcal{W}=D(A^{\theta(\alpha_0)})$, where
%\begin{eqnarray*}
%\theta(\alpha_0)=\left\{
%                   \begin{array}{ll}
%                     \frac{1}{2}+\frac{\alpha_0}{2}, & \rm{if}\  0<\alpha_0\leq \frac{1}{2} ;\\
%                     \frac{1}{4}+\alpha_0, &  \rm{if}\  \alpha_0>\frac{1}{2}.
%                   \end{array}
%                 \right.
%\end{eqnarray*}
To study the strong Feller property of our equations, we will follow the idea of Theorem 5.4 in \cite{F-M}. Firstly, we introduce an auxiliary Ornstein-Uhlenbeck process $Z$ which is a stationary ergodic solution to a stochastic Stokes equation, then, the approximation process $X^{(R)}_t$ is obtained. To achieve our goal, two steps are needed: the approximation process $X^{(R)}_t$ coincides with the original process on a small positive random time interval and $X^{(R)}_t$ is $\mathcal{W}$-strong Feller. The second step is more challenging, where we have to control the power of $Z$ to be less than 2 as we want to apply the Fernique's theorem to $Z$. We overcome this difficulty by making use of $\Lambda \frac{\partial Z}{\partial z}$, which is a key idea in our proof, in that case, $Z\in C([0,T]; D(A))$ is needed, then the noise has to be chosen as $A^{-\frac{5}{4}-\varepsilon_0}$ with $\varepsilon_0>0$. Thus, the corresponding $\mathcal{W}$ is equal to $D(A^{\frac{3}{4}+\varepsilon_0})$. In comparison with 3D Navier-Stokes equations, the regularity of the noise here is required to be higher than $A^{-\frac{11}{12}-\varepsilon_0}$ with $\varepsilon_0>0$ for 3D Navier-Stokes equations because of strong nonlinear terms $(\int^z_{-1}\nabla_H \cdot vdz')\frac{\partial v}{\partial z}$ and $(\int^z_{-1}\nabla_H \cdot vdz')\frac{\partial T}{\partial z}$ in the primitive equations. Also, the advective structure of the primitive equations leads to a delicate asymmetry in the nonlinear terms, which requires a more refined calculation.
\begin{theorem}\label{thm-4}
Assume \textbf{Hypothesis H1} holds. Let $(P_x)_{x\in H}$ be the Markov solution of (\ref{eq5-1})-(\ref{eqq20}) and $(P_t)_{t\geq 0}$ be the associated operators on $B_b(H)$ defined as (\ref{eq-30}), then $(P_t)_{t\geq 0}$ is $\mathcal{W}$-strong Feller.
\end{theorem}
The definition of $\mathcal{W}$-strong Feller is defined in Sect 5.1.
%Here, $\mathcal{W}=D(\Lambda^{1+2\varepsilon_0})$ is a Sobolev space, which is a little stronger than $H^1$.

\vskip0.3cm
This paper is organized as follows. In Sects. 2 and 3, we introduce the 3D stochastic primitive equations and make formulation of those equations. The abstract Markov selection principle and concrete proof are given in Sect.4.  Finally,  $\mathcal{W}$-strong Feller is proved in Sect. 5.
%Finally, we make some further explorations on the invariant measures in Sect. 7.

\section{Preliminaries}
The 3D stochastic primitive equations of the large-scale ocean under a stochastic forcing, in a Cartesian system, are written as
\begin{eqnarray}\label{eq-1}
\frac{\partial v}{\partial t}+(v\cdot \nabla_H)v+\theta\frac{\partial v}{\partial z}+f{k}\times v +\nabla_H P - \Delta v &=&\sigma_1(v,T)\frac{dW_1}{dt},\\
\label{eq-2}
\partial_{z}P+T&=&0,\\
\label{eq-3}
\nabla_H\cdot v+\partial_{z}\theta&=&0,\\
\label{eq-4}
\frac{\partial T}{\partial t}+(v\cdot\nabla_H)T+\theta\frac{\partial T}{\partial z}- \Delta T&=&\sigma_2(v,T)\frac{dW_2}{dt},
\end{eqnarray}
where the horizontal velocity field $v=(v^{(1)},v^{(2)})$, the three-dimensional velocity field\ $(v^{(1)},v^{(2)},\theta)$, the temperature\ $T$ and the pressure\ $P$ are unknown functions. $f$ is the Coriolis parameter. ${k}$ is vertical unit vector. Set $\nabla_H=(\partial x,\partial y)$  to be the horizontal  gradient operator and $\Delta=\partial^{2}_{x}+\partial^{2}_{y}+\partial^{2}_{z}$ to be the three dimensional Laplacian. $W_1$ and $W_2$ are two independent cylindrical Wiener processes on $H_1$ and $H_2$, respectively. $H_1$ and $H_2$ will be defined in Sect. 3.

%The viscosity and the heat diffusion operators $L_1$ and $L_2$ are given by
%\[
%L_1 v=-\frac{1}{Re_1}\Delta v -\frac{1}{Re_2}\frac{\partial^2 v}{\partial z^2},
%\]
%\[
%L_2 T=-\frac{1}{Rt_1}\Delta T -\frac{1}{Rt_2}\frac{\partial^2 T}{\partial z^2},
%\]
%where $Re_1$, $Re_2$ are positive constants representing the horizontal and vertical Reynolds numbers, respectively, and $Rt_1$, $Rt_2$ are positive constants which stand for the horizontal and vertical heat diffusivity, respectively.
%
%\textbf{In the whole article, we assume that the constants $Re_1$, $Re_2$, $Rt_1$, $Rt_2$ are all equal to 1. Without this assumption, we can also obtain the results of this paper. }

The spatial variable $(x,y,z)$ belongs to $\mathcal{M}:= \mathbb{T}^2\times (-1,0)$. For simplicity of the presentation, all the physical parameters (height, viscosity, size of periodic box) are set to 1.

Refer to \cite{C-L-T}, the boundary value conditions for (\ref{eq-1})-(\ref{eq-4}) are given by
\begin{eqnarray}\label{eqq1}
v, \ \theta\ and \ T\ are\ periodic\ in\ x\ and\ y,\\
\label{eqq2}
(\partial_{z}v,\theta)\mid _{z=-1,0}=(0,0),\ T\mid_{z=-1}=1, \ T\mid_{z=0}=0.\\
\label{eqq3}
(v,T)\mid_{t=0}=(v_0,T_0).
\end{eqnarray}
Replacing $T$ and $P$ by $T+z$ and $P-\frac{z^2}{2}$, respectively, then (\ref{eq-1})-(\ref{eq-4}) with (\ref{eqq1})-(\ref{eqq3}) is equivalent to the following system
\begin{eqnarray}\label{eqq4}
\frac{\partial v}{\partial t}+(v\cdot \nabla_H)v+\theta\frac{\partial v}{\partial z}+f{k}\times v +\nabla_H P -\Delta v &=&\sigma_1(v,T+z)\frac{dW_1}{dt},\\
\label{eqq5}
\partial_{z}P+T&=&0,\\
\label{eqq6}
\nabla_H\cdot v+\partial_{z}\theta&=&0,\\
\label{eqq7}
\frac{\partial T}{\partial t}+(v\cdot\nabla_H)T+\theta(\frac{\partial T}{\partial z}+1)-\Delta T&=&\sigma_2(v,T+z)\frac{dW_2}{dt},
\end{eqnarray}
subject to the boundary and initial conditions
\begin{eqnarray}\label{eqq8}
v, \ \theta\ and \ T\ are\ periodic\ in\ x\ and\ y,\\
\label{eqq9}
(\partial_{z}v,\theta)\mid _{z=-1,0}=(0,0),\ T\mid_{z=-1,z=0}=0,\\
 \label{eqq10}
(v,T)\mid_{t=0}=(v_0,T_0).
\end{eqnarray}
Here, for simplicity, we still denote by $T_0$ the initial temperature in (\ref{eqq10}), though it is now different from that in $(\ref{eqq3})$.
%Here, , we still denote $(v},{T},{P},{\theta})$ by $(v,T,P,\theta)$.

Inherent symmetries in the equations show that the solution of the primitive equations on $\mathbb{T}^2\times (-1,0)$ with boundaries (\ref{eqq8})-(\ref{eqq10}) may be recovered by solving the equations with periodic boundary conditions in $x, y $ and $z$ variables on the extended domain $\mathbb{T}^2\times (-1,1):= \mathbb{T}^3$, and restricting to $z\in (-1,0)$.

To see this, consider any solution of (\ref{eqq4})-(\ref{eqq7}) with boundaries (\ref{eqq8})-(\ref{eqq10}), we perform that
\begin{eqnarray*}
v(x,y,z)&=&v(x,y,-z), \ \rm{for}\ (x,y,z)\in \mathbb{T}^2\times (0,1),\\
T(x,y,z)&=&-T(x,y,-z), \ \rm{for}\ (x,y,z)\in \mathbb{T}^2\times (0,1),\\
P(x,y,z)&=&P(x,y,-z), \ \rm{for}\ (x,y,z)\in \mathbb{T}^2\times (0,1),\\
\theta(x,y,z)&=&-\theta(x,y,-z),\ \rm{for}\ (x,y,z)\in \mathbb{T}^2\times (0,1).
\end{eqnarray*}
 We also extend $\sigma_1$ in the even fashion and  $\sigma_2$ in the odd fashion across $\mathbb{T}^2\times \{0\}$. Hence, we consider the primitive equations on the extended domain $\mathbb{T}^3=\mathbb{T}^2\times (-1,1)$,
 \begin{eqnarray}\label{eqq14}
\frac{\partial v}{\partial t}+(v\cdot \nabla_H)v+\theta\frac{\partial v}{\partial z}+f{k}\times v +\nabla_H P -\Delta v &=&\phi({v},{T})\frac{dW_1}{dt},\\
\label{eqq15}
\partial_{z}P+T&=&0,\\
\label{eqq16}
\nabla_H\cdot v+\partial_{z}\theta&=&0,\\
\label{eqq17}
\frac{\partial T}{\partial t}+(v\cdot\nabla_H)T+\theta(\frac{\partial T}{\partial z}+1)-\Delta T&=&\varphi({v},{T})\frac{dW_2}{dt},
\end{eqnarray}
 subject to the boundary and initial conditions
 \begin{eqnarray}\label{eqq11}
v, \ \theta\ , P\ and \ T\ are\ periodic\ in\ x\ ,\ y,\ z,\\
\label{eqq12}
v\ and\ P\ are\ even\ in \ z, \  \theta\ and\ T\ are\ odd\ in\ z,\\
\label{eqq13}
(v,T)\mid_{t=0}=(v_0,T_0),
\end{eqnarray}
where
\[
\phi(v,T)=\sigma_1(v,T+z)\quad \varphi({v},{T})=\sigma_2(v,T+z).
\]
Because of the equivalent of the above two kinds of boundary and initial conditions, we consider, throughout this paper, the system (\ref{eqq14})-(\ref{eqq13}) defined on $\mathbb{T}^3$. Note that condition (\ref{eqq12}) is a symmetry condition, which is preserved by system (\ref{eqq14})-(\ref{eqq17}), that is if a smooth solution to system (\ref{eqq14})-(\ref{eqq17}) exists and is unique, then it must satisfy the symmetry condition (\ref{eqq12}), as long as it is initially satisfied.

Note that the vertical velocity $\theta$ can be expressed in terms of the horizonal velocity $v$, through the incompressibility condition (\ref{eqq16}) and the symmetry condition (\ref{eqq12}), as
\begin{equation}
\theta(t,x,y,z)=\Phi(v)(t,x,y,z)=-\int^{z}_{-1}\nabla_H\cdot v(t,x,y,z')dz',
\end{equation}
moreover,
\[
\int^{1}_{-1}\nabla_H\cdot v  dz=0.
\]
Supposing that $p_{b}$ is a certain unknown function at $\Gamma_{b}:=\mathbb{T}^2\times\{-1\}$, and integrating (\ref{eqq15}) from $-1$ to $z$, we have
\[
P(x,y,z,t)= p_{b}(x,y,t)-\int^{z}_{-1} T(x,y,z',t) dz'.
\]
 Now, (\ref{eqq14})-(\ref{eqq13}) can be rewritten as
 \begin{eqnarray}\label{eq5-1}
\frac{\partial v}{\partial t}+(v\cdot \nabla_H)v+\Phi(v)\frac{\partial v}{\partial z}+f{k}\times v +\nabla_H p_{b}-\int^{z}_{-1}\nabla_H T dz' -\Delta v &=&\phi(v,T),\\
\label{eq-6-1}
\frac{\partial T}{\partial t}+(v\cdot\nabla_H)T+\Phi(v)\frac{\partial T}{\partial z}+ \Phi(v) -\Delta T&=&\varphi(v,T),\\
\label{eq-7-1}
\int^{1}_{-1}\nabla_H\cdot v  dz&=&0.
\end{eqnarray}
The boundary and initial conditions for (\ref{eq5-1})-(\ref{eq-7-1}) are given by
\begin{eqnarray}\label{eqq18}
v\ and \ T\ are\ periodic\ in\ x, \ y \ and \ z,\\
\label{eqq19}
v\ and\ P\ are\ even\ in \ z, \  \theta\ and\ T\ are\ odd\ in\ z,\\
\label{eqq20}
(v, T)\mid_{t=0}=(v_0,T_0).
\end{eqnarray}
It is easy to know that  Markov Selection and $\mathcal{W}$-strong Feller property for $(v, T)$ of (\ref{eq5-1})-(\ref{eqq20}) implies the same results of the original solution $(v,T)$ of the system (\ref{eqq14})-(\ref{eqq13}). In the following, we will focus on  (\ref{eq5-1})-(\ref{eqq20}).
\section{Formulation of (\ref{eq5-1})-(\ref{eqq20})}
\subsection{ Functional Spaces }
  Let $\mathcal{L}(K_1;K_2)$ (resp. $\mathcal{L}_2(K_1;K_2)$) be the space of bounded (resp. Hilbert-Schmidt) linear operators from the Hilbert space $K_1$ to $K_2$, the norm is denoted by $\|\cdot\|_{\mathcal{L}(K_1;K_2)}(\|\cdot\|_{\mathcal{L}_2(K_1;K_2)})$. Denote by $|\cdot|_{L^p(\mathbb{T}^2)}$ the norm of $L^p(\mathbb{T}^2)$ and $|\cdot|_{p}$ the norm of $L^p(\mathbb{T}^3)$ for $p\in \mathbb{N}_{+}$. In particular, $|\cdot|$ and $(\cdot,\cdot)$ represent the norm and inner product of $L^2(\mathbb{T}^3)$. For the classical Sobolev space $W^{m,2}(\mathbb{T}^3)$, $m\in \mathbb{N}_+$,
\begin{equation}\notag
\left\{
  \begin{array}{ll}
    W^{m,2}(\mathbb{T}^3)=\Big\{U\in L^2(\mathbb{T}^3)\Big| \partial_{\alpha}U\in L^2(\mathbb{T}^3)\ {\rm for} \ |\alpha|\leq m\Big\},&  \\
    |U|^2_{W^{m,2}(\mathbb{T}^3)}=\sum_{0\leq|\alpha|\leq m}|\partial_{\alpha}U|^2. &
  \end{array}
\right.
\end{equation}
It's known that $(W^{m,2}(\mathbb{T}^3), |\cdot|_{W^{m,2}(\mathbb{T}^3)})$ is a Hilbert space.
%\begin{equation}\notag
%\left\{
%  \begin{array}{ll}
%    H^p(\mathbb{T}^3)=\Big\{U\in L^2(\mathbb{T}^3)\Big| \partial_{\alpha}U\in L^2(\mathbb{T}^3)\ {\rm for} \ |\alpha|\leq p\Big\},&  \\
%    |U|^2_{H^p(\mathbb{T}^3)}=\sum_{0\leq|\alpha|\leq p}|\partial_{\alpha}U|^2. &
%  \end{array}
%\right.
%\end{equation}
%It's  known that $(H^p(\mathbb{T}^3), |\cdot|_{H^p(\mathbb{T}^3)})$ is a Hilbert space.
  %By the Riesz representation theorem, we can identify the dual space $H'$ of $H$ (respectively $H'_{i}= H_i$, i=1,2). Then we have
%\[
%V\subset H = H'\subset V',
%\]
%where the two inclusions are compact continuous.

Define working spaces for equations (\ref{eq5-1})-(\ref{eqq20}). Let
 \begin{eqnarray}\notag
 &&\mathcal{V}_1:=\left\{v\in (C^{\infty}(\mathbb{T}^3))^2;\ \int^1_{-1}\nabla_H\cdot v dz=0, v \ is\ periodic\ in\ x, \ y \ and \ even\ in\ z, \int_{\mathbb{T}^3}vdxdydz=0\right\},\\ \notag
 &&\mathcal{V}_2:=\left\{T\in C^{\infty}(\mathbb{T}^3);\  \ T\ is\ periodic\ in\ x, \ y \ and \ odd\ in\ z, \int_{\mathbb{T}^3}Tdxdydz=0 \right\},
 \end{eqnarray}
$V_1$= the closure of $ \mathcal{V}_1$ with respect to the norm $|\cdot|_{W^{1,2}(\mathbb{T}^3)}\times |\cdot|_{W^{1,2}(\mathbb{T}^3)}$,\\
$V_2$= the closure of $ \mathcal{V}_2$ with respect to the norm $|\cdot|_{W^{1,2}(\mathbb{T}^3)}$,\\
$H_1$= the closure of $ \mathcal{V}_1$ with respect to the norm $|\cdot|\times |\cdot|$,\\
$H_2$= the closure of $ \mathcal{V}_2$ with respect to the norm $|\cdot|$,\\
\[
V=V_1\times V_2, \quad H=H_1\times H_2.
\]
The inner products and norms on $V$, $H$ are given by
\begin{eqnarray*}
(U,U_1)_{V}&=&(v,v_1)_{V_1}+(T,T_1)_{V_2},\\
(U,U_1)&=&(v,v_1)+(T,T_1)=(v^{(1)},v^{(1)}_1)+(v^{(2)},v^{(2)}_1)+(T,T_1),\\
(U,U)^{\frac{1}{2}}_{V}&=&(v,v)^{\frac{1}{2}}_{V_1}+(T,T)^{\frac{1}{2}}_{V_2}, \quad \|U\|_{V}=(U,U)^{\frac{1}{2}}_V.
\end{eqnarray*}
where $U=(v,T), U_1=(v_1,T_1), v=(v^{(1)}, v^{(2)})$ and $v_1=(v^{(1)}_1, v^{(2)}_1)$ .

On the periodic domain $\mathbb{T}^3$, it's known that $-\Delta$ is a self-adjoint compact operator, denote by $\{e_n\}_{n=1,2,\cdot\cdot\cdot}$ an eigenbasis and $\{\lambda_n\}_{n=1,2,\cdot\cdot\cdot}$ the corresponding increasing eigenvalue sequence of $-\Delta$. For $s\in \mathbb{R}^+$, define
\[
\|f\|^2_s=\sum^{\infty}_{k=1}|\lambda_k|^{s}|(f,e_k)|^2
\]
and let $H^s(\mathbb{T}^3)$ denote the Sobolev space of all $f\in H$ for which $\|f\|_s$ is finite. It is easy to know that $\|f\|_0=|f|$ and $\|f\|_1=|f|_{W^{1,2}(\mathbb{T}^3)}$. For simplicity, denote $\|\cdot\|_1=\|\cdot\|$.
For $s<0$, define $H^s(\mathbb{T}^3)$ to be the dual of $H^{-s}(\mathbb{T}^3)$. Set $\Lambda=(-\Delta)^{\frac{1}{2}}$, then
\[
\|f\|^2_s=|\Lambda^s f|^2, \quad |\Lambda^s v|^2=|\Lambda^s v^{(1)}|^2+|\Lambda^s v^{(2)}|^2, \quad \|U\|^2_s=|\Lambda^s v|^2+|\Lambda^s T|^2.
\]
%For $s\geq0$, $p\in[1,+\infty]$, we use $H^{s,p}$ to denote a subspace of $L^p(\mathbb{T}^3)$, consisting of all $f$ which can be written in the form $f=\Lambda^{-s}g$, for $g\in L^p(\mathbb{T}^3)$ and the $H^{s,p}$ norm of all $f$ is defined to be the $L^p$ norm of $g$, that is $\|f\|_{H^{s,p}}=|\Lambda^s f|_p$.
\subsection{Functionals}
Define three bilinear forms $a:V\times V\rightarrow \mathbb{R}$,\ $a_1:V_1\times V_1\rightarrow \mathbb{R}$,\ $a_2:V_2\times V_2\rightarrow \mathbb{R}$,
 and their corresponding linear operators $A: V\rightarrow V^{'}$, $A_1: V_1\rightarrow V^{'}_1$, $A_2: V_2\rightarrow V^{'}_2$ by setting
 \[
 a(U,U_1):=(AU,U_1)=   %该矩阵一共3 列，每一列都居中放置
   a_1(v,v_1)+  %第一行元素
  a_2(T,T_1)\\  %第二行元素
  ,
 \]
 where
\begin{eqnarray}\notag
a_1(v,v_1):=(A_1v, v_1)=\int_{\mathbb{T}^3}\left(\nabla_H v\cdot \nabla_H v_1+\frac{\partial v}{\partial z}\cdot\frac{\partial v_1}{\partial z}\right)dxdydz,
\end{eqnarray}
\begin{eqnarray}\notag
a_2(T,T_1):=(A_2T, T_1)=\int_{\mathbb{T}^3}\left(\nabla_H T\cdot \nabla_H T_1+\frac{\partial T}{\partial z}\frac{\partial T_1}{\partial z}\right)dxdydz,
\end{eqnarray}
for any $U=(v,T)$, $U_1=(v_1, T_1)\in V$.

%Define three bilinear forms $a:V\times V\rightarrow \mathbb{R}$,\ $a_1:V_1\times V_1\rightarrow \mathbb{R}$,\ $a_2:V_2\times V_2\rightarrow \mathbb{R}$,
% and their corresponding linear operators $A: V\rightarrow V^{'}$, $A_1: V_1\rightarrow V^{'}_1$, $A_2: V_2\rightarrow V^{'}_2$ by setting
% \[
% a(U,U_1):=(AY,U_1)=   %该矩阵一共3 列，每一列都居中放置
%   a_1(v,v_1)+  %第一行元素
%  a_2(T,T_1)\\  %第二行元素
%  ,
% \]
% where
%\begin{eqnarray}\notag
%a_1(v,v_1):=(A_1v, v_1)=\int_{\mathbb{T}^3}\left(\nabla_H v\cdot \nabla_H v_1+\frac{\partial v}{\partial z}\cdot\frac{\partial v_1}{\partial z}\right)dxdydz,
%\end{eqnarray}
%\begin{eqnarray}\notag
%a_2(T,T_1):=(A_2T, T_1)=\int_{\mathbb{T}^3}\left(\nabla_H T\cdot \nabla_H T_1+\frac{\partial T}{\partial z}\frac{\partial T_1}{\partial z}\right)dxdydz+\alpha \int_{\Gamma_u}SS_1dxdy,
%\end{eqnarray}
%for any $U=(v,T)$, $U_1=(v_1, T_1)\in V$. The following lemma follows readily.
\begin{lemma}
\begin{description}
  \item[(i)] The forms $a$, $a_i\ (i=1,2)$ are coercive, continuous, and therefore, the operators $A:V\rightarrow V'$ and $A_i: V_i\rightarrow V'_i\ (i=1,2)$ are isomorphisms. Moreover,
\begin{eqnarray*}
a(U,U_1)&\leq& C_1\|U\|_V\|U_1\|_V,\\
a(U,U)&\geq& C_2\|U\|^2_V,
\end{eqnarray*}
where $C_1$ and $C_2$ are two absolute constants (independent of the physically relevant constants $Re_i$, $Rt_i$, etc).
  \item[(ii)] The isomorphism  $A:V\rightarrow V'$ (respectively $A_i: V_i\rightarrow V'_i\ (i=1,2)$) can be extended to a self-adjoint unbounded linear  operator on $H$ (respectively on $H_i$, i=1,2), with compact inverse $A^{-1}: H\rightarrow H$ (respectively $A^{-1}_i: H_i\rightarrow H_i \ (i=1,2)$).
\end{description}
\end{lemma}

Now, we define three functionals $b: V\times V\times V\rightarrow \mathbb{R}$, $b_i: V_1\times V_i\times V_i\rightarrow \mathbb{R}\ (i=1,2)$  and the associated operators $B: V\times V\rightarrow V'$, $B_i: V_1\times V_i\rightarrow V'_i\ (i=1,2)$ by setting
\begin{eqnarray*}
b(U,U_1,U_2)&:=&(B(U,U_1),U_2)=   %该矩阵一共3 列，每一列都居中放置
   b_1(v,v_1,v_2)+  %第一行元素
   b_2(v, T_1, T_2) %第二行元素
 ,\\
b_1(v,v_1,v_2)&:=&(B_1(v,v_1), v_2)=\int_{\mathbb{T}^3}\left[(v\cdot \nabla_H)v_1+\Phi(v)\frac{\partial v_1}{\partial z}\right]\cdot v_2dxdydz,\\
b_2(v, T_1, T_2)&:=&(B_2(v,T_1), T_2)=\int_{\mathbb{T}^3}\left[(v\cdot \nabla_H)T_1+\Phi(v)\frac{\partial T_1}{\partial z}\right] T_2dxdydz,
\end{eqnarray*}
for any $U=(v, T)$, $U_i=(v_i,T_i)\in V$.

Moreover, we define another functional $g: V\times V \rightarrow \mathbb{R}$ and the associated linear operator $G: V\rightarrow V'$ by
\begin{eqnarray}\notag
g(U,U_1)&:=&(G(U), U_1)\notag
\\      &=&\int_{\mathbb{T}^3}\left[f(k\times v)\cdot v_1+(\nabla_H p_b-\int^z_{-1}\nabla_H Tdz')\cdot v_1+\Phi(v)\cdot T_1\right]dxdydz % 第一行元素
. \notag
\end{eqnarray}
Finally, using the functionals defined above
to obtain the following stochastic evolution equation
\begin{eqnarray}\label{aa}
\left\{
  \begin{array}{ll}
    dU(t)+AU(t)dt+B(U(t),U(t))dt+G(U(t))dt=\Psi(U(t))dW(t), \\
    U(0)=y,
  \end{array}
\right.
\end{eqnarray}
where
\begin{equation}\notag
W=\left(                 %左括号
  \begin{array}{c}   %该矩阵一共3列，每一列都居中放置
    W_1\\  %第一行元素
    W_2 \\  %第二行元素
  \end{array}
\right) ,\quad
\Psi(U)
=\left(                 %左括号
  \begin{array}{cc}   %该矩阵一共3列，每一列都居中放置
   \phi(v,T) & 0 \\  %第一行元素
   0 & \varphi(v,T) \\  %第二行元素
  \end{array}
\right).
\end{equation}
\subsection{ Inequalities}
%We recall some interpolation inequalities used later (see \cite{Adams}).\\
%For $h\in H^1(\mathbb{T}^2)$,
%\begin{eqnarray*}
%|h|_{L^4(\mathbb{T}^2)}&\leq& c|h|^{\frac{1}{2}}_{L^2(\mathbb{T}^2)}|h|^{\frac{1}{2}}_{H^1(\mathbb{T}^2)},\\
%|h|_{L^5(\mathbb{T}^2)}&\leq& c|h|^{\frac{3}{5}}_{L^3(\mathbb{T}^2)}|h|^{\frac{2}{5}}_{H^1(\mathbb{T}^2)},\\
%|h|_{L^6(\mathbb{T}^2)}&\leq& c|h|^{\frac{2}{3}}_{L^4(\mathbb{T}^2)}|h|^{\frac{1}{3}}_{H^1(\mathbb{T}^2)}.
%\end{eqnarray*}
%For $h\in H^1(\mathbb{T}^3)$,
%\begin{eqnarray*}
%|h|_3&\leq& c|h|^{\frac{1}{2}}|h|^{\frac{1}{2}}_{H^1(\mathbb{T}^3)},\\
%%|h|_4&\leq& c|h|^{\frac{1}{4}}|h|^{\frac{3}{4}}_{H^1(\mathbb{T}^3)},\\
%|h|_6&\leq& c|h|_{H^1(\mathbb{T}^3)},\\
%|h|_{\infty}&\leq& c|h|^{\frac{1}{2}}_{H^1(\mathbb{T}^3)}|h|^{\frac{1}{2}}_{H^2(\mathbb{T}^3)}.
%\end{eqnarray*}

Firstly, we recall the integral version of Minkowshi inequality for the $L^p$ spaces, $p\geq 1$. Let $\mathcal{O}_1\subset \mathbb{R}^{m_1}$ and $\mathcal{O}_2\subset \mathbb{R}^{m_2}$ be two measurable sets, where $m_1$ and $m_2$ are two positive integers. Suppose that $f(\xi, \eta)$ is measurable over $\mathcal{O}_1\times \mathcal{O}_2$. Then
\begin{equation}\notag
\left[\int_{\mathcal{O}_1}\left(\int_{\mathcal{O}_2}|f(\xi,\eta)|d\eta\right)^p d\xi\right]^{1/p}\leq \int_{\mathcal{O}_2}\left(\int_{\mathcal{O}_1}|f(\xi, \eta)|^p d\xi\right)^{1/p}d\eta.
\end{equation}

\begin{lemma}(\cite{C-T-1})\label{le-2}
 If $v_1 \in H^1(\mathbb{T}^3),v_2 \in H^2(\mathbb{T}^3),v_3 \in H^1(\mathbb{T}^3)$, then
\begin{description}
 \item[(i)] $|\int_{\mathbb{T}^3}(v_1 \cdot \nabla_H)v_2\cdot v_3 dxdydz|\leq c|\nabla_H v_2||v_3|_3|v_1|_6\leq c |\nabla_H v_2||v_3|^{\frac{1}{2}}|\nabla_H v_3|^{\frac{1}{2}}|\nabla_H v_1|$,
 \item[(ii)] $|\int_{\mathbb{T}^3}\Phi(v_1)v_{2z}\cdot v_3 dxdydz|\leq c|\nabla_H v_1||v_3|^{\frac{1}{2}}|\nabla_H v_3|^{\frac{1}{2}}|\partial_{z} v_2|^{\frac{1}{2}}|\nabla_H \partial_{z} v_2|^{\frac{1}{2}}$.
\end{description}
\end{lemma}
\section{Markov Selection}
In the following, we will introduce Markov selection for stochastic evolution equations using the same notations as \cite{B-M-X}.
\subsection{Preliminaries}
Let $(\mathbb{X}, \rho_{\mathbb{X}})$ be a polish space and set $\Omega := C([0,\infty); \mathbb{X})$. Denote by $\mathcal{B}$ the Borel $\sigma$-field of $\Omega$ and  by $Pr(\Omega)$ the set of all probability measures on $(\Omega, \mathcal{B})$. Define the canonical process $\xi: \Omega\rightarrow \mathbb{X}$ as
\[
\xi_t(\omega)=\omega(t).
\]

For fixed $t\geq 0$, let $\Omega^t := C([t,\infty); \mathbb{X})$ be the space of all continuous functions from $[t,\infty)$ to $\mathbb{X}$ with the metric
\[
\rho^t(x,y):=\sum^{\infty}_{m=\lfloor t\rfloor+1}\frac{1}{2^m}\left(\sup_{s\in[t,m]}\rho_{\mathbb{X}}\Big(x(s),y(s)\Big)\wedge 1\right)
\]
where $\lfloor t\rfloor$ denotes the integer part of $t$. Then $(\Omega^t, \rho^t)$ is a Polish space. For $s\geq t$, define the $\sigma$-algebra $\mathcal{B}^t_s$ on $\Omega^t$ by $\mathcal{B}^t_s:=\sigma[\xi_r: t\leq r\leq s]$, and write $\mathcal{B}^t:=\bigcup_{s\geq t}\mathcal{B}^t_s$. Thus, we have a measurable space with filtration $(\Omega^t,\mathcal{B}^t, (\mathcal{B}^t_s)_{s\geq t})$. If $t=0$, we simply write $(\Omega, \mathcal{B},(\mathcal{B}_s)_{s\geq 0})$.
Finally,
define the map $\Phi_t: \Omega \rightarrow \Omega^t$ defined by
\[
\Phi_t(\omega)(s):= \omega(s-t), \quad s\geq t,
\]
which establishes a measurable isomorphism between $(\Omega,\mathcal{B}, (\mathcal{B}_s)_{s\geq 0})$ and $(\Omega^t,\mathcal{B}^t, (\mathcal{B}^t_s)_{s\geq t})$.

Given $P\in Pr(\Omega)$ and $t>0$, denote $\omega \mapsto P|^{\omega}_{\mathcal{B}_t}: \Omega \rightarrow Pr(\Omega^t)$ a regular conditional probability distribution of $P$ on $\mathcal{B}_t$. Since $\Omega$ is a Polish space and every $\sigma$-field $\mathcal{B}_t$ is finitely generated, such a function exists and is unique, up to $P$-null sets. In particular,
\[
P|^{\omega}_{\mathcal{B}_t}[\xi_t=\omega(t)]=1
\]
for all $\omega\in \Omega$, and if $A\in \mathcal{B}_t$ and $B\in \mathcal{B}^t$,
\[
P(A\cap B)=\int_A P|^{\omega}_{\mathcal{B}_t}(B)P(d\omega).
\]
Refer to \cite{F-M}, we introduce the following definitions.
\begin{dfn}
Given a family $(P_x)_{x\in H}$ of probability measures in $Pr(\Omega)$, the Markov property  can be stated as
\[
P_x|^{\omega}_{\mathcal{B}_t}=\Phi(t)P_{\omega(t)},\quad for\ P_x-a.e.\ \omega\in \Omega,
\]
for each $x\in H$ and for all $t\geq 0$.
\end{dfn}
\begin{dfn}
The family $(P_x)_{x\in H}$ has the almost sure Markov property if for each $x\in H$, there is a set $\Gamma\subset (0,\infty)$ with null Lebesgue measure, such that
\[
P_x|^{\omega}_{\mathcal{B}_t}=\Phi(t)P_{\omega(t)},\quad for\ P_x-a.e.\ \omega\in \Omega,
\]
for all $t\notin \Gamma$.
\end{dfn}

\subsection{A General Criterion}
%In the following, we introduce a sufficient condition to obtain Markov family for stochastic evolution equations in \cite{B-M-X}, for convenience, the same notations are used as \cite{B-M-X}.

Let $\mathbb{H}$ be a separable Hilbert space, with inner product $\langle\cdot,\cdot\rangle_{\mathbb{H}}$ and norm $\|\cdot\|_{\mathbb{H}}$. Let $\mathbb{X}, $ $\mathbb{U}$ be two separable and reflexive Banach spaces with norms $\|\cdot\|_{\mathbb{X}}$ and $\|\cdot\|_{\mathbb{U}}$, such that
\[
\mathbb{U}\subset \mathbb{H}\subset \mathbb{X}
\]
continuously and densely. If we identify the dual of $\mathbb{H}$ with itself, then we get
\[
\mathbb{X}^*\subset \mathbb{H}^*\backsimeq \mathbb{H}\subset \mathbb{X}.
\]
The dual pair between $\mathbb{X}$ and $\mathbb{X}^*$ is denoted by
\[
{}_{\mathbb{X}}\langle x,y\rangle_{\mathbb{X}^*}, \ x\in \mathbb{X},\ y\in  \mathbb{X}^*.
\]
We remark that if $x\in \mathbb{H}$, then
\[
{}_{\mathbb{X}}\langle x,y\rangle_{\mathbb{X}^*}= \langle x,y\rangle_{\mathbb{H}}.
\]
Let $\mathcal{E}$ be a fixed countable dense subset of $\mathbb{X}^*$ which will be chosen in each case and
$(W(t))_{t\geq 0}$ be a cylindrical Brownian motion in another separable Hilbert space $(\mathbb{Y},\|\cdot\|_{\mathbb{Y}})$ with identity covariance. Consider the following evolution equation:
\begin{eqnarray}\label{eq-12}
dX(t)=\mathcal{A}(X(t))dt+\mathcal{R}(X(t))dW(t),\ t\geq 0, \ X(0)=x_0 \in \mathbb{H},
\end{eqnarray}
where $\mathcal{A}: \mathbb{U}\rightarrow \mathbb{X}$ is $\mathcal{B}(\mathbb{U})/ \mathcal{B}(\mathbb{X})$-measurable and $\mathcal{R}: \mathbb{U}\rightarrow \mathcal{L}_2(\mathbb{Y}; \mathbb{H})$ is
$\mathcal{B}(\mathbb{U})/ \mathcal{B}(\mathcal{L}_2(\mathbb{Y}; \mathbb{H}))$-measurable.
\begin{dfn}\cite{B-M-X}\label{dfn-1}
Let $x_0\in \mathbb{H}$. A probability measure $P\in Pr(\Omega)$ is called a martingale solution of (\ref{eq-12}) with initial value $x_0$, if it satisfies
\begin{description}
  \item[(M1)] $P(X(0)=x_0)=1$ and for any $n\in \mathbb{N}_{+}$
  \[
  P\Big\{X\in \Omega: \int^n_0\|\mathcal{A}(X(s))\|_{\mathbb{X}}ds+\int^n_0 \|\mathcal{R}(X(s))\|^2_{\mathcal{L}_2(\mathbb{Y}; \mathbb{H})}ds< +\infty\Big\}=1;
  \]
  \item[(M2)] for every $l\in \mathcal{E}$, the process
  \[
  M_l(t,X):=\langle X(t), l\rangle_{\mathbb{X}^*}-\int^t_0\langle \mathcal{A}(X(s)), l\rangle_{\mathbb{X}^*}ds
  \]
  is a continuous square integrable $\mathcal{B}_t$-martingale with respect to  $P$, whose quadratic variation process is given by
  \[
   \langle M_{l} \rangle (t,X):= \int^t_0\|\mathcal{R}^*(X(s))(l)\|^2_{\mathbb{Y}}ds,
  \]
  where the asterisk denote the adjoint operator of $\mathcal{R}(X(s))$;
  \item[(M3)] for any $p\in \mathbb{N}$, there exists a continuous positive real function $t\mapsto C_{t,p}$ (only depending on $p$ and $\mathcal{A}, \mathcal{R}$), a lower semi-continuous functional $\mathcal{N}_p: \mathbb{U}\rightarrow [0,\infty]$, and a Lebesgue null set $\mathbb{T}_p\subset (0,\infty)$ such that for all $0\leq s\notin \mathbb{T}_p$ and all $t\geq s$
      \begin{eqnarray}
      \mathbb{E}^P\Big(  \sup_{r\in [s,t]}\|X(r)\|^{2p}_{\mathbb{H}}+\int^t_s\mathcal{N}_p(X(r))dr \Big| \mathcal{B}_s\Big)\leq  C_{t-s,p}\cdot (\|X(s)\|^{2p}_{\mathbb{H}}+1).
      \end{eqnarray}
\end{description}
\end{dfn}

\begin{remark}
The above definition of martingale solution is in the sense of Stroock and Varadhan's martingale problem in \cite{S-V}, which is weaker than that in \cite{RR}.
\end{remark}

In \cite{B-M-X}, the authors give a sufficient conditions on $\mathcal{A}$ and $\mathcal{R}$ to obtain Markov family $\{P_{x_0}\}_{x_0\in \mathbb{H}}$ for (\ref{eq-12}). For this purpose, they firstly introduced the following function class $\mathcal{U}^q$, $q\geq 1$ : A lower semi-continuous function $\mathcal{N}: \mathbb{U} \rightarrow [0,\infty]$ belong to $\mathcal{U}^q$ if $\mathcal{N}(y)=0$ implies $y=0$, and
\[
\mathcal{N}(cy)\leq c^q \mathcal{N}(y), \ \forall c\geq 0,\ y\in \mathbb{U}.
\]
and
\[
\Big\{y\in \mathbb{U}: \mathcal{N}(y)\leq 1\Big\} \ is\ relatively\ compact\ in\  \mathbb{U}.
\]
The assumptions on $\mathcal{A}$ and $\mathcal{R}$ are given as follows:
\begin{description}
  \item[(C1)] (Demi-Continuity) For any $x\in \mathbb{X}^*$, if $y_n$ strong converges to $y$ in $\mathbb{U}$, then
      \[
      \lim_{n\rightarrow \infty}{}_{\mathbb{X}}\langle\mathcal{A}(y_n),x\rangle_{\mathbb{X}^*}={}_{\mathbb{X}}\langle\mathcal{A}(y),x\rangle_{\mathbb{X}^*},
      \]
      and
      \[
      \lim_{n\rightarrow \infty}\|\mathcal{R}^*(y_n)(x)-\mathcal{R}^*(y)(x)\|_{\mathbb{Y}}=0.
      \]
  \item[(C2)] (Coercivity Condition) There exist $\lambda_1\geq 0$ and $\mathcal{N}_1\in \mathcal{U}^q$ for some $q\geq 2$ such that for all $x\in \mathbb{X}^*$
      \[
      {}_{\mathbb{X}}\langle\mathcal{A}(x),x\rangle_{\mathbb{X}^*}\leq -\mathcal{N}_1(x)+\lambda_1(1+\|x\|^2_{\mathbb{H}}).
      \]

  \item[(C3)] (Growth Condition) There exist $\lambda_2, \lambda_3, \lambda_4>0$ and $\gamma'\geq \gamma>1$ such that for all $x\in \mathbb{U}$
      \begin{eqnarray*}
      \|\mathcal{A}(x)\|^{\gamma}_{\mathbb{X}}&\leq& \lambda_2\mathcal{N}_1(x)+\lambda_3(1+\|x\|^{\gamma'}_{\mathbb{H}}),\\
      \|\mathcal{R}(x)\|^2_{\mathcal{L}_2(\mathbb{Y};\mathbb{H})}&\leq& \lambda_4(1+\|x\|^2_{\mathbb{H}}),
      \end{eqnarray*}
      where $\mathcal{N}_1$ is as in $(C2)$.
\end{description}
\begin{theorem}(\cite{B-M-X})\label{thm-6}
Assume $(C1)-(C3)$ hold, for each $x_0\in\mathbb{H}$, there exists a martingale solution $P\in Pr(\Omega)$ starting from $x_0$ to (\ref{eq-12}) in the sense of Definition \ref{dfn-1}.
\end{theorem}
Then, the main result is
\begin{theorem}\label{thm-1}
( \cite{B-M-X}) Under $(C1)-(C3)$, there exists an almost sure Markov family $(P_x)_{x\in\mathbb{H}}$ to (\ref{eq-12}).
\end{theorem}

\subsection{Proof of Theorem \ref{thm-2}}
In this part, we will use Theorem \ref{thm-1} to get an almost surely Markov family  $\{P_{x}\}_{x\in H}$ for (\ref{aa}).  Firstly, define the operator $\mathcal{A}$ and $\mathcal{R}$ as follows:
\begin{eqnarray*}
\mathcal{A}(y)&:=&  Ay+B(y,y)+G(y),\\
\mathcal{R}(y)&:=&  \Psi(y)\ \rm{for} \ y\in C^{\infty}(\mathbb{T}^3).
\end{eqnarray*}
Here, for our equation (\ref{aa}), we choose
\[
\mathbb{U}=H^1(\mathbb{T}^3),\ \mathbb{Y}=\mathbb{H}=H,\ \mathbb{X}=(H^3(\mathbb{T}^3))^*, \ \mathbb{X}^*=H^3(\mathbb{T}^3),
\]
then $\mathbb{X}$ is a Hilbert space and $\mathbb{X}^* \subset \mathbb{U}$ compactly. Moreover, the covariance operator $\Psi$ is assumed to satisfy
 %We can deduce \textbf{Hypothesis H1} from \textbf{Hypothesis H0} that
\begin{description}
   \item[\textbf{Hypothesis H0}] (i) $\Psi:$ $H^1(\mathbb{T}^3)\rightarrow \mathcal{L}_2(H;H)$ is a continuous and bounded Lipschitz mapping, i.e.
       \[
       \|\Psi(y)\|^2_{\mathcal{L}_2(H;H)}\leq \lambda_0|y|^2 +\rho   \quad y\in H^1(\mathbb{T}^3),
       \]
for some constants $\lambda_0\geq0,\ \rho\geq0.$

(ii) If $y,y_n\in H^1(\mathbb{T}^3)$, such that $y_n$ strongly converges to $y$ in $H^1(\mathbb{T}^3)$, then for any $x\in C^{\infty}(\mathbb{T}^3)$,
\[
|\Psi(y_n)^*(x)-\Psi(y)^*(x)|_H\rightarrow 0 \quad  n\rightarrow \infty.
\]
\end{description}

%Let $\mathcal{E}=\{e_i, i\in \mathbb{N}\}$ be the orthonormal basis of $H$ introduced in Sect. 2.
By Lemma \ref{lem-1} below, $\mathcal{A}$ can be extended to an operator $\mathcal{A}: H^1(\mathbb{T}^3) \rightarrow \mathbb{X}$. For $y \notin H^1(\mathbb{T}^3)$, $\mathcal{A}(y):= \infty$.
\begin{lemma}\label{lem-1}
For any $y_1$, $y_2\in C^{\infty}(\mathbb{T}^3)$,
\begin{eqnarray*}
\|Ay_1-Ay_2\|_{\mathbb{X}}&\leq& C_1|y_1-y_2|,\\
\|B(y_1,y_1)-B(y_2,y_2)\|_{\mathbb{X}}&\leq&  C_2(\|y_1\|+\|y_2\|)\|y_1-y_2\|,\\
\|G(y_1)-G(y_2)\|_{\mathbb{X}}&\leq& C_3|y_1-y_2|.
\end{eqnarray*}
for constants $C_1, C_2, C_3$. In particular, the operator $\mathcal{A}: C^{\infty}(\mathbb{T}^3)\rightarrow \mathbb{X}$  extends to an operator $\mathcal{A}: H^1(\mathbb{T}^3)\rightarrow \mathbb{X}$ by continuity.
\end{lemma}
\textbf{Proof of Lemma \ref{lem-1}} \quad We only prove the second assertion, the first and third estimates can be obtained by H\"{o}lder inequality. Refer to \cite{L-T-W-2},
\[
\|B(y,y_1)\|_{-3}\leq C|y|\|y_1\|.
\]
Then,
\begin{eqnarray*}
\|B(y_1,y_1)-B(y_2,y_2)\|_{\mathbb{X}}
&=&\sup_{x\in C^{\infty}(\mathbb{T}^3): \|x\|_{3}\leq 1 }|\langle B(y_1,y_1)-B(y_2,y_2),x\rangle|\\
&=& \sup_{x\in C^{\infty}(\mathbb{T}^3): \|x\|_{3}\leq 1 }|\langle B(y_1,y_1-y_2)-B(y_1-y_2,y_2),x\rangle|\\
&\leq & C_2(\|y_1\|+\|y_2\|)\|y_1-y_2\|.
\end{eqnarray*}
$\hfill\blacksquare$

In order to use Theorem \ref{thm-1}, define the functional $\mathcal{N}_1$ on $\mathbb{U}$ as follows:
\begin{eqnarray*}
\mathcal{N}_1(y):= \left\{
                             \begin{array}{ll}
                               \|y\|^2, & if\ y\in H^1(\mathbb{T}^3), \\
                               +\infty, & otherwise.
                             \end{array}
                           \right.
\end{eqnarray*}
It is obvious that $\mathcal{N}_1\in \mathcal{U}^2$.

%\begin{theorem}\label{thm-3}
%Assume \textbf{Hypothesis H0} holds, then for each $x\in H$, there exists a martingale solution $P_x\in Pr(\Omega)$ to equation (\ref{eq-12}) in the sense of Definition \ref{dfn-1}.
%\end{theorem}
\textbf{Proof of Theorem \ref{thm-2}} \quad By Theorem \ref{thm-1}, we only need to check (C1)-(C3) for $\mathcal{A}$ and $\mathcal{R}$.
\begin{description}
  \item[(I)] The demi-continuity condition (C1) holds by Lemma \ref{lem-1} and \textbf{Hypothesis H0}.
  \item[(II)] The coercivity condition (C2) follows since
  \[
  \langle B(y,y),y\rangle=0,
  \]
  then, by Young inequality, we have
  \begin{eqnarray*}
  {}_{H^{-3}(\mathbb{T}^3)}\langle \mathcal{A}(y),y\rangle_{H^3(\mathbb{T}^3)}
  &=&{}_{H^{-3}(\mathbb{T}^3)}\langle Ay+B(y,y)+G(y),y\rangle_{H^3(\mathbb{T}^3)}\\
  &=&-\|y\|^2+C|y|\|y\|\\
  &\leq& -\|y\|^2+ \frac{1}{2}\|y\|^2+C|y|^2\\
  &\leq& -\frac{1}{2}\|y\|^2+\lambda_1(1+|y|^2).
  \end{eqnarray*}
  \item[(III)] The growth condition (C3) is clear since by Lemma \ref{lem-1}, it gives
  \[
  \|\mathcal{A}(y)\|^2_{-3}\leq \lambda_1 \|y\|^2 +\lambda_2 (1+|y|^2),
  \]
  and by \textbf{Hypothesis H0}, we have
  \[
  \|\mathcal{R}(y)\|^2_{\mathcal{L}_2(H;{H})}\leq \lambda_3(1+|y|^2).
  \]
\end{description}
$\hfill\blacksquare$
\begin{remark}
By Theorem \ref{thm-6}, for any $x_0\in H$, there exists a martingale solution $P_{x_0}\in Pr(\Omega)$ to (\ref{aa}) in the sense of Definition \ref{dfn-1}. Refer to \cite{B-M-X} and \cite{ROMITO}, we know that $P_{x_0}$ is obtained by means of maximisation.
\end{remark}
\section{$\mathcal{W}$-strong Feller}
In this section, we apply the abstract result ( Theorem 5.4 in \cite{F-M}) to obtain that every Markov selection in Sect. 4 has $\mathcal{W}$-strong Feller property.
\subsection{Preliminaries}
Firstly, we recall the following important lemma (\cite{Resnick}, Lemma A.4):
\begin{lemma}\label{lem-2}
Suppose that $s>0$ and $p\in(1,\infty)$. If $f,g\in C^{\infty}(\mathbb{T}^3)$, then
\[
|\Lambda^s(fg)|_p\leq C(|f|_{p_1}|\Lambda^sg|_{p_2}+|g|_{p_3}|\Lambda^sf|_{p_4}),
\]
with $p_i\in(1,\infty]$, $i=1,\cdot\cdot\cdot,4$ such that
\[
\frac{1}{p}=\frac{1}{p_1}+\frac{1}{p_2}=\frac{1}{p_3}+\frac{1}{p_4}.
\]
\end{lemma}
\begin{remark}
$|g|_{p_3}|\Lambda^sf|_{p_4}$ can be equal to $|f|_{p_1}|\Lambda^sg|_{p_2}$ by choosing suitable parameters $p_3$ and $p_4$, in that case, we only write one of them.
\end{remark}
We will also use the following Sobolev inequality (\cite{Stein}, Chapter V ):
\begin{lemma}\label{lem-3}
Suppose that $q>1$,$p\in[q,\infty)$ and
\[
\frac{1}{p}+\frac{\sigma}{3}=\frac{1}{q}.
\]
If $\Lambda^{\sigma}f\in L^q(\mathbb{T}^3)$, then $f\in L^p(\mathbb{T}^3)$ and there is a constant $C\geq0$ independent of $f$ such that
\[
|f|_p\leq C|\Lambda^{\sigma}f|_q.
\]
\end{lemma}
%Remark: During the following proof, $q=2$ and $q=3$ are used.
We shall use as well the following interpolation inequality ( \cite{Ju}, (5.5)).
\begin{lemma}\label{lem-4}
For $f\in C^{\infty}(\mathbb{T}^3)$, we have
\[
\|f\|_s\leq C\|f\|^{\frac{s_2-s}{s_2-s_1}}_{s_1}\|f\|^{\frac{s-s_1}{s_2-s_1}}_{s_2},\quad s_1<s<s_2.
\]
\end{lemma}
Refer to the appendix of \cite{C-T-2}, we have the following lemma.
\begin{lemma} For any $v,T $ and $\omega\in C^{\infty}(\mathbb{T}^3)$,
\begin{eqnarray*}
&&\Big|\int_{\mathbb{T}^3}\Phi(v)\frac{\partial T}{\partial z}\omega\Big|\\
&\leq& \int_{\mathbb{T}^2}\Big(\int^1_{-1}|\nabla_H v|dz'\Big)\Big(\int^1_{-1}|\frac{\partial T}{\partial z}\omega|dz\Big)dxdy\\
&\leq &C\int_{\mathbb{T}^2}\Big(\int^1_{-1}|\nabla_H v|dz'\Big)\Big(\int^1_{-1}|\frac{\partial T}{\partial z}|^2dz\Big)^\frac{1}{2}\Big(\int^1_{-1}|\omega|^2dz\Big)^\frac{1}{2}dxdy\\
&\leq& C|\omega|\|\nabla_H v\|_{s_1}\|T\|_{s_2},
\end{eqnarray*}
\end{lemma}
where $s_1+s_2=1$.

At last, we introduce the definition of $\mathcal{W}$-strong Feller.
\begin{dfn}
($\mathcal{W}$-strong Feller) A given semigroup $(P_t)_{t\geq 0}$ on $B_b(H)$ is $\mathcal{W}$-strong Feller, if for any $t>0$ and $\psi\in B_b(H)$, $P_t\psi\in C_b(\mathcal{W})$.
\end{dfn}
\subsection{$\mathcal{W}$ Space and Hypothesis }
For any $\varepsilon_0>0$ and set
\[
\mathcal{W}=D(A^{\frac{3}{4}+\varepsilon_0}),\ |U|_{\mathcal{W}}=|v|_{\mathcal{W}}+|T|_{\mathcal{W}}=|\Lambda^{\frac{3}{2}+2\varepsilon_0} v|+|\Lambda^{\frac{3}{2}+2\varepsilon_0} T|.
\]
In this section, we choose
\[
\Omega := C([0,\infty); H^{-\beta})
\]
for some $\beta >3$ and $\mathcal{B}$ denote the Borel $\sigma$-algebra on $\Omega$.
We assume the noise is additive, nondegenerate and regular. Concretely,
\begin{description}
   \item[\textbf{Hypothesis H1}] There are an isomorphism $Q_0$ of $H$ and a number $\alpha_0 = \frac{1}{2}+\varepsilon_0$ such that
       \[
       \Psi=Q^{\frac{1}{2}}=A^{-\frac{3}{4}-\alpha_0}Q^{\frac{1}{2}}_0=A^{-\frac{5}{4}-\varepsilon_0}Q^{\frac{1}{2}}_0,
       \]
       where the covariant $Q$: $H\rightarrow H$ is a symmetric non-negative trace-class operator on $H$.
 \end{description}
 \begin{remark}
 Firstly, we notice that \textbf{Hypothesis H1} implies \textbf{Hypothesis H0}. Indeed, the operator $A^{-\frac{3}{4}-\varepsilon}$ is Hilbert-Schmidt in $H$, for every $\varepsilon >0$. Moreover, $A^{-\frac{3}{4}-\varepsilon}Q^{\frac{1}{2}}_0W(t)$
 is a Brownian motion in H, for every $\varepsilon>0$ and every isomorphism $Q_0$ of H, where $W(t)$ is a cylindrical Wiener process on H. In conclusion, $A^{-\frac{3}{4}-\alpha_0}Q^{\frac{1}{2}}_0W(t)$ is a Brownian motion in $D(A^{\alpha})$ for every $\alpha_0>\alpha >0$.
\end{remark}
In the following, we will consider equations (\ref{aa}) in the following abstract form:
\begin{eqnarray}\label{bb}
\left\{
  \begin{array}{ll}
    dU(t)+AU(t)dt+B(U(t),U(t))dt+G(U(t))dt=Q^{\frac{1}{2}}dW(t),  \\
    U(0)=y.
  \end{array}
\right.
\end{eqnarray}
\begin{remark}
Under \textbf{Hypothesis H1}, in \cite{Guo}, the authors have proved that for $y\in V$, there exists a unique strong solution $U=(v,T)$. However, for $y\in H$, the uniqueness of the weak solution is still open, hence, we have to deal with the selected Markov process.
\end{remark}
For $y\in H$, let $P_y$ denote the law of the corresponding solution $U(\cdot, y)$ to (\ref{bb}). Since \textbf{Hypothesis H1} implies \textbf{Hypothesis H0}, by Theorem \ref{thm-2}, the measures $P_y$, $y\in H$ form a Markov process. Let $(P_t)_{t\geq 0}$ be the associated transition semigroup on $B_b(H)$, defined as
\begin{eqnarray}\label{eq-30}
P_t(\varphi)(y):= \mathbb{E}[\varphi(U(t,y))] \quad \forall y\in H,\ \forall \varphi\in B_b(H).
\end{eqnarray}

\subsection{Proof of Theorem \ref{thm-4}}
For the proof, we shall use Theorem 5.4 in \cite{F-M}, which is an abstract result to prove the strong Feller property of Markov selection. In order to achieve this, we follow the idea of Theorem 5.11 in \cite{F-M} to construct $P^{(R)}_y$. We introduce an equation which differs from the original one by a cut-off only, so that with large probability they have the same trajectories on a small random time interval. Consider
\begin{eqnarray}\label{eq-16}
dU(t)+AU(t)dt+\chi_R(|U|^2_{\mathcal{W}})\Big[B(U(t),U(t))+G(U(t))\Big]dt=Q^{\frac{1}{2}}dW(t),
\end{eqnarray}
where $\chi_R: \mathbb{R}\rightarrow [0,1]$ is of class $C^{\infty}$ such that $\chi_R(|U|)=1$ if $|U|\leq R$, and $\chi_R(|U|)=0$ if $|U|\geq R+1$ and its first derivative bounded by 1.
\begin{theorem}\label{thm-5}
\textbf{(Weak-strong uniqueness)} Suppose \textbf{Hypothesis H1} holds. Then for every $y\in \mathcal{W}$, equation (\ref{eq-16}) has a unique martingale solution $P^{(R)}_y$, with
\[
P^{(R)}_y[C([0,\infty);\mathcal{W})]=1.
\]
Let $\tau_R:\Omega\rightarrow [0,\infty]$ be defined by
\[
\tau_R(\omega):= \inf\{t\geq 0: |\omega(t)|^2_{\mathcal{W}}\geq R\}
\]
and $\tau_R(\omega):= \infty$ if this set is empty. If $y\in \mathcal{W}$ and $|y|^2_{\mathcal{W}}< R$, then
\begin{eqnarray}\label{eq-39}
\lim_{\varepsilon \rightarrow 0} P^{(R)}_{y+h}[\tau_R\geq \varepsilon]=1, \ uniformly \ in \ h\in \mathcal{W}, \ |h|_{\mathcal{W}}<1.
\end{eqnarray}
Moreover,
\begin{eqnarray}\label{eq-40}
\mathbb{E}^{P^{(R)}_y}[\varphi(\omega_t)I_{[\tau_R \geq t]}]=\mathbb{E}^{P_y}[\varphi(\omega_t)I_{[\tau_R \geq t]}]
\end{eqnarray}
 for every $t\geq 0$ and $\varphi \in B_b(H)$.
\end{theorem}
\textbf{Proof of Theorem \ref{thm-5}} \quad Let $Z$ be the solution to
\[
d Z(t)+A Z(t) dt = Q^{\frac{1}{2}}dW(t),
\]
with the initial data $Z(0)=0$ and let $X^{(R)}_y$ be the solution to the auxilary problem
\begin{eqnarray}\label{eq-17}
\frac{dX^{(R)}(t) }{dt}+AX^{(R)}(t)+\chi_R(|X^{(R)}+Z|^2_{\mathcal{W}})\Big[B(X^{(R)}+Z,X^{(R)}+Z)+G(X^{(R)}+Z)\Big]=0,
\end{eqnarray}
with $X^{(R)}(0)=y$. Moreover, define $U^{(R)}(t)=X^{(R)}(t)+Z(t)$, which is a weak solution to equation (\ref{eq-16}). We denote its law on $\Omega$ by $P^{(R)}_y$.
For the noise, by \textbf{Hypothesis H1}, the trajectories of the noise belong to
\[
\Omega^*:= \bigcap_{\beta\in (0,\frac{1}{2}),\ \alpha\in [0,\frac{1}{2}+\varepsilon_0)}C^{\beta}([0,\infty);D(A^{\alpha}))
\]
with probability one. Hence, the analyticity of the semigroup generated by $A$ implies that for each $\omega\in \Omega^*$,
 \begin{eqnarray}\label{eq-22}
 Z(\omega)\in C([0,\infty);D(A^{1+\varepsilon_0-\varepsilon}))\subseteq C([0,\infty);\mathcal{W}),
 \end{eqnarray}
 for every $\varepsilon\in(0,\frac{1}{4})$.

 Now, fix $\omega\in \Omega^*$, we will prove that equation (\ref{eq-17}) has a unique global weak solution in $C([0,\infty);\mathcal{W})$.

Denoting by $X=(\kappa^{(R)},g^{(R)})$, $Z=(Z_1,Z_2)$, $\kappa^{(R)}=v^{(R)}-Z_1$, $g^{(R)}=T-Z_2$, then (\ref{eq-17}) can be rewritten as
\begin{eqnarray}\label{eq-18}
\frac{\partial \kappa^{(R)}}{\partial t}&+&\Big((\kappa^{(R)}+ Z_1)\cdot \nabla_H\Big)(\kappa^{(R)}+ Z_1)+\Phi(\kappa^{(R)}+Z_1)\frac{\partial( \kappa^{(R)}+ Z_1)}{\partial z} +f{k}\times (\kappa^{(R)}+Z_1)\\ \notag
&+&\nabla_H p_{b}-\int^{z}_{-1}\nabla_H (g^{(R)}+Z_2)dz'-\Delta\kappa^{(R)}=0,
\end{eqnarray}
\begin{eqnarray}\label{eq-19}
\frac{\partial g^{(R)}}{\partial t}+\Big((\kappa^{(R)}+ Z_1)\cdot \nabla_H\Big)(g^{(R)}+ Z_2)+\Phi(\kappa^{(R)}+Z_1)\frac{\partial (g^{(R)}+Z_2)}{\partial z}+\Phi(\kappa^{(R)}+Z_1)-\Delta g^{(R)} =0.
\end{eqnarray}

 \textbf{(Existence of weak solution)} \quad  Multiplying (\ref{eq-18}) by $-\Lambda^{3+4\varepsilon_0} \kappa^{(R)}$, integrating over $\mathbb{T}^3$, it follows that
 \begin{eqnarray*}
 &&\frac{1}{2}\frac{d}{dt}|\Lambda^{\frac{3}{2}+2\varepsilon_0}\kappa^{(R)}|^2+|\Lambda^{\frac{5}{2}+2\varepsilon_0}\kappa^{(R)}|^2\\
 &=&\chi_R(|U^{(R)}|^2_{\mathcal{W}})\int_{\mathbb{T}^3}\Lambda^{3+4\varepsilon_0}\kappa^{(R)} \Big[\Big((\kappa^{(R)}+ Z_1)\cdot \nabla_H\Big)(\kappa^{(R)}+ Z_1)\Big]dxdydz\\ \notag
 && +\chi_R(|U^{(R)}|^2_{\mathcal{W}})\int_{\mathbb{T}^3}\Lambda^{3+4\varepsilon_0}\kappa^{(R)} \Big[\Phi(\kappa^{(R)}+Z_1)\frac{\partial( \kappa^{(R)}+ Z_1)}{\partial z}\Big]dxdydz\\ \notag
 &&+\chi_R(|U^{(R)}|^2_{\mathcal{W}})\int_{\mathbb{T}^3}\Lambda^{3+4\varepsilon_0}\kappa^{(R)} \Big[f{k}\times (\kappa^{(R)}+Z_1)+\nabla_H p_{b}-\int^{z}_{-1}\nabla_H (g^{(R)}+Z_2)dz'\Big] dxdydz\\
 &:=& \chi_R(|U^{(R)}|^2_{\mathcal{W}})(I_1+I_2+I_3).
 \end{eqnarray*}
For $I_1$, we have
\begin{eqnarray*}
&&\Big|\int_{\mathbb{T}^3}\Lambda^{3+4\varepsilon_0}\kappa^{(R)}[(\kappa^{(R)}\cdot \nabla_H)\kappa^{(R)}] dxdydz\Big|\\
&=& \Big|\int_{\mathbb{T}^3}\Lambda^{\frac{5}{2}+2\varepsilon_0}\kappa^{(R)}\Lambda^{\frac{1}{2}+2\varepsilon_0}[(\kappa^{(R)}\cdot \nabla_H)\kappa^{(R)}]  dxdydz \Big|\\
%&\leq & \int_{\mathbb{T}^3}|\Lambda^{2+2\varepsilon_0}\kappa^{(R)}||\Lambda^{2\varepsilon_0}\kappa^{(R)}||\nabla_H \kappa^{(R)}|dxdydz+\int_{\mathbb{T}^3}|\Lambda^{2+2\varepsilon_0}\kappa^{(R)}||\kappa^{(R)}||\Lambda^{2\varepsilon_0}\nabla_H \kappa^{(R)}|dxdydz\\
&\leq & C|\Lambda^{\frac{5}{2}+2\varepsilon_0}\kappa^{(R)}||\Lambda^{\frac{1}{2}+2\varepsilon_0+\sigma_1}\kappa^{(R)}||\Lambda^{1+\sigma_2}\kappa^{(R)}|+
C|\Lambda^{\frac{5}{2}+2\varepsilon_0}\kappa^{(R)}||\Lambda^{s_1}\kappa^{(R)}||\Lambda^{\frac{3}{2}+2\varepsilon_0+s_2}\kappa^{(R)}|\\
&\leq & C|\Lambda^{\frac{5}{2}+2\varepsilon_0}\kappa^{(R)}||\Lambda^{\frac{3}{2}}\kappa^{(R)}||\Lambda^{\frac{3}{2}+2\varepsilon_0}\kappa^{(R)}|\\
&\leq & \varepsilon | \Lambda^{\frac{5}{2}+2\varepsilon_0}\kappa^{(R)}|^2+C|\Lambda^{\frac{3}{2}+2\varepsilon_0}\kappa^{(R)}|^4,
\end{eqnarray*}
where the first equality follows from Lemma \ref{lem-2}. In the first inequality, $\sigma_1+\sigma_2=\frac{3}{2}$, $s_1+s_2=\frac{3}{2}$, we choose
\[
\sigma_1=1-2\varepsilon_0, \quad \sigma_2=\frac{1}{2}+2\varepsilon_0,\quad
s_1=\frac{3}{2}, \quad s_2=0.
\]
The second inequality follows from Lemma \ref{lem-4}. The Young inequality is used in the last inequality. By the same argument, we have
\begin{eqnarray*}
&&\Big|\int_{\mathbb{T}^3}\Lambda^{3+4\varepsilon_0}\kappa^{(R)}[(\kappa^{(R)}\cdot \nabla_H) Z_1)] dxdydz\Big|\\
&=& \Big|\int_{\mathbb{T}^3}\Lambda^{\frac{5}{2}+2\varepsilon_0}\kappa^{(R)}\Lambda^{\frac{1}{2}+2\varepsilon_0}[(\kappa^{(R)}\cdot \nabla_H) Z_1] dxdydz \Big|\\
%&\leq & \int_{\mathbb{T}^3}|\Lambda^{2+2\varepsilon_0}\kappa^{(R)}||\Lambda^{2\varepsilon_0}\kappa^{(R)}||\nabla_H Z_1|dxdydz+\int_{\mathbb{T}^3}|\Lambda^{2+2\varepsilon_0}\kappa^{(R)}||\kappa^{(R)}||\Lambda^{2\varepsilon_0}\nabla_H Z_1|dxdydz\\
&\leq & C|\Lambda^{\frac{5}{2}+2\varepsilon_0}\kappa^{(R)}||\Lambda^{\frac{1}{2}+2\varepsilon_0+\sigma_1}\kappa^{(R)}||\Lambda^{1+\sigma_2}Z_1|+
C|\Lambda^{\frac{5}{2}+2\varepsilon_0}\kappa^{(R)}||\Lambda^{s_1}\kappa^{(R)}||\Lambda^{\frac{3}{2}+2\varepsilon_0+s_2}Z_1|\\
&\leq& C|\Lambda^{\frac{5}{2}+2\varepsilon_0}\kappa^{(R)}|\Big(|\Lambda^{\frac{3}{2}+2\varepsilon_0}\kappa^{(R)}||\Lambda^{\frac{3}{2}}Z_1|
+|\Lambda^{\frac{3}{2}+2\varepsilon_0}Z_1||\Lambda^{\frac{3}{2}}\kappa^{(R)}|\Big)\\
&\leq &\varepsilon | \Lambda^{\frac{5}{2}+2\varepsilon_0}\kappa^{(R)}|^2+C|\Lambda^{\frac{3}{2}+2\varepsilon_0}\kappa^{(R)}|^2|\Lambda^{\frac{3}{2}}Z_1|^2
+C|\Lambda^{\frac{3}{2}+2\varepsilon_0}Z_1|^2|\Lambda^{\frac{3}{2}}\kappa^{(R)}|^2,
\end{eqnarray*}
where $\sigma_1+\sigma_2=\frac{3}{2}$, $s_1+s_2=\frac{3}{2}$, we choose
\[
\sigma_1=1, \quad \sigma_2=\frac{1}{2},\quad
s_1=\frac{3}{2}, \quad s_2=0.
\]
\begin{eqnarray*}
&&\Big|\int_{\mathbb{T}^3}\Lambda^{3+4\varepsilon_0}\kappa^{(R)}[(Z_1\cdot \nabla_H) \kappa^{(R)})] dxdydz\Big|\\
&=& \Big|\int_{\mathbb{T}^3}\Lambda^{\frac{5}{2}+2\varepsilon_0}\kappa^{(R)}\Lambda^{\frac{1}{2}+2\varepsilon_0}[(Z_1\cdot \nabla_H)\kappa^{(R)}] dxdydz \Big|\\
%&\leq & \int_{\mathbb{T}^3}|\Lambda^{2+2\varepsilon_0}\kappa^{(R)}||\Lambda^{2\varepsilon_0}Z_1||\nabla_H \kappa^{(R)}|dxdydz+\int_{\mathbb{T}^3}|\Lambda^{2+2\varepsilon_0}\kappa^{(R)}||Z_1||\Lambda^{2\varepsilon_0}\nabla_H \kappa^{(R)}|dxdydz\\
&\leq & C|\Lambda^{\frac{5}{2}+2\varepsilon_0}\kappa^{(R)}||\Lambda^{\frac{1}{2}+2\varepsilon_0+\sigma_1}Z_1||\Lambda^{1+\sigma_2}\kappa^{(R)}|+
C|\Lambda^{\frac{5}{2}+2\varepsilon_0}\kappa^{(R)}||\Lambda^{s_1}Z_1||\Lambda^{\frac{3}{2}+2\varepsilon_0+s_2}\kappa^{(R)}|\\
&\leq& C|\Lambda^{\frac{5}{2}+2\varepsilon_0}\kappa^{(R)}||\Lambda^{\frac{3}{2}}Z_1||\Lambda^{\frac{3}{2}+2\varepsilon_0}\kappa^{(R)}|\\
&\leq &\varepsilon | \Lambda^{\frac{5}{2}+2\varepsilon_0}\kappa^{(R)}|^2+C|\Lambda^{\frac{3}{2}}Z_1|^2|\Lambda^{\frac{3}{2}+2\varepsilon_0}\kappa^{(R)}|^2,
\end{eqnarray*}
where $\sigma_1+\sigma_2=\frac{3}{2}$, $s_1+s_2=\frac{3}{2}$, we choose
\[
\sigma_1=1-2\varepsilon_0, \quad \sigma_2=\frac{1}{2}+2\varepsilon_0,\quad
s_1=\frac{3}{2}, \quad s_2=0.
\]
\begin{eqnarray*}
&&\Big|\int_{\mathbb{T}^3}\Lambda^{3+4\varepsilon_0}\kappa^{(R)}[(Z_1\cdot \nabla_H) Z_1)] dxdydz\Big|\\
&=& \Big|\int_{\mathbb{T}^3}\Lambda^{\frac{5}{2}+2\varepsilon_0}\kappa^{(R)}\Lambda^{\frac{1}{2}+2\varepsilon_0}[(Z_1\cdot \nabla_H)Z_1] dxdydz \Big|\\
%&\leq & \int_{\mathbb{T}^3}|\Lambda^{2+2\varepsilon_0}\kappa^{(R)}||\Lambda^{2\varepsilon_0}Z_1||\nabla_H Z_1|dxdydz+\int_{\mathbb{T}^3}|\Lambda^{2+2\varepsilon_0}\kappa^{(R)}||Z_1||\Lambda^{2\varepsilon_0}\nabla_H Z_1|dxdydz\\
&\leq & C|\Lambda^{\frac{5}{2}+2\varepsilon_0}\kappa^{(R)}||\Lambda^{\frac{1}{2}+2\varepsilon_0+\sigma_1}Z_1||\Lambda^{1+\sigma_2}Z_1|+
C|\Lambda^{\frac{5}{2}+2\varepsilon_0}\kappa^{(R)}||\Lambda^{s_1}Z_1||\Lambda^{\frac{3}{2}+2\varepsilon_0+s_2}Z_1|\\
&\leq & C|\Lambda^{\frac{5}{2}+2\varepsilon_0}\kappa^{(R)}||\Lambda^{\frac{3}{2}+2\varepsilon_0}Z_1||\Lambda^{\frac{3}{2}}Z_1|\\
&\leq &\varepsilon | \Lambda^{\frac{5}{2}+2\varepsilon_0}\kappa^{(R)}|^2+C|\Lambda^{\frac{3}{2}}Z_1|^2|\Lambda^{\frac{3}{2}+2\varepsilon_0}Z_1|^2,
\end{eqnarray*}
where $\sigma_1+\sigma_2=\frac{3}{2}$, $s_1+s_2=\frac{3}{2}$, we choose
\[
\sigma_1=1, \quad \sigma_2=\frac{1}{2},\quad
s_1=\frac{3}{2}, \quad s_2=0.
\]
For $I_2$, we have
\begin{eqnarray*}
&&\Big|\int_{\mathbb{T}^3}\Lambda^{3+4\varepsilon_0}\kappa^{(R)} [\Phi(\kappa^{(R)})\frac{\partial \kappa^{(R)}}{\partial z}]dxdydz\Big|\\
&=& \Big|\int_{\mathbb{T}^3}\Lambda^{\frac{5}{2}+2\varepsilon_0}\kappa^{(R)}\Lambda^{\frac{1}{2}+2\varepsilon_0} [\Phi(\kappa^{(R)})\frac{\partial \kappa^{(R)}}{\partial z}]dxdydz\Big|\\
%&\leq &\int_{\mathbb{T}^3}|\Lambda^{2+2\varepsilon_0}\kappa^{(R)}||\Lambda^{2\varepsilon_0}\Phi(\kappa^{(R)})||\nabla_H \kappa^{(R)}|dxdydz+\int_{\mathbb{T}^3}|\Lambda^{2+2\varepsilon_0}\kappa^{(R)}||\Phi(\kappa^{(R)})||\Lambda^{2\varepsilon_0}\nabla_H \kappa^{(R)}|dxdydz\\
&\leq & C|\Lambda^{\frac{5}{2}+2\varepsilon_0}\kappa^{(R)}||\Lambda^{\frac{3}{2}+2\varepsilon_0+\sigma_1}\kappa^{(R)}||\Lambda^{1+\sigma_2} \kappa^{(R)}|+C|\Lambda^{\frac{5}{2}+2\varepsilon_0}\kappa^{(R)}||\Lambda^{1+s_1}\kappa^{(R)}||\Lambda^{\frac{3}{2}+2\varepsilon_0+s_2} \kappa^{(R)}|\\
&\leq & C|\Lambda^{\frac{5}{2}+2\varepsilon_0}\kappa^{(R)}|
|\Lambda^{\frac{3}{2}+2\varepsilon_0}\kappa^{(R)}||\Lambda^2\kappa^{(R)}|\\
&\leq& C|\Lambda^{\frac{5}{2}+2\varepsilon_0}\kappa^{(R)}|^{\frac{3}{2}-2\varepsilon_0}
|\Lambda^{\frac{3}{2}+2\varepsilon_0}\kappa^{(R)}|^{\frac{1}{2}+2\varepsilon_0}|\Lambda^{\frac{3}{2}+2\varepsilon_0}\kappa^{(R)}|\\
&\leq & \varepsilon|\Lambda^{\frac{5}{2}+2\varepsilon_0}\kappa^{(R)}|^2+C|\Lambda^{\frac{3}{2}+2\varepsilon_0}\kappa^{(R)}|^{2+\frac{4}{1+4\varepsilon_0}},
\end{eqnarray*}
where $\sigma_1+\sigma_2=1$, $s_1+s_2=1$, we choose
\[
\sigma_1=s_2=\frac{1}{2}-2\varepsilon_0, \quad \sigma_2=s_1=\frac{1}{2}+2\varepsilon_0.
\]
\begin{eqnarray*}
&&\Big|\int_{\mathbb{T}^3}\Lambda^{3+4\varepsilon_0}\kappa^{(R)} [\Phi(\kappa^{(R)})\frac{\partial Z_1}{\partial z}]dxdydz\Big|\\
&=& \Big|\int_{\mathbb{T}^3}\Lambda^{\frac{5}{2}+2\varepsilon_0}\kappa^{(R)}\Lambda^{\frac{1}{2}+2\varepsilon_0} [\Phi(\kappa^{(R)})\frac{\partial Z_1}{\partial z}]dxdydz\Big|\\
&\leq & C|\Lambda^{\frac{5}{2}+2\varepsilon_0}\kappa^{(R)}||\Lambda^{\frac{3}{2}+2\varepsilon_0+\sigma_1}\kappa^{(R)}||\Lambda^{1+\sigma_2} Z_1|+C|\Lambda^{\frac{5}{2}+2\varepsilon_0}\kappa^{(R)}||\Lambda^{1+s_1}\kappa^{(R)}||\Lambda^{\frac{3}{2}+2\varepsilon_0+s_2} Z_1|\\
&\leq &
|\Lambda^{\frac{5}{2}+2\varepsilon_0}\kappa^{(R)}||\Lambda^{\frac{3}{2}+2\varepsilon_0}Z_1||\Lambda^{2}\kappa^{(R)}|\\
&\leq& |\Lambda^{\frac{5}{2}+2\varepsilon_0}\kappa^{(R)}|^{\frac{3}{2}-2\varepsilon_0}|\Lambda^{\frac{3}{2}+2\varepsilon_0}\kappa^{(R)}|^{\frac{1}{2}+2\varepsilon_0}
|\Lambda^{\frac{3}{2}+2\varepsilon_0}Z_1|^2\\
&\leq&
\varepsilon|\Lambda^{\frac{5}{2}+2\varepsilon_0}\kappa^{(R)}|^2+C|\Lambda^{\frac{3}{2}+2\varepsilon_0}\kappa^{(R)}|^2|\Lambda^{\frac{3}{2}+2\varepsilon_0}Z_1|^{\frac{4}{1+4\varepsilon_0}},
\end{eqnarray*}
where $\sigma_1+\sigma_2=1$, $s_1+s_2=1$, we choose
\[
\sigma_1=\frac{1}{2}-2\varepsilon_0, \quad \sigma_2=\frac{1}{2}+2\varepsilon_0,\quad s_1=1,\quad s_2=0.
\]

%By the changeable of $\kappa^{(R)}$ and $Z_1$ in $\int_{\mathbb{T}^3}\Lambda^{3+4\varepsilon_0}\kappa^{(R)} [\Phi(\kappa^{(R)})\frac{\partial Z_1}{\partial z}]dxdydz$ and $\int_{\mathbb{T}^3}\Lambda^{2+4\varepsilon_0}\kappa^{(R)} [\Phi(Z_1)\frac{\partial \kappa^{(R)}}{\partial z}]dxdydz$, they have the same estimation.
\begin{eqnarray*}
&&\Big|\int_{\mathbb{T}^3}\Lambda^{3+4\varepsilon_0}\kappa^{(R)} [\Phi(Z_1)\frac{\partial \kappa^{(R)}}{\partial z}]dxdydz\Big|\\
&=& \Big|\int_{\mathbb{T}^3}\Lambda^{\frac{5}{2}+2\varepsilon_0}\kappa^{(R)}\Lambda^{\frac{1}{2}+2\varepsilon_0} [\Phi(Z_1)\frac{\partial \kappa^{(R)}}{\partial z}]dxdydz\Big|\\
&\leq & C|\Lambda^{\frac{5}{2}+2\varepsilon_0}\kappa^{(R)}||\Lambda^{\frac{3}{2}+2\varepsilon_0+\sigma_1}Z_1||\Lambda^{1+\sigma_2} \kappa^{(R)}|+C|\Lambda^{\frac{5}{2}+2\varepsilon_0}\kappa^{(R)}||\Lambda^{1+s_1}Z_1||\Lambda^{\frac{3}{2}+2\varepsilon_0+s_2} \kappa^{(R)}|\\
&\leq &
|\Lambda^{\frac{5}{2}+2\varepsilon_0}\kappa^{(R)}||\Lambda^{\frac{3}{2}+2\varepsilon_0}Z_1||\Lambda^{2}\kappa^{(R)}|\\
&\leq& |\Lambda^{\frac{5}{2}+2\varepsilon_0}\kappa^{(R)}|^{\frac{3}{2}-2\varepsilon_0}|\Lambda^{\frac{3}{2}+2\varepsilon_0}\kappa^{(R)}|^{\frac{1}{2}+2\varepsilon_0}
|\Lambda^{\frac{3}{2}+2\varepsilon_0}Z_1|^2\\
&\leq&
\varepsilon|\Lambda^{\frac{5}{2}+2\varepsilon_0}\kappa^{(R)}|^2+C|\Lambda^{\frac{3}{2}+2\varepsilon_0}\kappa^{(R)}|^2|\Lambda^{\frac{3}{2}+2\varepsilon_0}Z_1|^{\frac{4}{1+4\varepsilon_0}},
\end{eqnarray*}
where $\sigma_1+\sigma_2=1$, $s_1+s_2=1$, we choose
\[
\sigma_1=0, \quad \sigma_2=1,\quad s_1=\frac{1}{2}+2\varepsilon_0 ,\quad s_2=\frac{1}{2}-2\varepsilon_0.
\]
\begin{eqnarray*}
&&\Big|\int_{\mathbb{T}^3}\Lambda^{3+4\varepsilon_0}\kappa^{(R)} [\Phi(Z_1)\frac{\partial Z_1}{\partial z}]dxdydz\Big|\\
&=& \Big|\int_{\mathbb{T}^3}\Lambda^{\frac{5}{2}+2\varepsilon_0}\kappa^{(R)}\Lambda^{\frac{1}{2}+2\varepsilon_0} [\Phi(Z_1)\frac{\partial Z_1}{\partial z}]dxdydz\Big|\\
&\leq& C|\Lambda^{\frac{5}{2}+2\varepsilon_0}\kappa^{(R)}||\Lambda^{\frac{3}{2}+2\varepsilon_0+\sigma_1}Z_1||\Lambda^{1+\sigma_2} Z_1|+C|\Lambda^{\frac{5}{2}+2\varepsilon_0}\kappa^{(R)}||\Lambda^{1+s_1}Z_1||\Lambda^{\frac{3}{2}+2\varepsilon_0+s_2} Z_1|\\
&\leq& |\Lambda^{\frac{5}{2}+2\varepsilon_0}\kappa^{(R)}||\Lambda^{\frac{3}{2}+2\varepsilon_0}Z_1||\Lambda^{2} Z_1| \\
&\leq & \varepsilon|\Lambda^{2+2\varepsilon_0}\kappa^{(R)}|^2+C|\Lambda^{\frac{3}{2}+2\varepsilon_0}Z_1|^2|\Lambda^{2} Z_1|^2,
\end{eqnarray*}
where $\sigma_1+\sigma_2=1$, $s_1+s_2=1$, we choose
\[
\sigma_1=0, \quad \sigma_2=1,\quad s_1=1 ,\quad s_2=0.
\]
For $I_3$,
\begin{eqnarray*}
&&\Big|\int_{\mathbb{T}^3}\Lambda^{3+4\varepsilon_0}\kappa^{(R)} [f{k}\times (\kappa^{(R)}+Z_1)+\nabla_H p_{b}-\int^{z}_{-1}\nabla_H (g^{(R)}+Z_2)dz'] dxdydz\Big|\\
&=&\Big|\int_{\mathbb{T}^3}\Lambda^{\frac{5}{2}+2\varepsilon_0}\kappa^{(R)} \Lambda^{\frac{1}{2}+2\varepsilon_0}[f{k}\times (\kappa^{(R)}+Z_1)+\nabla_H p_{b}-\int^{z}_{-1}\nabla_H (g^{(R)}+Z_2)dz'] dxdydz\Big|\\
&\leq & \varepsilon|\Lambda^{\frac{5}{2}+2\varepsilon_0}\kappa^{(R)}|^2+C|\Lambda^{\frac{1}{2}+2\varepsilon_0}\kappa^{(R)}|^2+C|\Lambda^{\frac{1}{2}+2\varepsilon_0}Z_1|^2
+C|\Lambda^{\frac{3}{2}+2\varepsilon_0}g^{(R)}|^2+C|\Lambda^{\frac{3}{2}+2\varepsilon_0}Z_2|^2,
\end{eqnarray*}
thus,
\begin{eqnarray}\label{eq-20}
&&\frac{1}{2}\frac{d}{dt}|\Lambda^{\frac{3}{2}+2\varepsilon_0}\kappa^{(R)}|^2+|\Lambda^{\frac{5}{2}+2\varepsilon_0}\kappa^{(R)}|^2\\ \notag
 &\leq& \varepsilon|\Lambda^{\frac{5}{2}+2\varepsilon_0}\kappa^{(R)}|^2
 +C(R+|\Lambda^{\frac{3}{2}+2\varepsilon_0}Z_1|^{2+\frac{4}{1+4\varepsilon_0}})+C(R+|\Lambda^2 Z_1|^2)|\Lambda^{\frac{3}{2}+2\varepsilon_0}Z_1|^2+|\Lambda^{\frac{3}{2}+2\varepsilon_0}Z_2|^2. \notag
\end{eqnarray}
 Multiplying (\ref{eq-19}) by $-\Lambda^{3+4\varepsilon_0} g^{(R)}$, integrating over $\mathbb{T}^3$, it follows that
\begin{eqnarray*}
 &&\frac{1}{2}\frac{d}{dt}|\Lambda^{\frac{3}{2}+2\varepsilon_0} g^{(R)}|^2+|\Lambda^{\frac{5}{2}+2\varepsilon_0} g^{(R)}|^2\\
& =&\chi_R(|U^{(R)}|^2_{\mathcal{W}})\int_{\mathbb{T}^3}\Lambda^{3+4\varepsilon_0}g^{(R)}[\Big((\kappa^{(R)}+ Z_1)\cdot \nabla_H\Big)(g^{(R)}+ Z_2)]dxdydz\\ &+&\chi_R(|U^{(R)}|^2_{\mathcal{W}})\int_{\mathbb{T}^3}\Lambda^{3+4\varepsilon_0}g^{(R)}[\Phi(\kappa^{(R)}+Z_1)\frac{\partial( g^{(R)}+Z_2)}{\partial z}]dxdydz\\
&+&\chi_R(|U^{(R)}|^2_{\mathcal{W}})\int_{\mathbb{T}^3}\Lambda^{3+4\varepsilon_0}g^{(R)}\Phi(\kappa^{(R)}+Z_1)dxdydz\\
&:=& \chi_R(|U^{(R)}|^2_{\mathcal{W}})(I_4+I_5+I_6).
\end{eqnarray*}
For $I_4$, we have
\begin{eqnarray*}
&&\Big|\int_{\mathbb{T}^3}\Lambda^{3+4\varepsilon_0}g^{(R)}[(\kappa^{(R)}\cdot \nabla_H)g^{(R)}]dxdydz\Big|\\
&=& \Big|\int_{\mathbb{T}^3}\Lambda^{\frac{5}{2}+2\varepsilon_0}g^{(R)}\Lambda^{\frac{1}{2}+2\varepsilon_0}[(\kappa^{(R)}\cdot \nabla_H)g^{(R)}]dxdydz\Big|\\
&\leq& C|\Lambda^{\frac{5}{2}+2\varepsilon_0}g^{(R)}||\Lambda^{\frac{1}{2}+2\varepsilon_0+\sigma_1}\kappa^{(R)}||\Lambda^{1+\sigma_2}g^{(R)}|
+C|\Lambda^{\frac{5}{2}+2\varepsilon_0}g^{(R)}||\Lambda^{s_1}\kappa^{(R)}||\Lambda^{\frac{3}{2}+2\varepsilon_0+s_2}g^{(R)}|\\
&\leq& C|\Lambda^{\frac{5}{2}+2\varepsilon_0}g^{(R)}||\Lambda^{\frac{3}{2}+2\varepsilon_0}g^{(R)}||\Lambda^{\frac{3}{2}}\kappa^{(R)}|\\
&\leq&  \varepsilon |\Lambda^{\frac{5}{2}+2\varepsilon_0}g^{(R)}|^2+C|\Lambda^{\frac{3}{2}+2\varepsilon_0}g^{(R)}|^2|\Lambda^{\frac{3}{2}}\kappa^{(R)}|^2,
\end{eqnarray*}
where $\sigma_1+\sigma_2=\frac{3}{2}$, $s_1+s_2=\frac{3}{2}$, we choose
\[
\sigma_1=1-2\varepsilon_0, \quad \sigma_2=\frac{1}{2}+2\varepsilon_0, \quad s_1=\frac{3}{2}, \quad s_2=0.
\]
 By the same argument, we have
\begin{eqnarray*}
\Big|\int_{\mathbb{T}^3}\Lambda^{3+4\varepsilon_0}g^{(R)}[(\kappa^{(R)}\cdot \nabla_H)Z_2]dxdydz\Big|
&\leq& \varepsilon |\Lambda^{\frac{5}{2}+2\varepsilon_0}g^{(R)}|^2+C|\Lambda^{\frac{3}{2}+2\varepsilon_0}\kappa^{(R)}|^2|\Lambda^{\frac{3}{2}}Z_2|^2
+C|\Lambda^{\frac{3}{2}+2\varepsilon_0}Z_2|^2|\Lambda^{\frac{3}{2}}\kappa^{(R)}|^2,\\
\Big|\int_{\mathbb{T}^3}\Lambda^{3+4\varepsilon_0}g^{(R)}[(Z_1\cdot \nabla_H)g^{(R)}]dxdydz\Big|
&\leq& \varepsilon |\Lambda^{\frac{5}{2}+2\varepsilon_0}g^{(R)}|^2+C|\Lambda^{\frac{3}{2}+2\varepsilon_0}g^{(R)}|^2|\Lambda^{\frac{3}{2}}Z_1|^2,\\
\Big|\int_{\mathbb{T}^3}\Lambda^{3+4\varepsilon_0}g^{(R)}[(Z_1\cdot \nabla_H)Z_2]dxdydz\Big|
&\leq& \varepsilon |\Lambda^{\frac{5}{2}+2\varepsilon_0}g^{(R)}|^2+C|\Lambda^{\frac{3}{2}+2\varepsilon_0}Z_2|^2|\Lambda^{\frac{3}{2}}Z_1|^2.
\end{eqnarray*}
For $I_5$, similar to  $I_4$,
\begin{eqnarray*}
\Big|\int_{\mathbb{T}^3}\Lambda^{3+4\varepsilon_0}g^{(R)}[\Phi(\kappa^{(R)})\frac{\partial g^{(R)}}{\partial z}]dxdydz\Big|
&\leq& \varepsilon |\Lambda^{\frac{5}{2}+2\varepsilon_0}g^{(R)}|^2+C|\Lambda^{\frac{3}{2}+2\varepsilon_0}g^{(R)}|^2|\Lambda^{\frac{3}{2}+2\varepsilon_0}\kappa^{(R)}|^{\frac{4}{1+4\varepsilon_0}},\\
\Big|\int_{\mathbb{T}^3}\Lambda^{3+4\varepsilon_0}g^{(R)}[\Phi(\kappa^{(R)})\frac{\partial Z_2}{\partial z}]dxdydz\Big|
&\leq& \varepsilon |\Lambda^{\frac{5}{2}+2\varepsilon_0}g^{(R)}|^2+C|\Lambda^{\frac{3}{2}+2\varepsilon_0}\kappa^{(R)}|^2|\Lambda^{2}Z_2|^{2},\\
\Big|\int_{\mathbb{T}^3}\Lambda^{3+4\varepsilon_0}g^{(R)}[\Phi(Z_1)\frac{\partial g^{(R)}}{\partial z}]dxdydz\Big|
&\leq& \varepsilon |\Lambda^{\frac{5}{2}+2\varepsilon_0}g^{(R)}|^2+C|\Lambda^{\frac{3}{2}+2\varepsilon_0}g^{(R)}|^2|\Lambda^{\frac{3}{2}+2\varepsilon_0}Z_1|^{\frac{4}{1+4\varepsilon_0}},\\
\Big|\int_{\mathbb{T}^3}\Lambda^{3+4\varepsilon_0}g^{(R)}[\Phi(Z_1)\frac{\partial Z_2}{\partial z}]dxdydz\Big|
&\leq& \varepsilon |\Lambda^{\frac{5}{2}+2\varepsilon_0}g^{(R)}|^2+C|\Lambda^{\frac{3}{2}+2\varepsilon_0}Z_1|^2|\Lambda^{2}Z_2|^{2}.
\end{eqnarray*}
For $I_6$, we have
\begin{eqnarray*}
&&\Big|\int_{\mathbb{T}^3}\Lambda^{3+4\varepsilon_0}g^{(R)}\Phi(\kappa^{(R)+Z_1})dxdydz\Big|\\
&\leq& C|\Lambda^{\frac{5}{2}+2\varepsilon_0}g^{(R)}||\Lambda^{\frac{3}{2}+2\varepsilon_0}(\kappa^{(R)}+Z_1)|\\
&\leq& \varepsilon |\Lambda^{\frac{5}{2}+2\varepsilon_0}g^{(R)}|^2
+C|\Lambda^{\frac{3}{2}+2\varepsilon_0}\kappa^{(R)}|^2+C|\Lambda^{\frac{3}{2}+2\varepsilon_0}Z_1|^{2}.
\end{eqnarray*}
Thus,
\begin{eqnarray}\label{eq-21}
&&\frac{1}{2}\frac{d}{dt}|\Lambda^{\frac{3}{2}+2\varepsilon_0}g^{(R)}|^2+|\Lambda^{\frac{5}{2}+2\varepsilon_0}g^{(R)}|^2
\\ \notag
& \leq& \varepsilon|\Lambda^{\frac{5}{2}+2\varepsilon_0}g^{(R)}|^2
 +C(R+|\Lambda^{\frac{3}{2}+2\varepsilon_0}Z_1|^{\frac{4}{1+4\varepsilon_0}}|\Lambda^{\frac{3}{2}+2\varepsilon_0}Z_2|^2
 +|\Lambda^{\frac{3}{2}+2\varepsilon_0}Z_1|^{2}|\Lambda^{2}Z_2|^2). \notag
\end{eqnarray}
Since $|X^{(R)}|^2_{\mathcal{W}}=|\kappa^{(R)}|^2_{\mathcal{W}}+|g^{(R)}|^2_{\mathcal{W}}$, combining (\ref{eq-20}) and (\ref{eq-21}), we have
\begin{eqnarray}\label{eq-26}
\frac{d|X^{(R)}|^2_{\mathcal{W}}}{dt}+\|X^{(R)}\|^2_{\frac{5}{2}+2\varepsilon_0}\leq C(R+|\Lambda^{2}Z|^{2+\frac{4}{1+4\varepsilon_0}}),
\end{eqnarray}
by the property of $Z$ in (\ref{eq-22}), (\ref{eq-17}) has a weak solution $X^{(R)}$ in $L^{\infty}([0,T];\mathcal{W})\cap L^2([0,T]; D(A^{\frac{5}{4}+\varepsilon_0}))$.

\textbf{(Continuity of weak solution)} \quad
Multiplying (\ref{eq-18}) by $-\Lambda^{1+4\varepsilon_0} \frac{ d \kappa^{(R)}}{dt}$, integrating over $\mathbb{T}^3$, it follows that
\begin{eqnarray*}
 &&\frac{1}{2}\frac{d}{dt}\int_{\mathbb{T}^3}|\Lambda^{\frac{3}{2}+2\varepsilon_0}\kappa^{(R)}|^2dxdydz+| \Lambda^{\frac{1}{2}+2\varepsilon_0} \dot{\kappa}^{(R)}|^2\\
 &=&\chi_R(|U^{(R)}|^2_{\mathcal{W}})\int_{\mathbb{T}^3}\Big[\Big((\kappa^{(R)}+ Z_1)\cdot \nabla_H\Big)(\kappa^{(R)}+ Z_1)\Big] \Lambda^{1+4\varepsilon_0} \dot{\kappa}^{(R)}dxdydz\\ \notag
 &+&\chi_R(|U^{(R)}|^2_{\mathcal{W}})\int_{\mathbb{T}^3}\Big[\Phi(\kappa^{(R)}+Z_1)\frac{\partial( \kappa^{(R)}+ Z_1)}{\partial z}\Big] \Lambda^{1+4\varepsilon_0} \dot{\kappa}^{(R)}dxdydz\\ \notag
 &+&\chi_R(|U^{(R)}|^2_{\mathcal{W}})\int_{\mathbb{T}^3}\Big(f{k}\times (\kappa^{(R)}+Z_1)+\nabla_H p_{b}-\int^{z}_{-1}\nabla_H (g^{(R)}+Z_2)dz'\Big) \Lambda^{1+4\varepsilon_0} \dot{\kappa}^{(R)}dxdydz\\
 &:=& \chi_R(|U^{(R)}|^2_{\mathcal{W}})(J_1+J_2+J_3),
 \end{eqnarray*}
 where $\dot{\kappa}^{(R)}$ denotes $\frac{d \kappa^{(R)}}{dt}$, and this symbol always denotes the deviation with respect to $t$.

For $J_1$,
By H\"{o}lder inequality and the Young inequality, we have
\begin{eqnarray*}
&&\Big|\int_{\mathbb{T}^3}[(\kappa^{(R)}\cdot \nabla_H)\kappa^{(R)}] \Lambda^{1+4\varepsilon_0} \dot{\kappa}^{(R)}dxdydz\Big|\\
&\leq&C|\Lambda^{\frac{1}{2}+2\varepsilon_0} \dot{\kappa}^{(R)}||\Lambda^{\frac{1}{2}+2\varepsilon_0+\sigma_1} \kappa^{(R)}||\Lambda^{1+\sigma_2}\kappa^{(R)}|+C|\Lambda^{\frac{1}{2}+2\varepsilon_0} \dot{\kappa}^{(R)}||\Lambda^{s_1} \kappa^{(R)}||\Lambda^{\frac{3}{2}+2\varepsilon_0+s_2}\kappa^{(R)}|\\
&\leq&C|\Lambda^{\frac{1}{2}+2\varepsilon_0} \dot{\kappa}^{(R)}||\Lambda^{\frac{3}{2}+2\varepsilon_0} \kappa^{(R)}||\Lambda^{\frac{3}{2}}\kappa^{(R)}|\\
&\leq&\varepsilon |\Lambda^{\frac{1}{2}+2\varepsilon_0} \dot{\kappa}^{(R)}|^2+C|\Lambda^{\frac{3}{2}+2\varepsilon_0}\kappa^{(R)}|^4,
\end{eqnarray*}
similarly,
\begin{eqnarray*}
&&\Big|\int_{\mathbb{T}^3}[(\kappa^{(R)}\cdot \nabla_H) Z_1] \Lambda^{1+4\varepsilon_0} \dot{\kappa}^{(R)}dxdydz\Big|\\
&\leq&C|\Lambda^{\frac{1}{2}+2\varepsilon_0} \dot{\kappa}^{(R)}||\Lambda^{\frac{1}{2}+2\varepsilon_0+\sigma_1} \kappa^{(R)}||\Lambda^{1+\sigma_2}Z_1|+C|\Lambda^{\frac{1}{2}+2\varepsilon_0} \dot{\kappa}^{(R)}||\Lambda^{s_1} \kappa^{(R)}||\Lambda^{\frac{3}{2}+2\varepsilon_0+s_2}Z_1|\\
&\leq&C|\Lambda^{\frac{1}{2}+2\varepsilon_0} \dot{\kappa}^{(R)}|\Big(|\Lambda^{\frac{3}{2}+2\varepsilon_0} \kappa^{(R)}||\Lambda^{\frac{3}{2}}Z_1|+|\Lambda^{\frac{3}{2}+2\varepsilon_0} Z_1||\Lambda^{\frac{3}{2}}\kappa^{(R)}|\Big)\\
&\leq&\varepsilon |\Lambda^{\frac{1}{2}+2\varepsilon_0} \dot{\kappa}^{(R)}|^2+C|\Lambda^{\frac{3}{2}+2\varepsilon_0}\kappa^{(R)}|^2|\Lambda^{\frac{3}{2}+2\varepsilon_0}Z_1|^2,
\end{eqnarray*}
\begin{eqnarray*}
&&\Big|\int_{\mathbb{T}^3}[(Z_1\cdot \nabla_H)\kappa^{(R)}] \Lambda^{1+4\varepsilon_0} \dot{\kappa}^{(R)}dxdydz\Big|\\
&\leq&C|\Lambda^{\frac{1}{2}+2\varepsilon_0} \dot{\kappa}^{(R)}||\Lambda^{\frac{1}{2}+2\varepsilon_0+\sigma_1} Z_1||\Lambda^{1+\sigma_2}\kappa^{(R)}|+C|\Lambda^{\frac{1}{2}+2\varepsilon_0} \dot{\kappa}^{(R)}||\Lambda^{s_1} Z_1||\Lambda^{\frac{3}{2}+2\varepsilon_0+s_2}\kappa^{(R)}|\\
&\leq&C|\Lambda^{\frac{1}{2}+2\varepsilon_0} \dot{\kappa}^{(R)}| \Big(|\Lambda^{\frac{3}{2}+2\varepsilon_0}Z_1||\Lambda^{\frac{3}{2}}\kappa^{(R)}|
+|\Lambda^{\frac{3}{2}}Z_1||\Lambda^{\frac{3}{2}+2\varepsilon_0}\kappa^{(R)}|\Big) \\
&\leq&\varepsilon |\Lambda^{\frac{1}{2}+2\varepsilon_0} \dot{\kappa}^{(R)}|^2+C|\Lambda^{\frac{3}{2}+2\varepsilon_0}Z_1|^2|\Lambda^{\frac{3}{2}}\kappa^{(R)}|^2
+C|\Lambda^{\frac{3}{2}}Z_1|^2|\Lambda^{\frac{3}{2}+2\varepsilon_0}\kappa^{(R)}|^2,
\end{eqnarray*}
and
\begin{eqnarray*}
&&\Big|\int_{\mathbb{T}^3}[(Z_1\cdot \nabla_H)Z_1] \Lambda^{1+4\varepsilon_0} \dot{\kappa}^{(R)}dxdydz\Big|\\
&\leq&C|\Lambda^{\frac{1}{2}+2\varepsilon_0} \dot{\kappa}^{(R)}||\Lambda^{\frac{1}{2}+2\varepsilon_0+\sigma_1} Z_1||\Lambda^{1+\sigma_2}Z_1|+C|\Lambda^{\frac{1}{2}+2\varepsilon_0} \dot{\kappa}^{(R)}||\Lambda^{s_1} Z_1||\Lambda^{\frac{3}{2}+2\varepsilon_0+s_2}Z_1|\\
&\leq& C|\Lambda^{\frac{1}{2}+2\varepsilon_0} \dot{\kappa}^{(R)}||\Lambda^{\frac{3}{2}+2\varepsilon_0}Z_1||\Lambda^{\frac{3}{2}}Z_1|\\
&\leq&\varepsilon |\Lambda^{\frac{1}{2}+2\varepsilon_0} \dot{\kappa}^{(R)}|^2+C|\Lambda^{\frac{3}{2}+2\varepsilon_0}Z_1|^2|\Lambda^{\frac{3}{2}}Z_1|^2.
\end{eqnarray*}

For $J_2$, we obtain
\begin{eqnarray*}
&&\Big|\int_{\mathbb{T}^3}[\Phi(\kappa^{(R)})\frac{\partial \kappa^{(R)}}{\partial z}] \Lambda^{1+4\varepsilon_0} \dot{\kappa}^{(R)}dxdydz\Big|\\
&\leq&C|\Lambda^{\frac{1}{2}+2\varepsilon_0} \dot{\kappa}^{(R)}|\Big(|\Lambda^{\frac{3}{2}+2\varepsilon_0+\sigma_1} \kappa^{(R)}||\Lambda^{1+\sigma_2}\kappa^{(R)}|+|\Lambda^{1+s_1}\kappa^{(R)}||\Lambda^{\frac{3}{2}+2\varepsilon_0+s_2} \kappa^{(R)}|\Big)\\
&\leq&C|\Lambda^{\frac{1}{2}+2\varepsilon_0} \dot{\kappa}^{(R)}||\Lambda^{2} \kappa^{(R)}||\Lambda^{\frac{3}{2}+2\varepsilon_0}\kappa^{(R)}|\\
&\leq& C|\Lambda^{\frac{1}{2}+2\varepsilon_0} \dot{\kappa}^{(R)}||\Lambda^{\frac{3}{2}+2\varepsilon_0} \kappa^{(R)}|^{\frac{1}{2}+2\varepsilon_0}|\Lambda^{\frac{5}{2}+2\varepsilon_0} \kappa^{(R)}|^{\frac{1}{2}-2\varepsilon_0}|\Lambda^{\frac{3}{2}+2\varepsilon_0}\kappa^{(R)}|\\
&\leq&\varepsilon |\Lambda^{\frac{1}{2}+2\varepsilon_0} \dot{\kappa}^{(R)}|^2+C|\Lambda^{\frac{5}{2}+2\varepsilon_0}\kappa^{(R)}|^2
+C|\Lambda^{\frac{3}{2}+2\varepsilon_0}\kappa^{(R)}|^2
|\Lambda^{\frac{3}{2}+2\varepsilon_0}\kappa^{(R)}|^{\frac{4}{1+4\varepsilon_0}},
\end{eqnarray*}
similarly,
\begin{eqnarray*}
&&\Big|\int_{\mathbb{T}^3}[\Phi(\kappa^{(R)})\frac{\partial Z_1}{\partial z}] \Lambda^{1+4\varepsilon_0} \dot{\kappa}^{(R)}dxdydz\Big|\\
&\leq&C|\Lambda^{\frac{1}{2}+2\varepsilon_0} \dot{\kappa}^{(R)}|\Big(|\Lambda^{\frac{3}{2}+2\varepsilon_0+\sigma_1} \kappa^{(R)}||\Lambda^{1+\sigma_2}Z_1|+|\Lambda^{1+s_1}\kappa^{(R)}||\Lambda^{\frac{3}{2}+2\varepsilon_0+s_2} Z_1|\Big)\\
&\leq&C|\Lambda^{\frac{1}{2}+2\varepsilon_0} \dot{\kappa}^{(R)}||\Lambda^{2} Z_1||\Lambda^{\frac{3}{2}+2\varepsilon_0}\kappa^{(R)}|\\
&\leq&\varepsilon |\Lambda^{\frac{1}{2}+2\varepsilon_0} \dot{\kappa}^{(R)}|^2
+C|\Lambda^{\frac{3}{2}+2\varepsilon_0}\kappa^{(R)}|^2|\Lambda^{2}Z_1|^2,\\
&&\Big|\int_{\mathbb{T}^3}[\Phi(Z_1)\frac{\partial \kappa^{(R)}}{\partial z}] \Lambda^{1+4\varepsilon_0} \dot{\kappa}^{(R)}dxdydz\Big|\\
&\leq&C|\Lambda^{\frac{1}{2}+2\varepsilon_0} \dot{\kappa}^{(R)}||\Lambda^{\frac{3}{2}+2\varepsilon_0} \kappa^{(R)}||\Lambda^{2}Z_1|\\
&\leq&\varepsilon |\Lambda^{\frac{1}{2}+2\varepsilon_0} \dot{\kappa}^{(R)}|^2
+C|\Lambda^{\frac{3}{2}+2\varepsilon_0}\kappa^{(R)}|^2|\Lambda^{2}Z_1|^2,\\
&&\Big|\int_{\mathbb{T}^3}[\Phi(Z_1)\frac{\partial Z_1}{\partial z}] \Lambda^{1+4\varepsilon_0} \dot{\kappa}^{(R)}dxdydz\Big|\\
&\leq&C|\Lambda^{\frac{1}{2}+2\varepsilon_0} \dot{\kappa}^{(R)}||\Lambda^{\frac{3}{2}+2\varepsilon_0} Z_1||\Lambda^{2}Z_1|\\
&\leq&\varepsilon |\Lambda^{\frac{1}{2}+2\varepsilon_0} \dot{\kappa}^{(R)}|^2+C|\Lambda^{\frac{3}{2}
+2\varepsilon_0}Z_1|^2|\Lambda^{2}Z_1|^2.
\end{eqnarray*}
For $J_3$, we have
\begin{eqnarray*}
&&\Big|\int_{\mathbb{T}^3}[f{k}\times (\kappa^{(R)}+Z_1)+\nabla_H p_{b}-\int^{z}_{-1}\nabla_H (g^{(R)}+Z_2)dz'] \Lambda^{1+4\varepsilon_0} \dot{\kappa}^{(R)}dxdydz\Big|\\
&\leq & \varepsilon |\Lambda^{\frac{1}{2}+2\varepsilon_0} \dot{\kappa}^{(R)}|^2+C|\Lambda^{\frac{1}{2}+2\varepsilon_0}(\kappa^{(R)}+Z_1)|^2+C|\Lambda^{\frac{3}{2}+2\varepsilon_0}(g^{(R)}+Z_2)|^2.
\end{eqnarray*}
Multiplying (\ref{eq-19}) by $-\Lambda^{1+4\varepsilon_0} \frac{d}{dt} g^{(R)}$, integrating over $\mathbb{T}^3$, it follows that
\begin{eqnarray*}
 &&\frac{1}{2}\frac{d}{dt}\int_{\mathbb{T}^3}|\Lambda^{\frac{3}{2}+2\varepsilon_0} g^{(R)}|^2dxdydz+|\Lambda^{\frac{1}{2}+2\varepsilon_0}\dot{g}^{(R)}|^2\\
 &=&\chi_R(|U^{(R)}|^2_{\mathcal{W}})\int_{\mathbb{T}^3}\Big[\Big((\kappa^{(R)}+ Z_1)\cdot \nabla_H\Big)(g^{(R)}+ Z_2)\Big] \Lambda^{1+4\varepsilon_0}\dot{g}^{(R)}dxdydz,\\
 &+&\chi_R(|U^{(R)}|^2_{\mathcal{W}})\int_{\mathbb{T}^3}\Big[\Phi(\kappa^{(R)}+Z_1)\frac{\partial( g^{(R)}+Z_2)}{\partial z}\Big]\Lambda^{1+4\varepsilon_0}\dot{g}^{(R)}dxdydz,\\
&+&\chi_R(|U^{(R)}|^2_{\mathcal{W}})\int_{\mathbb{T}^3}\Phi(\kappa^{(R)}+Z_1)\Lambda^{1+4\varepsilon_0}\dot{g}^{(R)}dxdydz,\\
 &:=& \chi_R(|U^{(R)}|^2_{\mathcal{W}})(J_4+J_5+J_6),
 \end{eqnarray*}
 where $\dot{g}^{(R)}$ denotes $\frac{d g^{(R)}}{dt}$.

 For $J_4$, we have
 \begin{eqnarray*}
&&\Big|\int_{\mathbb{T}^3}[(\kappa^{(R)}\cdot \nabla_H)g^{(R)}] \Lambda^{1+4\varepsilon_0}\dot{g}^{(R)}dxdydz\Big|\\
&\leq& C|\Lambda^{\frac{1}{2}+2\varepsilon_0} \dot{g}^{(R)}|\Big(|\Lambda^{\frac{1}{2}+2\varepsilon_0+\sigma_1}\kappa^{(R)}||\Lambda^{1+\sigma_2}g^{(R)}|
+|\Lambda^{s_1}\kappa^{(R)}||\Lambda^{\frac{3}{2}+2\varepsilon_0+s_2}g^{(R)}|\Big)\\
&\leq& C|\Lambda^{\frac{1}{2}+2\varepsilon_0} \dot{g}^{(R)}||\Lambda^{\frac{3}{2}+2\varepsilon_0} g^{(R)}||\Lambda^{\frac{3}{2}}\kappa^{(R)}|\\
&\leq & \varepsilon |\Lambda^{\frac{1}{2}+2\varepsilon_0} \dot{g}^{(R)}|^2+C|\Lambda^{\frac{3}{2}+2\varepsilon_0} g^{(R)}|^2|\Lambda^{\frac{3}{2}}\kappa^{(R)}|^2,
\end{eqnarray*}
similarly,
\begin{eqnarray*}
\Big|\int_{\mathbb{T}^3}[(\kappa^{(R)}\cdot \nabla_H) Z_2] \Lambda^{1+4\varepsilon_0}\dot{g}^{(R)}dxdydz\Big|
&\leq & \varepsilon |\Lambda^{\frac{1}{2}+2\varepsilon_0} \dot{g}^{(R)}|^2+C|\Lambda^{\frac{3}{2}+2\varepsilon_0}Z_2|^2|\Lambda^{\frac{3}{2}}\kappa^{(R)}|^2,\\
\Big|\int_{\mathbb{T}^3}[(Z_1\cdot \nabla_H) g^{(R)}] \Lambda^{1+4\varepsilon_0}\dot{g}^{(R)}dxdydz\Big|
&\leq& \varepsilon |\Lambda^{\frac{1}{2}+2\varepsilon_0} \dot{g}^{(R)}|^2+C|\Lambda^{\frac{3}{2}+2\varepsilon_0}g^{(R)}|^2|\Lambda^{\frac{3}{2}}Z_1|^2,\\
\Big|\int_{\mathbb{T}^3}[(Z_1\cdot \nabla_H) Z_2] \Lambda^{4\varepsilon_0}\dot{g}^{(R)}dxdydz\Big|
&\leq& \varepsilon |\Lambda^{\frac{1}{2}+2\varepsilon_0} \dot{g}^{(R)}|^2+C|\Lambda^{\frac{3}{2}+2\varepsilon_0}Z_2|^2|\Lambda^{\frac{3}{2}}Z_1|^2.
\end{eqnarray*}
%\begin{eqnarray*}
%&&|\int_{\mathbb{T}^3}[(Z_1\cdot \nabla_H) g^{(R)}] %\Lambda^{4\varepsilon_0}\dot{g}^{(R)}dxdydz|\\
%&\leq& C|\Lambda^{2\varepsilon_0} \dot{\kappa}^{(R)}||\Lambda^{\frac{3}{2}+2\varepsilon_0} g^{(R)}||\Lambda Z_1|\\
%&\leq & \varepsilon |\Lambda^{2\varepsilon_0} \dot{g^{(R)}}^{R}|^2+C|\Lambda^{2+2\varepsilon_0} g^{(R)}|^2+C|\Lambda^{1+2\varepsilon_0}g^{(R)}|^2|\Lambda Z_1|^4,
%\end{eqnarray*}
%\begin{eqnarray*}
%&&|\int_{\mathbb{T}^3}[(Z_1\cdot \nabla_H) Z_2] \Lambda^{4\varepsilon_0}\dot{g}^{(R)}dxdydz|\\
%&\leq & \varepsilon |\Lambda^{2\varepsilon_0} \dot{g^{(R)}}^{R}|^2+C|\Lambda^{\frac{5}{4}+\varepsilon_0}\kappa^{(R)}|^2|\Lambda^{\frac{5}{4}+\varepsilon_0}Z_2|^2,
%\end{eqnarray*}
For $J_5$, we have
\begin{eqnarray*}
&&\Big|\int_{\mathbb{T}^3}[\Phi(\kappa^{(R)})\frac{\partial g^{(R)}}{\partial z}] \Lambda^{1+4\varepsilon_0}\dot{g}^{(R)}dxdydz\Big|\\
&\leq &|\Lambda^{\frac{1}{2}+2\varepsilon_0} \dot{g}^{(R)}||\Lambda^{\frac{3}{2}+2\varepsilon_0+\sigma_1} \kappa^{(R)}||\Lambda^{1+\sigma_2} g^{(R)}|^2+C|\Lambda^{2\varepsilon_0} \dot{g}^{(R)}||\Lambda^{1+s_1} \kappa^{(R)}||\Lambda^{\frac{3}{2}+2\varepsilon_0+s_2} g^{(R)}|^2\\
&\leq & C|\Lambda^{\frac{1}{2}+2\varepsilon_0} \dot{g}^{(R)}||\Lambda^{\frac{3}{2}+2\varepsilon_0} \kappa^{(R)}||\Lambda^{2} g^{(R)}|^2\\
&\leq & \varepsilon |\Lambda^{\frac{1}{2}+2\varepsilon_0} \dot{g}^{(R)}|^2+C|\Lambda^{\frac{5}{2}+2\varepsilon_0} g^{(R)}|^2+C|\Lambda^{\frac{3}{2}+2\varepsilon_0}g^{(R)}|^2
|\Lambda^{\frac{3}{2}+2\varepsilon_0}\kappa^{(R)}|^{\frac{4}{1+4\varepsilon_0}},
\end{eqnarray*}
similarly, we get
 \begin{eqnarray*}
\Big|\int_{\mathbb{T}^3}[\Phi(\kappa^{(R)})\frac{\partial Z_2}{\partial z}] \Lambda^{1+4\varepsilon_0}\dot{g}^{(R)}dxdydz\Big|
&\leq &\varepsilon |\Lambda^{\frac{1}{2}+2\varepsilon_0} \dot{g}^{(R)}|^2+C|\Lambda^{\frac{5}{2}+2\varepsilon_0} \kappa^{(R)}|^2+C|\Lambda^{\frac{3}{2}+2\varepsilon_0}\kappa^{(R)}|^2
|\Lambda^{\frac{3}{2}+2\varepsilon_0}Z_2|^{\frac{4}{1+4\varepsilon_0}},\\
|\int_{\mathbb{T}^3}[\Phi(Z_1)\frac{\partial Z_2}{\partial z}] \Lambda^{1+4\varepsilon_0}\dot{g}^{(R)}dxdydz|
&\leq & \varepsilon |\Lambda^{\frac{1}{2}+2\varepsilon_0} \dot{g}^{(R)}|^2+C|\Lambda^{\frac{5}{2}+2\varepsilon_0} g^{(R)}|^2+C|\Lambda^{\frac{3}{2}+2\varepsilon_0}g^{(R)}|^2
|\Lambda^{\frac{3}{2}+2\varepsilon_0}Z_1|^{\frac{4}{1+4\varepsilon_0}},\\
|\int_{\mathbb{T}^3}[\Phi(Z_1)\frac{\partial Z_2}{\partial z}] \Lambda^{1+4\varepsilon_0}\dot{g}^{(R)}dxdydz|
&\leq & \varepsilon |\Lambda^{\frac{1}{2}+2\varepsilon_0} \dot{g}^{(R)}|^2+C|\Lambda^{2} Z_1|^2|\Lambda^{\frac{3}{2}+2\varepsilon_0}Z_2|^2.
\end{eqnarray*}
For $J_6$, we have
\begin{eqnarray*}
&&\Big|\int_{\mathbb{T}^3}\Phi(\kappa^{(R)}+Z_1) \Lambda^{1+4\varepsilon_0}\dot{g}^{(R)}dxdydz\Big|\\
&\leq & C|\Lambda^{\frac{1}{2}+2\varepsilon_0} \dot{g}^{(R)}||\Lambda^{\frac{3}{2}+2\varepsilon_0} (\kappa^{(R)}+Z_1)|\\
&\leq & \varepsilon |\Lambda^{\frac{1}{2}+2\varepsilon_0} \dot{g}^{(R)}|^2+C|\Lambda^{\frac{3}{2}+2\varepsilon_0}\kappa^{(R)}|^2
+C|\Lambda^{\frac{3}{2}+2\varepsilon_0}Z_1|^{2}.
\end{eqnarray*}
  %\begin{eqnarray*}
%&&|\int_{\mathbb{T}^3}[\Phi(Z_1)\frac{\partial g^{(R)}}{\partial z}] \Lambda^{4\varepsilon_0}\dot{g}^{(R)}dxdydz|\\
%&\leq &|\Lambda^{2\varepsilon_0} \dot{g^{(R)}}^{R}||\Lambda^{1+2\varepsilon_0+\sigma_1} Z_1||\Lambda^{1+\sigma_2} g^{(R)}|\\
%&\leq & C|\Lambda^{2\varepsilon_0} \dot{g^{(R)}}^{R}||\Lambda^{\frac{3}{2}+2\varepsilon_0} g^{(R)}||\Lambda^{\frac{3}{2}} z_1|\\
%&\leq & \varepsilon |\Lambda^{2\varepsilon_0} \dot{g^{(R)}}^{R}|^2+C|\Lambda^{2+2\varepsilon_0} g^{(R)}|^2+C|\Lambda^{1+2\varepsilon_0}\kappa^{(R)}|^2|\Lambda^{\frac{3}{2}} Z_1|^4,
%\end{eqnarray*}
%and
%\begin{eqnarray*}
%&&|\int_{\mathbb{T}^3}[\Phi(Z_1)\frac{\partial Z_2}{\partial z}] \Lambda^{4\varepsilon_0}\dot{g}^{(R)}dxdydz|\\
%&\leq & \varepsilon |\Lambda^{2\varepsilon_0} \dot{g^{(R)}}^{R}|^2+C|\Lambda^{\frac{3}{2}+\varepsilon_0}Z_1|^2|\Lambda^{\frac{3}{2}+\varepsilon_0} Z_2|^4,
%\end{eqnarray*}
Combing the above estimations, we obtain
\begin{eqnarray}\label{eq-25}
&&\frac{d}{dt}\|X^{(R)}\|^2_{\frac{3}{2}+2\varepsilon_0}+\|\dot{X}^{(R)}\|^2_{\frac{1}{2}+2\varepsilon_0}\\ \notag
&\leq&
C\chi_R(|U^{(R)}|^2_{\mathcal{W}})\Big(C(R)+|\Lambda^{\frac{5}{2}+2\varepsilon_0} \kappa^{(R)}|^2+|\Lambda^{\frac{5}{2}+2\varepsilon_0} g^{(R)}|^2+|\Lambda^{2} Z_1|^{2+\frac{4}{1+4\varepsilon_0}}+|\Lambda^{2} Z_2|^{2+\frac{4}{1+4\varepsilon_0}}\Big).
\end{eqnarray}
Integrating (\ref{eq-25}) on $t$ from $0$ to $T$, as $\int^T_0|\Lambda^{\frac{5}{2}+2\varepsilon_0} \kappa^{(R)}|^2dt$ and $\int^T_0|\Lambda^{\frac{5}{2}+2\varepsilon_0} g^{(R)}|^2dt$ can be dominated by (\ref{eq-26}), by the property of $Z$ in (\ref{eq-22}), we get the time derivative $\frac{d{X}^{(R)}}{dt} \in L^2([0,T]; D(A^{\frac{1}{4}+\varepsilon_0}))$. Then by (\ref{eq-26}) and \cite{T-R}, we obtain $X^{(R)}\in C([0,T]; \mathcal{W})$.

\textbf{(Uniqueness of weak solution)}\quad Let $X_1=(\kappa^{(R)}_1, g^{(R)}_1)$, $X_2=(\kappa^{(R)}_2, g^{(R)}_2)$ be two solutions of (\ref{eq-17}) in $C([0,T]; \mathcal{W})$  and set
\begin{eqnarray*}
Y=X_1-X_2=(\kappa,g)=(\kappa^{(R)}_1-\kappa^{(R)}_2,g^{(R)}_1-g^{(R)}_2),\\
v_1=\kappa^{(R)}_1+Z_1, \quad T_1=g^{(R)}_1+Z_2,\quad v_2=\kappa^{(R)}_2+Z_1,\quad T_2=g^{(R)}_2+Z_2.
\end{eqnarray*}
Firstly, from (\ref{eq-17}), we obtain
\begin{eqnarray*}
\frac{dY}{dt}+AY&+&B(U_1,U_1)\Big(\chi_R(|U_1|^2_{\mathcal{W}})-\chi_R(|U_2|^2_{\mathcal{W}})\Big)+B(U_1,Y)\chi_R(|U_2|^2_{\mathcal{W}})
+B(Y,U_2)\chi_R(|U_2|^2_{\mathcal{W}})\\
&&+G(U_1)\Big(\chi_R(|U_1|^2_{\mathcal{W}})-\chi_R(|U_2|^2_{\mathcal{W}})\Big)
+\tilde{G}(Y)\chi_R(|U_2|^2_{\mathcal{W}})=0,
\end{eqnarray*}
where
\[
\tilde{G}(Y)=\left(                 %左括号
  \begin{array}{c}   %该矩阵一共3列，每一列都居中放置
   fk\times \kappa-\int^{z}_{-1}\nabla_H g dz'  \\  % 第一行元素
  \Phi(\kappa)  \\  %第二行元素
  \end{array}
\right).
\]
That is,
\begin{eqnarray}\label{eq-27}
\frac{d\kappa}{dt}&+&A_1\kappa+\Big((v_1\cdot\nabla_H)v_1+\Phi(v_1)\frac{\partial v_1}{\partial z}\Big)\Big(\chi_R(|U_1|^2_{\mathcal{W}})-\chi_R(|U_2|^2_{\mathcal{W}})\Big)
\\ \notag
 &&+\Big((v_1\cdot\nabla_H)\kappa+\Phi(v_1)\frac{\partial \kappa}{\partial z}\Big)\chi_R(|U_2|^2_{\mathcal{W}})+\Big((\kappa\cdot\nabla_H)v_2+\Phi(\kappa)\frac{\partial v_2}{\partial z})\Big)\chi_R(|U_2|^2_{\mathcal{W}})
\\ \notag
&&+(fk\times v_1+\nabla_H p_b-\int^z_{-1}\nabla_HT_1dz')\Big(\chi_R(|U_1|^2_{\mathcal{W}})-\chi_R(|U_2|^2_{\mathcal{W}})\Big)
\\ \notag
&& +(fk\times \kappa-\int^z_{-1}\nabla_H gdz')\chi_R(|U_2|^2_{\mathcal{W}})=0,
\end{eqnarray}
and
\begin{eqnarray} \label{eq-28}
\frac{dg}{dt}&+&A_2g+\Big((v_1\cdot\nabla_H)T_1+\Phi(v_1)\frac{\partial T_1}{\partial z}\Big)\Big(\chi_R(|U_1|^2_{\mathcal{W}})-\chi_R(|U_2|^2_{\mathcal{W}})\Big)\\ \notag
&&+\Big((v_1\cdot\nabla_H)g+\Phi(v_1)\frac{\partial g}{\partial z}\Big)\chi_R(|U_2|^2_{\mathcal{W}})
 +\Big((\kappa\cdot\nabla_H)T_2+\Phi(\kappa)\frac{\partial T_2}{\partial z}\Big)\chi_R(|U_2|^2_{\mathcal{W}})\\ \notag
 &&+\Phi(\kappa)\chi_R(|U_2|^2_{\mathcal{W}})=0.
\end{eqnarray}
Multiplying (\ref{eq-27}) by $-\Lambda^{1+4\varepsilon_0} \kappa$, then integrating over $\mathbb{T}^3$, we have
\begin{eqnarray*}
&&\frac{1}{2}\frac{d}{dt}\int_{\mathbb{T}^3}|\Lambda^{\frac{1}{2}+2\varepsilon_0}  \kappa|^2dxdydz+|\Lambda^{\frac{3}{2}+2\varepsilon_0}  \kappa|^2\\
&=&\Big(\chi_R(|U_1|^2_{\mathcal{W}})-\chi_R(|U_2|^2_{\mathcal{W}})\Big)\int_{\mathbb{T}^3}[(v_1\cdot\nabla_H)v_1+\Phi(v_1)\frac{\partial v_1}{\partial z}]\Lambda^{1+4\varepsilon_0}  \kappa dxdydz\\ \notag
\quad \quad &&+ \chi_R(|U_2|^2_{\mathcal{W}})\int_{\mathbb{T}^3}[(v_1\cdot\nabla_H)\kappa+\Phi(v_1)\frac{\partial \kappa}{\partial z}]\Lambda^{1+4\varepsilon_0}  \kappa dxdydz\\ \notag
\quad \quad &&+\chi_R(|U_2|^2_{\mathcal{W}})\int_{\mathbb{T}^3}[(\kappa\cdot\nabla_H)v_2+\Phi(\kappa)\frac{\partial v_2}{\partial z}]\Lambda^{1+4\varepsilon_0}  \kappa dxdydz\\ \notag
\quad \quad &&+ \Big(\chi_R(|U_1|^2_{\mathcal{W}})-\chi_R(|U_2|^2_{\mathcal{W}})\Big)\int_{\mathbb{T}^3}(fk\times v_1+\nabla_H p_b-\int^z_{-1}\nabla_HT_1dz')\Lambda^{1+4\varepsilon_0}  \kappa dxdydz\\ \notag
\quad \quad &&+\chi_R(|U_2|^2_{\mathcal{W}})\int_{\mathbb{T}^3} \Big(fk\times \kappa -\int^z_{-1}\nabla_H gdz'\Big)\Lambda^{1+4\varepsilon_0}  \kappa dxdydz\\ \notag
&&:=I+II+III+IV+V.
\end{eqnarray*}
It's easy to know that
\[
|\chi_R(|U_1|^2_{\mathcal{W}})-\chi_R(|U_2|^2_{\mathcal{W}})|\leq C(R)|Y|_{\mathcal{W}}[I_{[0,R+1]}(|U_1|^2_{\mathcal{W}})+I_{[0,R+1]}(|U_2|^2_{\mathcal{W}})].
\]
For $I$, since $\int_{\mathbb{T}^3}((v_1\cdot\nabla_H)v_1\Lambda^{1+4\varepsilon_0}  \kappa dxdydz$ is weaker than $\int_{\mathbb{T}^3}(\Phi(v_1)\frac{\partial v_1}{\partial z})\Lambda^{1+4\varepsilon_0}  \kappa dxdydz$, we only need to estimate the term involved $\Phi$. For $\varepsilon_1\in(0,2\varepsilon_0)$, we have
\begin{eqnarray*}
&&C(R)|Y|_{\mathcal{W}}\Big|\int_{\mathbb{T}^3}[\Phi(v_1)\frac{\partial v_1}{\partial z}]\Lambda^{1+4\varepsilon_0}  \kappa dxdydz\Big|\\
&\leq& C(R)|Y|_{\mathcal{W}}\Big|\int_{\mathbb{T}^3}\Lambda^{2\varepsilon_0-\frac{1}{2}+\varepsilon_1}[\Phi(v_1)\frac{\partial v_1}{\partial z}]\Lambda^{\frac{3}{2}+2\varepsilon_0-\varepsilon_1}  \kappa dxdydz\Big|\\
&\leq&
C(R)|Y|_{\mathcal{W}}
\Big(|\Lambda^{2\varepsilon_0+\frac{1}{2}+\varepsilon_1+\sigma_1}v_1||\Lambda^{1+\sigma_2}v_1||\Lambda^{\frac{3}{2}+2\varepsilon_0-\varepsilon_1}\kappa|+
|\Lambda^{1+s_1}v_1||\Lambda^{2\varepsilon_0+\frac{1}{2}+\varepsilon_1+s_2}v_1||\Lambda^{\frac{3}{2}+2\varepsilon_0-\varepsilon_1}\kappa|
\Big)\\
&\leq& C(R)|\Lambda^{\frac{3}{2}+2\varepsilon_0}\kappa||\Lambda^{1+\varepsilon_1}v_1||\Lambda^{\frac{3}{2}+2\varepsilon_0}v_1||\Lambda^{\frac{3}{2}+2\varepsilon_0-\varepsilon_1}\kappa|
+|\Lambda^{\frac{3}{2}+2\varepsilon_0}g||\Lambda^{1+\varepsilon_1}v_1||\Lambda^{\frac{3}{2}+2\varepsilon_0}v_1||\Lambda^{\frac{3}{2}+2\varepsilon_0-\varepsilon_1}\kappa|
\\
&\leq& C(R)\Big(|\Lambda^{\frac{3}{2}+2\varepsilon_0}\kappa|+|\Lambda^{\frac{3}{2}+2\varepsilon_0}g|\Big)|\Lambda^{\frac{1}{2}+2\varepsilon_0}\kappa|^{\varepsilon_1}
|\Lambda^{\frac{3}{2}+2\varepsilon_0}\kappa|^{1-\varepsilon_1}|\Lambda^{\frac{3}{2}+2\varepsilon_0}v_1|^2
\\
&\leq& \varepsilon|\Lambda^{\frac{3}{2}+2\varepsilon_0}\kappa|^2+\varepsilon|\Lambda^{\frac{3}{2}+2\varepsilon_0}g|^2
+C(R)|\Lambda^{\frac{3}{2}+2\varepsilon_0}v_1|^{\frac{4}{\varepsilon_1}}|\Lambda^{\frac{1}{2}+2\varepsilon_0}\kappa|^2\\
&\leq& \varepsilon|\Lambda^{\frac{3}{2}+2\varepsilon_0}\kappa|^2+\varepsilon|\Lambda^{\frac{3}{2}+2\varepsilon_0}g|^2
+C(R,|U_1|^2_{\mathcal{W}})|\Lambda^{\frac{1}{2}+2\varepsilon_0}\kappa|^2,
\end{eqnarray*}
where $\sigma_1+\sigma_2=1$, $s_1+s_2=1$, and we choose
\[
\sigma_1=s_2=\frac{1}{2}-2\varepsilon_0,\quad \sigma_2=s_1=\frac{1}{2}+2\varepsilon_0.
\]

For $II$,
\begin{eqnarray*}
&&\Big|\int_{\mathbb{T}^3}[\Phi(v_1)\frac{\partial \kappa}{\partial z}]\Lambda^{1+4\varepsilon_0}  \kappa dxdydz\Big|\\
&\leq&\Big|\int_{\mathbb{T}^3}\Lambda^{2\varepsilon_0-\frac{1}{2}}  [\Phi(v_1)\frac{\partial \kappa}{\partial z}]\Lambda^{\frac{3}{2}+2\varepsilon_0}  \kappa dxdydz\Big|\\
&\leq&C|\Lambda^{\frac{3}{2}+2\varepsilon_0}\kappa||\Lambda^{2\varepsilon_0+\frac{1}{2}+\sigma_1}v_1||\Lambda^{1+\sigma_2}\kappa|
+C|\Lambda^{\frac{3}{2}+2\varepsilon_0}\kappa||\Lambda^{1+s_1}v_1||\Lambda^{2\varepsilon_0+\frac{1}{2}+s_2}\kappa|\\
&\leq& C|\Lambda^{\frac{3}{2}+2\varepsilon_0}\kappa||\Lambda^{\frac{3}{2}+2\varepsilon_0}v_1||\Lambda \kappa|\\
&\leq& C|\Lambda^{\frac{3}{2}+2\varepsilon_0}\kappa|
|\Lambda^{\frac{3}{2}+\varepsilon_0}v_1||\Lambda^{\frac{1}{2}+2\varepsilon_0}\kappa|^{\frac{1}{2}+2\varepsilon_0}|\Lambda^{\frac{3}{2}+2\varepsilon_0}\kappa|^{\frac{1}{2}-2\varepsilon_0}\\
&\leq& \varepsilon|\Lambda^{\frac{3}{2}+2\varepsilon_0}\kappa|^2+C|\Lambda^{\frac{3}{2}+2\varepsilon_0}v_1|^{\frac{4}{1+4\varepsilon_0}}|\Lambda^{\frac{1}{2}+2\varepsilon_0}\kappa|^2\\
&\leq& \varepsilon|\Lambda^{\frac{3}{2}+2\varepsilon_0}\kappa|^2+C(R,|U_1|^2_{\mathcal{W}})|\Lambda^{\frac{1}{2}+2\varepsilon_0}\kappa|^2,
\end{eqnarray*}
where $\sigma_1+\sigma_2=1$, $s_1+s_2=1$, and we choose
\[
\sigma_1=1,\quad \sigma_2=0,\quad s_1=\frac{1}{2}+2\varepsilon_0,\quad s_2=\frac{1}{2}-2\varepsilon_0.
\]

%\begin{eqnarray*}
%&&\chi_{R}(|U_2|^2_{\mathcal{W}})\int_{\mathbb{T}^3}(\Phi(Z_1)\frac{\partial \kappa}{\partial z})\Lambda^{4\varepsilon_0}  \kappa dxdydz\\
%&\leq& \chi_{R}(|U_2|^2_{\mathcal{W}})\int_{\mathbb{T}^3}\Lambda^{2\varepsilon_0-1}  [\Phi(Z_1)\frac{\partial \kappa}{\partial z}]\Lambda^{1+2\varepsilon_0}  \kappa dxdydz\\
%&\leq& \varepsilon|\Lambda^{1+2\varepsilon_0}\kappa|^2+C|\Lambda^{1+2\varepsilon_0}Z_1|^{\frac{2}{\varepsilon_0}}|\Lambda^{2\varepsilon_0}\kappa|^2,
%\end{eqnarray*}
For $III$,
\begin{eqnarray*}
&&\Big|\int_{\mathbb{T}^3}(\Phi(\kappa)\frac{\partial v_2}{\partial z})\Lambda^{1+4\varepsilon_0}  \kappa dxdydz\Big|\\
&\leq&\int_{\mathbb{T}^3}\Lambda^{2\varepsilon_0-\frac{1}{2}}  [\Phi(\kappa)\frac{\partial v_2}{\partial z}]\Lambda^{\frac{3}{2}+2\varepsilon_0}  \kappa dxdydz\\
&\leq&C|\Lambda^{\frac{3}{2}+2\varepsilon_0}\kappa||\Lambda^{2\varepsilon_0+\frac{1}{2}+\sigma_1}\kappa||\Lambda^{1+\sigma_2}v_2|
+C|\Lambda^{\frac{3}{2}+2\varepsilon_0}\kappa||\Lambda^{1+s_1}\kappa||\Lambda^{2\varepsilon_0+\frac{1}{2}+s_2}v_2|\\
&\leq& C|\Lambda^{\frac{3}{2}+2\varepsilon_0}\kappa||\Lambda \kappa||\Lambda^{\frac{3}{2}+2\varepsilon_0}v_2|\\
&\leq& \varepsilon|\Lambda^{\frac{3}{2}+2\varepsilon_0}\kappa|^2+C(R,|U_2|^2_{\mathcal{W}})|\Lambda^{\frac{1}{2}+2\varepsilon_0}\kappa|^2,
\end{eqnarray*}
where $\sigma_1+\sigma_2=1$, $s_1+s_2=1$, and we choose
\[
\sigma_1=\frac{1}{2}-2\varepsilon_0,\quad \sigma_2=\frac{1}{2}+2\varepsilon_0,\quad s_1=0,\quad s_2=1.
\]

For $IV$,
\begin{eqnarray*}
&&\Big(\chi_R(|U_1|^2_{\mathcal{W}})-\chi_R(|U_2|^2_{\mathcal{W}})\Big)\Big|\int_{\mathbb{T}^3}[fk\times v_1+\nabla_H p_b-\int^z_{-1}\nabla_HT_1dz']\Lambda^{1+4\varepsilon_0}  \kappa dxdydz\Big|\\
&\leq& C(R)[I_{[0,R+1]}(|U_1|^2_{\mathcal{W}})+I_{[0,R+1]}(|U_2|^2_{\mathcal{W}})]|Y|_{\mathcal{W}}\Big|\int_{\mathbb{T}^3}\Lambda^{2\varepsilon_0-\frac{1}{2}+\varepsilon_1}[fk\times v_1\\
&& +\nabla_H p_b-\int^z_{-1}\nabla_HT_1dz']\Lambda^{\frac{3}{2}+2\varepsilon_0-\varepsilon_1}  \kappa dxdydz\Big|\\
&\leq&C(R)[I_{[0,R+1]}(|U_1|^2_{\mathcal{W}})+I_{[0,R+1]}(|U_2|^2_{\mathcal{W}})]|Y|_{\mathcal{W}}
|\Lambda^{\frac{3}{2}+2\varepsilon_0-\varepsilon_1}\kappa||\Lambda^{2\varepsilon_0-\frac{1}{2}+\varepsilon_1}v_1|\\
&& +C(R)[I_{[0,R+1]}(|U_1|^2_{\mathcal{W}})+I_{[0,R+1]}(|U_2|^2_{\mathcal{W}})]|Y|_{\mathcal{W}}
|\Lambda^{\frac{3}{2}+2\varepsilon_0-\varepsilon_1}\kappa||\Lambda^{2\varepsilon_0+\frac{1}{2}+\varepsilon_1}T_1|\\
&:=&C(R)[I_{[0,R+1]}(|U_1|^2_{\mathcal{W}})+I_{[0,R+1]}(|U_2|^2_{\mathcal{W}})](IV^{(1)}+IV^{(2)}),
\end{eqnarray*}
for $IV^{(1)}$,  we have
\begin{eqnarray*}
IV^{(1)}&=&|Y|_{\mathcal{W}}
|\Lambda^{\frac{3}{2}+2\varepsilon_0-\varepsilon_1}\kappa||\Lambda^{2\varepsilon_0-\frac{1}{2}+\varepsilon_1}v_1|\\
&\leq& C|\Lambda^{\frac{3}{2}+2\varepsilon_0}\kappa|^{2-\varepsilon_1}|\Lambda^{\frac{1}{2}+2\varepsilon_0}\kappa|^{\varepsilon_1}|\Lambda^{2\varepsilon_0-\frac{1}{2}+\varepsilon_1}v_1|
+|\Lambda^{\frac{3}{2}+2\varepsilon_0}g||\Lambda^{\frac{1}{2}+2\varepsilon_0}\kappa|^{\varepsilon_1}|\Lambda^{\frac{3}{2}+2\varepsilon_0}\kappa|^{1-\varepsilon_1}|\Lambda^{2\varepsilon_0-\frac{1}{2}+\varepsilon_1}v_1|\\
&\leq& \varepsilon|\Lambda^{\frac{3}{2}+2\varepsilon_0}\kappa|^2+\varepsilon|\Lambda^{\frac{3}{2}+2\varepsilon_0}\kappa|^2+C|\Lambda^{\frac{1}{2}+2\varepsilon_0}\kappa|^2|\Lambda^{2\varepsilon_0-\frac{1}{2}+\varepsilon_1}v_1|^{\frac{2}{\varepsilon_1}}\\
&\leq& \varepsilon|\Lambda^{1+2\varepsilon_0}\kappa|^2+\varepsilon|\Lambda^{1+2\varepsilon_0}\kappa|^2+
C(R,|U_1|^2_{\mathcal{W}})|\Lambda^{\frac{1}{2}+2\varepsilon_0}\kappa|^{2}
\end{eqnarray*}
$IV^{(2)}$ is similar to $IV^{(1)}$, we have
\begin{eqnarray*}
IV^{(2)}&=&|Y|_{\mathcal{W}}
|\Lambda^{\frac{3}{2}+2\varepsilon_0-\varepsilon_1}\kappa||\Lambda^{\frac{1}{2}+2\varepsilon_0+\varepsilon_1}T_1|\\
%&\leq& C|\Lambda^{1+2\varepsilon_0}\kappa|^{2-\varepsilon_1}|\Lambda^{2\varepsilon_0}\kappa|^{\varepsilon_1}|\Lambda^{2\varepsilon_0+\varepsilon_1}(g^{(R)}_1+Z_2)|
%+|\Lambda^{1+2\varepsilon_0}g||\Lambda^{2\varepsilon_0}\kappa|^{\varepsilon_1}|\Lambda^{1+2\varepsilon_0}\kappa|^{1-\varepsilon_1}|\Lambda^{2\varepsilon_0+\varepsilon_1}(g^{(R)}_1+Z_2)|\\
&\leq& \varepsilon|\Lambda^{\frac{3}{2}+2\varepsilon_0}\kappa|^2+\varepsilon|\Lambda^{\frac{3}{2}+2\varepsilon_0}\kappa|^2+C|\Lambda^{\frac{1}{2}+2\varepsilon_0}\kappa|^2|\Lambda^{2\varepsilon_0+\frac{1}{2}+\varepsilon_1}T_1|^{\frac{2}{\varepsilon_1}}\\
&\leq& \varepsilon|\Lambda^{\frac{3}{2}+2\varepsilon_0}\kappa|^2+\varepsilon|\Lambda^{\frac{3}{2}+2\varepsilon_0}\kappa|^2+C(R,|U_1|^2_{\mathcal{W}})|\Lambda^{\frac{1}{2}+2\varepsilon_0}\kappa|^{2}.
\end{eqnarray*}
For $V$, we have
\begin{eqnarray*}
&&\chi_{R}(|U_2|^2_{\mathcal{W}})\Big|\int_{\mathbb{T}^3}\Big(fk\times \kappa-\int^z_{-1}\nabla_H g dz'\Big)\Lambda^{1+4\varepsilon_0}  \kappa dxdydz\Big|\\
&\leq& \chi_{R}(|U_2|^2_{\mathcal{W}})\int_{\mathbb{T}^3}\Lambda^{2\varepsilon_0-\frac{1}{2}}[fk\times \kappa-\int^z_{-1}\nabla_H g dz']\Lambda^{\frac{3}{2}+2\varepsilon_0}  \kappa dxdydz\\
&\leq& \varepsilon|\Lambda^{\frac{3}{2}+2\varepsilon_0}\kappa|^2+C|\Lambda^{2\varepsilon_0-\frac{1}{2}}\kappa|^2+C|\Lambda^{\frac{1}{2}+2\varepsilon_0}g|^2,
\end{eqnarray*}
thus, we have
\begin{eqnarray}\label{eq-33}
&&\frac{1}{2}\frac{d}{dt}\int_{\mathbb{T}^3}|\Lambda^{\frac{1}{2}+2\varepsilon_0} \kappa|^2dxdydz+|\Lambda^{\frac{3}{2}+2\varepsilon_0} \kappa|^2\\ \notag
&\leq& \varepsilon|\Lambda^{\frac{3}{2}+2\varepsilon_0}\kappa|^2
+C(R,|U_1|_{\mathcal{W}},|U_2|_{\mathcal{W}})(1+|\Lambda^{\frac{1}{2}+2\varepsilon_0} \kappa|^2+|\Lambda^{\frac{1}{2}+2\varepsilon_0} g|^2).\notag
\end{eqnarray}
Multiplying (\ref{eq-28}) by $-\Lambda^{1+4\varepsilon_0} g$, then integrating over $\mathbb{T}^3$, we have
 \begin{eqnarray*}
&&\frac{1}{2}\frac{d}{dt}\int_{\mathbb{T}^3}|\Lambda^{\frac{1}{2}+2\varepsilon_0} g|^2dxdydz+|\Lambda^{\frac{3}{2}+2\varepsilon_0} g|^2\\
&&= \Big(\chi_R(|U_1|^2_{\mathcal{W}})-\chi_R(|U_2|^2_{\mathcal{W}})\Big)\int_{\mathbb{T}^3} \Big[(v_1\cdot\nabla_H)T_1+\Phi(v_1)\frac{\partial T_1}{\partial z}\Big]\Lambda^{1+4\varepsilon_0} g dxdydz\\
&&+\chi_R(|U_2|^2_{\mathcal{W}})\int_{\mathbb{T}^3}\Big[(v_1\cdot\nabla_H)g+\Phi(v_1)\frac{\partial g}{\partial z}\Big]\Lambda^{1+4\varepsilon_0}g dxdydz\\
&&+\chi_R(|U_2|^2_{\mathcal{W}})\int_{\mathbb{T}^3}\Big[(\kappa\cdot\nabla_H)T_2+\Phi(\kappa)\frac{\partial T_2}{\partial z}\Big]\Lambda^{1+4\varepsilon_0} g dxdydz\\
&&+\chi_R(|U_2|^2_{\mathcal{W}})\int_{\mathbb{T}^3}\Phi(\kappa)\Lambda^{1+4\varepsilon_0} g dxdydz\\
&&:= VI+VII+VIII+IX.
 \end{eqnarray*}
For $VI$, similar to $I$, we only need to estimate
\begin{eqnarray*}
&&|Y|_{\mathcal{W}}\Big|\int_{\mathbb{T}^3}\Big [\Phi(v_1)\frac{\partial T_1}{\partial z}\Big]\Lambda^{1+4\varepsilon_0} g dxdydz\Big| \\
&=& |Y|_{\mathcal{W}}\Big|\int_{\mathbb{T}^3} \Lambda^{2\varepsilon_0-\frac{1}{2}+\varepsilon_1}[\Phi(v_1)\frac{\partial T_1}{\partial z}]\Lambda^{\frac{3}{2}+2\varepsilon_0-\varepsilon_1} g dxdydz\Big|\\
&\leq& |Y|_{\mathcal{W}}|\Lambda^{\frac{3}{2}+2\varepsilon_0-\varepsilon_1} g|\Big(|\Lambda^{\frac{1}{2}+2\varepsilon_0+\varepsilon_1+\sigma_1} v_1||\Lambda^{1+\sigma_2} T_1|+
|\Lambda^{1+s_1} v_1||\Lambda^{\frac{1}{2}+2\varepsilon_0+\varepsilon_1+s_2}T_1|\Big)\\
&\leq& |Y|_{\mathcal{W}}|\Lambda^{\frac{3}{2}+2\varepsilon_0-\varepsilon_1} g||\Lambda^{1+\varepsilon_1} v_1||\Lambda^{\frac{3}{2}+2\varepsilon_0} T_1|\\
&\leq& C\Big(|\Lambda^{\frac{3}{2}+2\varepsilon_0} \kappa|+|\Lambda^{\frac{3}{2}+2\varepsilon_0} g|\Big)|\Lambda^{\frac{1}{2}+2\varepsilon_0} g|^{\varepsilon_1}|\Lambda^{\frac{3}{2}+2\varepsilon_0} g|^{1-\varepsilon_1}|\Lambda^{1+\varepsilon_1} v_1||\Lambda^{\frac{3}{2}+2\varepsilon_0} T_1|\\
&\leq& \varepsilon|\Lambda^{\frac{3}{2}+2\varepsilon_0} g|^2+\varepsilon|\Lambda^{\frac{3}{2}+2\varepsilon_0} \kappa|^2+C|\Lambda^{1+2\varepsilon_0} v_1|^{\frac{2}{\varepsilon_1}}|\Lambda^{\frac{3}{2}+2\varepsilon_0} T_1|^{\frac{2}{\varepsilon_1}}|\Lambda^{\frac{1}{2}+2\varepsilon_0}g|^2,
\end{eqnarray*}
where $\sigma_1+\sigma_2=1$, $s_1+s_2=1$, and we choose
\[
\sigma_1=\frac{1}{2}-2\varepsilon_0,\quad \sigma_2=\frac{1}{2}+2\varepsilon_0,\quad s_1=\varepsilon_1,\quad s_2=1-\varepsilon_1.
\]

For $VII$,
\begin{eqnarray*}
&&\Big|\int_{\mathbb{T}^3}[\Phi(v_1)\frac{\partial g}{\partial z}]\Lambda^{1+4\varepsilon_0}g dxdydz\Big|\\
&=&\Big|\int_{\mathbb{T}^3}\Lambda^{2\varepsilon_0-\frac{1}{2}}[\Phi(v_1)\frac{\partial g}{\partial z}]\Lambda^{\frac{3}{2}+2\varepsilon_0}g dxdydz\Big|\\
&\leq& C|\Lambda^{\frac{3}{2}+2\varepsilon_0} g|\Big(|\Lambda^{2\varepsilon_0+\frac{1}{2}+\sigma_1} v_1||\Lambda^{1+\sigma_2} g|+|\Lambda^{1+s_1} v_1||\Lambda^{2\varepsilon_0+\frac{1}{2}+s_2} g|\Big)\\
&\leq& C|\Lambda^{\frac{3}{2}+2\varepsilon_0} g||\Lambda^{1} g||\Lambda^{\frac{3}{2}+2\varepsilon_0} v_1|\\
&\leq&C|\Lambda^{\frac{3}{2}+2\varepsilon_0} g||\Lambda^{\frac{1}{2}+2\varepsilon_0} g|^{\frac{1}{2}+2\varepsilon_0}|\Lambda^{\frac{3}{2}+2\varepsilon_0} g|^{\frac{1}{2}-2\varepsilon_0}|\Lambda^{\frac{3}{2}+2\varepsilon_0} v_1|\\
&\leq& \varepsilon|\Lambda^{\frac{3}{2}+2\varepsilon_0} g|^2+C|\Lambda^{\frac{3}{2}+2\varepsilon_0} v_1|^{\frac{4}{1+4\varepsilon_0}}|\Lambda^{\frac{1}{2}+2\varepsilon_0}g|^2,
\end{eqnarray*}
where $\sigma_1+\sigma_2=1$, $s_1+s_2=1$, and we choose
\[
\sigma_1=1,\quad \sigma_2=0, \quad s_1=\frac{1}{2}+2\varepsilon_0 ,\quad s_2=\frac{1}{2}-2\varepsilon_0.
\]
For $VIII$,
\begin{eqnarray*}
&&\Big|\int_{\mathbb{T}^3}[\Phi(\kappa)\frac{\partial T_2}{\partial z}]\Lambda^{1+4\varepsilon_0} g dxdydz\Big|\\
&=& \Big|\int_{\mathbb{T}^3}\Lambda^{2\varepsilon_0-\frac{1}{2}}[\Phi(\kappa)\frac{\partial T_2}{\partial z}]\Lambda^{\frac{3}{2}+2\varepsilon_0} g dxdydz\Big|\\
&\leq& C|\Lambda^{\frac{3}{2}+2\varepsilon_0} g|(|\Lambda^{\frac{1}{2}+2\varepsilon_0+\sigma_1} \kappa||\Lambda^{1+\sigma_2} T_2|+|\Lambda^{1+s_1} \kappa||\Lambda^{\frac{1}{2}+2\varepsilon_0+s_2} T_2|)\\
&\leq& C|\Lambda^{\frac{3}{2}+2\varepsilon_0} g||\Lambda^1 \kappa||\Lambda^{\frac{3}{2}+2\varepsilon_0} T_2|\\
&\leq& C|\Lambda^{\frac{3}{2}+2\varepsilon_0} g||\Lambda^{\frac{1}{2}+2\varepsilon_0} \kappa|^{\frac{1}{2}+2\varepsilon_0}|\Lambda^{\frac{3}{2}+2\varepsilon_0} \kappa|^{\frac{1}{2}-2\varepsilon_0}|\Lambda^{\frac{3}{2}+2\varepsilon_0} T_2|\\
&\leq&\varepsilon|\Lambda^{\frac{3}{2}+2\varepsilon_0} g|^2+\varepsilon|\Lambda^{\frac{3}{2}+2\varepsilon_0} \kappa|^2+C|\Lambda^{\frac{3}{2}+2\varepsilon_0} T_2|^{\frac{4}{1+4\varepsilon_0}}|\Lambda^{\frac{1}{2}+2\varepsilon_0}\kappa|^2,
\end{eqnarray*}
where $\sigma_1+\sigma_2=1$, $s_1+s_2=1$, and we choose
\[
\sigma_1=\frac{1}{2}-2\varepsilon_0,\quad \sigma_2=\frac{1}{2}+2\varepsilon_0,\quad s_1=0,\quad s_2=1.
\]
For $IX$, we have
\begin{eqnarray*}
&&\Big|\int_{\mathbb{T}^3}\Phi(\kappa)\Lambda^{1+4\varepsilon_0}g dxdydz\Big|\\
&=&\Big|\int_{\mathbb{T}^3}\Lambda^{2\varepsilon_0-\frac{1}{2}}\Phi(\kappa)\Lambda^{\frac{3}{2}+2\varepsilon_0}g dxdydz\Big|\\
&\leq& C|\Lambda^{\frac{3}{2}+2\varepsilon_0} g||\Lambda^{\frac{1}{2}+2\varepsilon_0} \kappa|\\
&\leq& \varepsilon|\Lambda^{\frac{3}{2}+2\varepsilon_0} g|^2+C|\Lambda^{\frac{1}{2}+2\varepsilon_0} \kappa|^2,
\end{eqnarray*}
thus, we obtain
\begin{eqnarray}\label{eq-37}
&&\frac{1}{2}\frac{d}{dt}\int_{\mathbb{T}^3}|\Lambda^{\frac{1}{2}+2\varepsilon_0} g|^2dxdydz+|\Lambda^{\frac{3}{2}+2\varepsilon_0} g|^2\\ \notag
&\leq& \varepsilon|\Lambda^{\frac{3}{2}+2\varepsilon_0} g|^2+\varepsilon|\Lambda^{\frac{3}{2}+2\varepsilon_0} \eta|^2+ C(R,|U_1|_{\mathcal{W}},|U_2|_{\mathcal{W}})(1+|\Lambda^{\frac{1}{2}+2\varepsilon_0}\kappa|^2+|\Lambda^{\frac{1}{2}+2\varepsilon_0}g|^2). \notag
\end{eqnarray}
Combining (\ref{eq-33}) with (\ref{eq-37}), we have
\begin{eqnarray}\label{eq-38}
\frac{d\| Y\|^2_{\frac{1}{2}+2\varepsilon_0}}{dt}+\| Y\|^2_{\frac{3}{2}+2\varepsilon_0}
\leq C(R,\sup_{t\in [0,T]}|U_1|_{\mathcal{W}},\sup_{t\in [0,T]}|U_2|_{\mathcal{W}})(1+\| Y\|^2_{\frac{1}{2}+2\varepsilon_0}),
\end{eqnarray}
by Gronwall inequality, (\ref{eq-38}) yields that $\|Y\|_{\frac{1}{2}+2\varepsilon_0 }=0$, which implies $Y=0$.

Up to now, we have proved that equation (\ref{eq-17}) has a unique global weak solution in the space $C([0,T];\mathcal{W})$.

Next, we prove (\ref{eq-39}). In order to do so, it is sufficient to show that $P^{(R)}_y[\tau_R<\varepsilon]\leq C(\varepsilon, R)$ with $C(\varepsilon, R)\downarrow 0$ as $\varepsilon\downarrow 0$, for all $y\in \mathcal{W}$, with $|y|^2_{\mathcal{W}}\leq \frac{R}{8}$.
So, fix $\varepsilon$ small enough, let
\[
\Theta_{\varepsilon, R}:= \sup_{t\in[0,\varepsilon]}|AZ(t)|
\]
and assume that $\Theta^{6}_{\varepsilon, R}\leq \frac{R}{8}$.
Moreover, setting
\[
\varphi(t):= |X^{(R)}|^2_{\mathcal{W}}+\Theta^{6}_{\varepsilon, R},
\]
by (\ref{eq-26}), we get $\dot{\varphi}\leq C(R)$. This implies, together with the bounds on $y$ and $\Theta_{\varepsilon, R}$, that
\[
\sup_{t\in[0,\varepsilon]}|U^{(R)}(t)|^2_{\mathcal{W}}\leq R
\]
for $\varepsilon$ small enough. It follows that $\tau_R\geq \varepsilon$. Hence,
\[
P^{(R)}_y[\tau_R<\varepsilon]\leq P^{(R)}_y\Big[\sup_{t\in [0,\varepsilon]}|A Z|^{6}> \frac{R}{8}\Big],
\]
letting $\varepsilon\downarrow 0$, we have
\[
P^{(R)}_y[\tau_R<\varepsilon]\rightarrow 0.
\]
Since the probability above is independent of $y$, (\ref{eq-39}) is proved. Finally, since
\[
U(t\wedge \tau_R(U^{(R)}))=U^{(R)}(t\wedge \tau_R(U^{(R)})) \quad \forall t,\ P-a.s.
\]
and $U$ is $H$-valued weakly continuous, we obtain $\tau_R(U^{(R)})=\tau_R(U)$, thus (\ref{eq-40}) is proved.

$\hfill\blacksquare$

In order to apply Theorem 5.4 in \cite{F-M} to obtain Theorem \ref{thm-4}, we now only need to prove
\begin{prp}\label{prp-1}
Assume \textbf{Hypothesis H1} holds. For every $R>0$, the transition semigroup $(P^{(R)}_t)_{t\geq 0}$ associated to equation (\ref{eq-16}) is $\mathcal{W}$-strong Feller.
\end{prp}
%Before the proof of Proposition \ref{prp-1}, we give a Lemma firstly, which plays a key role in the following proof.
%\begin{lemma}\cite{Guo}\label{lem-5}
%For any initial data $y\in V$, for (\ref{eq-18}) and (\ref{eq-19}), we have
%\begin{eqnarray*}
%|\frac{\partial \kappa}{\partial z}|^2+|\frac{\partial g}{\partial z}|^2+\int^t_0|\Lambda\frac{\partial \kappa}{\partial z} |^2ds+\int^t_0|\Lambda\frac{\partial g}{\partial z} |^2ds\leq C(t,|\Lambda y|^2).
%\end{eqnarray*}
%\end{lemma}
%Proof of Lemma \ref{lem-5} \quad
%It's easy to verify that $Q^{\frac{1}{2}}=A^{-\frac{3}{4}-\alpha_0}Q^{\frac{1}{2}}_0$ is Hilbert-Schmidt operator from $H$ to $H^{1+2\varepsilon_0}$, which means that $Q^{\frac{1}{2}}$ satisfies the condition on the noise required by Guo and Huang \cite{Guo}, thus, the result of the $H^1$ estimate on $v$ in \cite{Guo} can be used directly here. For $T$, it is similar to $v$, we omit.
%
%$\hfill\blacksquare$
%\\

\textbf{Proof of Proposition \ref{prp-1}} \quad Let $(\Omega, \mathcal{B}, (\mathcal{B}_t)_{t\geq 0}, {P})$ be a filtered probability space, $(W_t)_{t\geq 0}$ a cylindrical Wiener process on $H$ and for every $y\in \mathcal{W}$, denote by $U^{(R)}_y$ the solution to (\ref{eq-16}) with initial value $y\in \mathcal{W}$. By the Bismut-Elworthy-Li formula,
\[
D_z(P^{(R)}_t\psi)(y)=\frac{1}{t}E^{\mathbb{P}}\Big[\psi(U^{(R)}_y(t))\int^t_0(Q^{-\frac{1}{2}}D_zY^{(R)}_y(s),dW(s))\Big],
\]
where $D_z(P^{(R)}_t\psi)$ denotes $(D(P^{(R)}_t\psi),y)$, for $y\in H$, $D_zY^{(R)}_y=DY^{(R)}_y\cdot z$, and $DY^{(R)}_y$ denotes the derivative of $U^{(R)}_y$ with respect to the initial value. Then
for $|\psi|_{\infty}\leq 1$, by the B-D-G inequality, we have
\begin{eqnarray}\label{eq-46}
|(P^{(R)}_t\psi)(y_0+h)-(P^{(R)}_t\psi)(y_0)|\leq \frac{C}{t}\sup_{\eta\in[0,1]}E^{\mathbb{P}}\Big[(\int^t_0|Q^{-\frac{1}{2}}D_hY^{(R)}_{y_0+\eta h}(s)|ds)^{\frac{1}{2}}\Big].
\end{eqnarray}
This proposition is proved once we prove that the right side of the  above inequality converges to 0 as $|h|_{\mathcal{W}}\rightarrow 0$.

For any $y\in \mathcal{W}, h\in H$, write $U=U^{(R)}_y, DY=D_h U=(D_h v, D_h T)=(\eta(t,y)\cdot h,\gamma(t,y)\cdot h)$ and denote $DY=\beta=(\eta,\gamma)$ for simplicity.
Refer to (\ref{eq-16}) and \cite{RR}, we have
\begin{eqnarray}\label{eq-41}
&&\frac{\partial \beta}{\partial t}+\chi_R(|U|^2_{\mathcal{W}})\Big[(v\cdot \nabla_H)\beta+\Phi(v)\frac{\partial \beta}{\partial z}+(\eta\cdot \nabla_H)U+\Phi(\eta)\frac{\partial U}{\partial z}+\tilde{G}(\beta) \Big]
+A\beta\\ \notag
&&\quad \quad +2\chi'_R(|U|^2_{\mathcal{W}})(U,\beta)_{\mathcal{W}}\Big(B(U,U)+G(U)\Big)=0
\end{eqnarray}
with initial value $\beta(0)=h$, and
\[
\tilde{G}(\beta)=\left(                 %左括号
  \begin{array}{c}   %该矩阵一共3列，每一列都居中放置
   fk\times\eta-\int^{z}_{-1}\nabla_H \gamma dz'  \\  %第一行元素
  \Phi(\eta)  \\  %第二行元素
  \end{array}
\right).
\]
We also can rewrite (\ref{eq-41}) in the following form
\begin{eqnarray}\label{eq-42}
&&\frac{\partial \eta}{\partial t}+\chi_R(|U|^2_{\mathcal{W}})[(v\cdot \nabla_H)\eta+\Phi(v)\frac{\partial \eta}{\partial z}+(\eta\cdot \nabla_H)v+\Phi(\eta)\frac{\partial v}{\partial z}+fk\times \eta-\int^z_{-1}\nabla_H\gamma dz' ] \\ \notag
&&\quad +2\chi'_R(|U|^2_{\mathcal{W}})(U,\beta)_{\mathcal{W}}[(v\cdot\nabla_H)v+\Phi(v)\frac{\partial v}{\partial z}+fk\times v+\nabla_H p_b-\int^z_{-1}\nabla_HTdz']-\Delta\eta=0,
\end{eqnarray}
and
\begin{eqnarray}\label{eq-43}
&&\frac{\partial \gamma}{\partial t}+\chi_R(|U|^2_{\mathcal{W}})[(v\cdot \nabla_H)\gamma+\Phi(v)\frac{\partial \gamma}{\partial z}+(\eta\cdot \nabla_H)T+\Phi(\eta)\frac{\partial T}{\partial z}+\Phi(\eta) ] \\ \notag
&&\quad +2\chi'_R(|U|^2_{\mathcal{W}})(U,\beta)_{\mathcal{W}}[(v\cdot\nabla_H)T+\Phi(v)\frac{\partial T}{\partial z}+\Phi(v)]-\Delta\gamma=0.
\end{eqnarray}
Multiplying (\ref{eq-42})  by $-\Lambda^{3+4\varepsilon_0} \eta$, then integrating over $\mathbb{T}^3$, we have
\begin{eqnarray*}
&&\frac{1}{2}\frac{d|\Lambda^{\frac{3}{2}+2\varepsilon_0} \eta|^2}{dt}+|\Lambda^{\frac{5}{2}+2\varepsilon_0} \eta|^2\\
&=&\chi_R(|U|^2_{\mathcal{W}})\int_{\mathbb{T}^3}[(v\cdot \nabla_H)\eta+\Phi(v)\frac{\partial \eta}{\partial z}+(\eta\cdot \nabla_H)v+\Phi(\eta)\frac{\partial v}{\partial z}]\Lambda^{3+4\varepsilon_0}\eta dxdydz\\
&& \quad +\chi_R(|U|^2_{\mathcal{W}})\int_{\mathbb{T}^3}[fk\times \eta-\int^z_{-1}\nabla_H\gamma dz' ]\Lambda^{3+4\varepsilon_0}\eta dxdydz\\ \notag
&&\quad+2\int_{\mathbb{T}^3}\chi'_R(|U|^2_{\mathcal{W}})(U,\beta)_{\mathcal{W}}[(v\cdot\nabla_H)v+\Phi(v)\frac{\partial v}{\partial z}]\Lambda^{3+4\varepsilon_0} \eta dxdydz
\\ \notag
&& \quad+2\int_{\mathbb{T}^3}\chi'_R(|U|^2_{\mathcal{W}})(U,\beta)_{\mathcal{W}}[fk\times v+\nabla_H p_b-\int^z_{-1}\nabla_HTdz']\Lambda^{3+4\varepsilon_0} \eta dxdydz\\
&:=&  K_1+K_2+K_3+K_4.
\end{eqnarray*}
For $K_1$, we only need to estimate
\begin{eqnarray*}
&&\Big|\int_{\mathbb{T}^3}[\Phi(v)\frac{\partial \eta}{\partial z}+\Phi(\eta)\frac{\partial v}{\partial z}]\Lambda^{3+4\varepsilon_0}\eta dxdydz\Big| \\
&=&\Big|\int_{\mathbb{T}^3}\Lambda^{\frac{1}{2}+2\varepsilon_0}[\Phi(v)\frac{\partial \eta}{\partial z}+\Phi(\eta)\frac{\partial v}{\partial z}]\Lambda^{\frac{5}{2}+2\varepsilon_0}\eta dxdydz\Big|\\
&\leq& |\Lambda^{\frac{5}{2}+2\varepsilon_0}\eta |\Big(|\Lambda^{\frac{3}{2}+2\varepsilon_0+\sigma_1}\eta||\Lambda^{1+\sigma_2}v|+|\Lambda^{1+s_1}\eta||\Lambda^{\frac{3}{2}+2\varepsilon_0+s_2}v|\Big)\\
&\leq& |\Lambda^{\frac{5}{2}+2\varepsilon_0}\eta ||\Lambda^2\eta||\Lambda^{\frac{3}{2}+2\varepsilon_0}v|\\
&\leq& |\Lambda^{\frac{5}{2}+2\varepsilon_0}\eta ||\Lambda^{\frac{3}{2}+2\varepsilon_0}\eta|^{\frac{1}{2}+2\varepsilon_0}
|\Lambda^{\frac{5}{2}+2\varepsilon_0}\eta|^{\frac{1}{2}-2\varepsilon_0}|\Lambda^{\frac{3}{2}+2\varepsilon_0}v|\\
&\leq& \varepsilon|\Lambda^{\frac{5}{2}+2\varepsilon_0}\eta |^2+C|\Lambda^{\frac{3}{2}+2\varepsilon_0}\eta |^2|\Lambda^{\frac{3}{2}+2\varepsilon_0}v|^{\frac{4}{1+4\varepsilon_0}},
\end{eqnarray*}
where $\sigma_1+\sigma_2=1$, $s_1+s_2=1$, and we choose
\[
\sigma_1=\frac{1}{2}-2\varepsilon_0,\quad \sigma_2=\frac{1}{2}+2\varepsilon_0,\quad s_1=1,\quad s_2=0.
\]
For $K_2$,
\begin{eqnarray*}
&&\Big|\int_{\mathbb{T}^3}[fk\times \eta-\int^z_{-1}\nabla_H\gamma dz' ]\Lambda^{3+4\varepsilon_0}\eta dxdydz \Big|\\
&\leq&\varepsilon|\Lambda^{\frac{5}{2}+2\varepsilon_0}\eta |^2+C|\Lambda^{\frac{1}{2}+2\varepsilon_0}\gamma |^2+C|\Lambda^{\frac{3}{2}+2\varepsilon_0}\eta|^2.
\end{eqnarray*}
For $K_3$,
\begin{eqnarray*}
&&\Big|\int_{\mathbb{T}^3}(U,\beta)_{\mathcal{W}}[\Phi(v)\frac{\partial v}{\partial z}]\Lambda^{3+4\varepsilon_0} \eta dxdydz\Big|\\
&\leq &(U,\beta)_{\mathcal{W}}|\Lambda^{\frac{5}{2}+2\varepsilon_0}\eta |\Big(|\Lambda^{\frac{3}{2}+2\varepsilon_0+\sigma_1}v ||\Lambda^{\sigma_2}\frac{\partial v}{\partial z} |+|\Lambda^{1+s_1}v ||\Lambda^{\frac{1}{2}+2\varepsilon_0+s_2}\frac{\partial v}{\partial z} |\Big)\\
&\leq& C(U,\beta)_{\mathcal{W}}|\Lambda^{\frac{5}{2}+2\varepsilon_0}\eta ||\Lambda^{\frac{3}{2}+2\varepsilon_0}v ||\Lambda \frac{\partial v}{\partial z}|\\
&\leq& C|\Lambda^{\frac{3}{2}+2\varepsilon_0}v ||\Lambda^{\frac{3}{2}+2\varepsilon_0}\eta ||\Lambda^{\frac{5}{2}+2\varepsilon_0}\eta ||\Lambda^{\frac{3}{2}+2\varepsilon_0}v ||\Lambda \frac{\partial v}{\partial z}|+C|\Lambda^{\frac{3}{2}+2\varepsilon_0}T ||\Lambda^{\frac{3}{2}+2\varepsilon_0}\gamma ||\Lambda^{\frac{5}{2}+2\varepsilon_0}\eta ||\Lambda^{\frac{3}{2}+2\varepsilon_0}v ||\Lambda \frac{\partial v}{\partial z}|\\
&\leq & \varepsilon|\Lambda^{\frac{5}{2}+2\varepsilon_0}\eta |^2+C|\Lambda^{\frac{3}{2}+2\varepsilon_0}v |^4|\Lambda^{\frac{3}{2}+2\varepsilon_0}\eta |^2(|\Lambda^{\frac{5}{2}+2\varepsilon_0}\kappa|^2+|\Lambda\frac{\partial Z_1}{\partial z}|^2)\\
&&\quad +C|\Lambda^{\frac{3}{2}+2\varepsilon_0}v |^2|\Lambda^{\frac{3}{2}+2\varepsilon_0}T |^2|\Lambda^{\frac{3}{2}+2\varepsilon_0}\gamma |^2(|\Lambda^{\frac{5}{2}+2\varepsilon_0}\kappa|^2+|\Lambda\frac{\partial Z_1}{\partial z}|^2),
\end{eqnarray*}
where $\sigma_1+\sigma_2=1$, $s_1+s_2=1$, and we choose
\[
\sigma_1=0,\quad \sigma_2=1,\quad s_1=\frac{1}{2}+2\varepsilon_0,\quad s_2=\frac{1}{2}-2\varepsilon_0.
\]
For $K_4$,
\begin{eqnarray*}
&&\Big|\int_{\mathbb{T}^3}(U,\beta)_{\mathcal{W}}[\int^z_{-1}\nabla_HTdz']\Lambda^{3+4\varepsilon_0} \eta dxdydz\Big|\\
&=& \Big|\int_{\mathbb{T}^3}(U,\beta)_{\mathcal{W}}\Lambda^{\frac{1}{2}+2\varepsilon_0}[\int^z_{-1}\nabla_HTdz']\Lambda^{\frac{5}{2}+2\varepsilon_0} \eta dxdydz\Big|\\
&\leq & (U,\beta)_{\mathcal{W}}|\Lambda^{\frac{5}{2}+2\varepsilon_0}\eta ||\Lambda^{\frac{3}{2}+2\varepsilon_0}T|\\
&\leq &C|\Lambda^{\frac{3}{2}+2\varepsilon_0}v ||\Lambda^{\frac{3}{2}+2\varepsilon_0}\eta ||\Lambda^{\frac{5}{2}+2\varepsilon_0}\eta ||\Lambda^{\frac{3}{2}+2\varepsilon_0}T|+C|\Lambda^{\frac{3}{2}+2\varepsilon_0}T ||\Lambda^{\frac{3}{2}+2\varepsilon_0}\gamma ||\Lambda^{\frac{5}{2}+2\varepsilon_0}\eta ||\Lambda^{\frac{3}{2}+2\varepsilon_0}T|\\
&\leq& \varepsilon|\Lambda^{\frac{5}{2}+2\varepsilon_0}\eta |^2+C|\Lambda^{\frac{3}{2}+2\varepsilon_0}v |^2|\Lambda^{\frac{3}{2}+2\varepsilon_0}T |^2|\Lambda^{\frac{3}{2}+2\varepsilon_0}\eta |^2+C|\Lambda^{\frac{3}{2}+2\varepsilon_0}T |^4|\Lambda^{\frac{3}{2}+2\varepsilon_0}\gamma |^2.
\end{eqnarray*}
Multiplying (\ref{eq-43})  by $-\Lambda^{3+4\varepsilon_0} \gamma$, then integrating over $\mathbb{T}^3$, we have
\begin{eqnarray*}
&&\frac{1}{2}\frac{d |\Lambda^{\frac{3}{2}+2\varepsilon_0}\gamma|^2}{d t}+|\Lambda^{\frac{5}{2}+2\varepsilon_0} \gamma|^2\\
&=&\chi_R(|U|^2_{\mathcal{W}})\int_{\mathbb{T}^3}\Big[(v\cdot \nabla_H)\gamma+\Phi(v)\frac{\partial \gamma}{\partial z}+(\eta\cdot \nabla_H)T+\Phi(\eta)\frac{\partial T}{\partial z} +\Phi(\eta)\Big]\Lambda^{3+4\varepsilon_0} \gamma dxdydz \\
&&\quad +2\int_{\mathbb{T}^3}\chi'_R(|U|^2_{\mathcal{W}})(U,\beta)_{\mathcal{W}}\int_{\mathbb{T}^3}\Big[(v\cdot\nabla_H)T+\Phi(v)\frac{\partial T}{\partial z}+\Phi(v)\Big]\Lambda^{3+4\varepsilon_0} \gamma dxdydz\\
&:=& K_5+K_6.
\end{eqnarray*}
For $K_5$, we only need to estimate the following term,
\begin{eqnarray*}
&&\Big|\int_{\mathbb{T}^3}[\Phi(v)\frac{\partial \gamma}{\partial z}+\Phi(\eta)\frac{\partial T}{\partial z} ]\Lambda^{3+4\varepsilon_0} \gamma dxdydz\Big| \\
&=& \Big|\int_{\mathbb{T}^3}\Lambda^{\frac{1}{2}+2\varepsilon_0} [\Phi(v)\frac{\partial \gamma}{\partial z}+\Phi(\eta)\frac{\partial T}{\partial z} ]\Lambda^{\frac{5}{2}+2\varepsilon_0} \gamma dxdydz\Big|\\
&\leq& C|\Lambda^{\frac{5}{2}+2\varepsilon_0} \gamma|\Big(|\Lambda^{\frac{3}{2}+2\varepsilon_0} v||\Lambda^{2} \gamma|+|\Lambda^{\frac{3}{2}+2\varepsilon_0} T||\Lambda^{2} \eta|\Big)\\
&\leq& \varepsilon|\Lambda^{\frac{5}{2}+2\varepsilon_0}\eta |^2+\varepsilon|\Lambda^{\frac{5}{2}+2\varepsilon_0}\gamma |^2
+C|\Lambda^{\frac{3}{2}+2\varepsilon_0}\gamma |^2|\Lambda^{\frac{3}{2}+2\varepsilon_0}v |^{\frac{4}{1+4\varepsilon_0}}+C|\Lambda^{\frac{3}{2}+2\varepsilon_0}\eta |^2|\Lambda^{\frac{3}{2}+2\varepsilon_0}T |^{\frac{4}{1+4\varepsilon_0}}.
\end{eqnarray*}
For $K_6$,
\begin{eqnarray*}
&&\Big|\int_{\mathbb{T}^3}\chi'_R(|U|^2_{\mathcal{W}})(U,\beta)_{\mathcal{W}}\int_{\mathbb{T}^3}[\Phi(v)\frac{\partial T}{\partial z}]\Lambda^{3+4\varepsilon_0} \gamma dxdydz\Big|\\
&=& \Big|\int_{\mathbb{T}^3}\chi'_R(|U|^2_{\mathcal{W}})(U,\beta)_{\mathcal{W}}\int_{\mathbb{T}^3}\Lambda^{\frac{1}{2}+2\varepsilon_0}[\Phi(v)\frac{\partial T}{\partial z}]\Lambda^{\frac{5}{2}+2\varepsilon_0} \gamma dxdydz\Big|\\
&\leq& C(U,\beta)_{\mathcal{W}}|\Lambda^{\frac{5}{2}+2\varepsilon_0} \gamma|\Big(|\Lambda^{\frac{3}{2}+2\varepsilon_0+\sigma_1}v||\Lambda^{1+\sigma_2} T|+|\Lambda^{1+s_1}v||\Lambda^{\frac{1}{2}+2\varepsilon_0+s_2} \frac{\partial T}{\partial z}|\Big)\\
&\leq& C(U,\beta)_{\mathcal{W}}|\Lambda^{\frac{5}{2}+2\varepsilon_0} \gamma||\Lambda^{\frac{3}{2}+2\varepsilon_0} v||\Lambda \frac{\partial T}{\partial z}|\\
&\leq& C|\Lambda^{\frac{3}{2}+2\varepsilon_0} v|^2|\Lambda^{\frac{3}{2}+2\varepsilon_0} \eta||\Lambda \frac{\partial T}{\partial z}||\Lambda^{\frac{5}{2}+2\varepsilon_0} \gamma|+C|\Lambda^{\frac{3}{2}+2\varepsilon_0} v||\Lambda^{\frac{3}{2}+2\varepsilon_0} T||\Lambda^{\frac{3}{2}+2\varepsilon_0} \gamma||\Lambda \frac{\partial T}{\partial z}||\Lambda^{\frac{5}{2}+2\varepsilon_0} \gamma|\\
&\leq& C|\Lambda^{\frac{3}{2}+2\varepsilon_0} v|^2|\Lambda^{\frac{3}{2}+2\varepsilon_0} \eta||\Lambda^{\frac{5}{2}+2\varepsilon_0} g||\Lambda^{\frac{5}{2}+2\varepsilon_0} \gamma|+C|\Lambda^{\frac{3}{2}+2\varepsilon_0} v||\Lambda^{\frac{3}{2}+2\varepsilon_0} T||\Lambda^{\frac{3}{2}+2\varepsilon_0} \gamma||\Lambda^{\frac{5}{2}+2\varepsilon_0} g||\Lambda^{\frac{5}{2}+2\varepsilon_0} \gamma|\\
&&\quad +C|\Lambda^{\frac{3}{2}+2\varepsilon_0} v|^2|\Lambda^{\frac{3}{2}+2\varepsilon_0} \eta||\Lambda \frac{\partial Z_2}{\partial z}||\Lambda^{\frac{5}{2}+2\varepsilon_0} \gamma|+C|\Lambda^{\frac{3}{2}+2\varepsilon_0} v||\Lambda^{\frac{3}{2}+2\varepsilon_0} T||\Lambda^{\frac{3}{2}+2\varepsilon_0} \gamma||\Lambda \frac{\partial Z_2}{\partial z}||\Lambda^{\frac{5}{2}+2\varepsilon_0} \gamma|\\
&\leq& \varepsilon |\Lambda^{\frac{5}{2}+2\varepsilon_0} \gamma|^2+C|\Lambda^{\frac{3}{2}+2\varepsilon_0} v|^4|\Lambda^{\frac{3}{2}+2\varepsilon_0} \eta|^2(|\Lambda^{\frac{5}{2}+2\varepsilon_0} g|^2+|\Lambda\frac{\partial Z_2}{\partial z}|^2)\\
&&\quad+C|\Lambda^{\frac{3}{2}+2\varepsilon_0} v|^2|\Lambda^{\frac{3}{2}+2\varepsilon_0} T|^2|\Lambda^{\frac{3}{2}+2\varepsilon_0} \gamma|^2(|\Lambda^{\frac{5}{2}+2\varepsilon_0} g|^2+|\Lambda\frac{\partial Z_2}{\partial z}|^2).
\end{eqnarray*}
Combing all the above estimations, we have for $|U|^2_{\mathcal{W}}\leq R$,
\begin{eqnarray*}
\frac{d| \beta|^2_{\mathcal{W}}}{dt}+\|\beta\|^2_{\frac{5}{2}+2\varepsilon_0}
\leq C\Big(C(R)+|\Lambda^{\frac{5}{2}+2\varepsilon_0} \kappa|^2+|\Lambda^{\frac{5}{2}+2\varepsilon_0} g|^2+|\Lambda\frac{\partial Z_1}{\partial z}|^2+|\Lambda\frac{\partial Z_2}{\partial z}|^2\Big)|\beta|^2_{\mathcal{W}},
\end{eqnarray*}
by Gronwall's inequality  and (\ref{eq-26}), we finally get
\begin{eqnarray*}
&&\int^t_0\|\beta(l)\|^2_{\frac{5}{2}+2\varepsilon_0}dl\\
&\leq& C| h|^2_{\mathcal{W}}+\exp\left(C\int^t_0\Big(C(R)+|\Lambda^{\frac{5}{2}+2\varepsilon_0} \kappa|^2+|\Lambda^{\frac{5}{2}+2\varepsilon_0} g|^2+|\Lambda\frac{\partial Z_1}{\partial z}|^2+|\Lambda\frac{\partial Z_2}{\partial z}|^2\Big)dl\right)| h|^2_{\mathcal{W}}\\
&\leq& C| h|^2_{\mathcal{W}}+\exp\left(C\left(| \Lambda^{\frac{3}{2}+2\varepsilon_0}y|^2+\int^t_0\Big(C(R)+|\Lambda\frac{\partial Z_1}{\partial z}|^2+|\Lambda\frac{\partial Z_2}{\partial z}|^2\Big)dl\right)\right)| h|^2_{\mathcal{W}}.
\end{eqnarray*}
Since $Z_i$ is a Gaussian random variable in $C([0,\infty), D(\Lambda^{2+2\varepsilon_0-2\varepsilon}))$, for every $\varepsilon>0$, by Fernique's theorem, we could choose $t_0$ small enough and obtain
\[
\mathbb{E}\int^{t_0}_0\|\beta(l)\|^2_{\frac{5}{2}+2\varepsilon_0}dl\leq C(t_0,R)| h|^2_{\mathcal{W}}.
\]
From \textbf{Hypothesis H1}, it follows that
\[
Q^{-\frac{1}{2}}=Q^{-\frac{1}{2}}_0A^{\frac{5}{4}+\varepsilon_0},
\]
thus, the assertion of (\ref{eq-46}) holds for $t_0$. For general $t$, by the semigroup property, the assertion follows easily.
$\hfill\blacksquare$
%\section{Ergodicity for Markov Selection}
\begin{remark}
In order to apply Fernique's theorem to Gaussian process $Z$, we have to make use of $\Lambda \frac{\partial Z}{\partial z}$ and control its power to be less than 2 during the estimate of $K_3$ and $K_6$.
\end{remark}

\textbf{Proof of Theorem \ref{thm-4}} \quad Theorem \ref{thm-5} and Proposition \ref{prp-1} imply the result by Theorem 5.4 in \cite{F-M}.

\vskip0.5cm {\small {\bf  Acknowledgements}\ \  This work was supported by National Natural Science Foundation of China (NSFC) (No. 11431014, No. 11371041, No. 11401557, No. 11271356)}, the Fundamental Research Funds for the Central Universities (No. 0010000048),
Key Laboratory of Random Complex Structures and Data Science, Academy of Mathematics and Systems Science, Chinese Academy of  Sciences(No. 2008DP173182), and the applied mathematical research for the important strategic demand of China in information science and related fields(No. 2011CB808000).

%%%%%%%%%%%%%%%%%%%%%%%%%%%%%%%%%%%%%%%%%%%%%%%%%%%%%%%%%%%%%%%%%%%%%%%%%%%%%%%%%%%%%%%%%%%%%%%%%%%%%%%%%%%%%%%
%\vskip 0.2cm {\small {\bf  Acknowledgements}\   This work was
%partly supported by  NSFC(No.11071008, No.11101419) and the Doctoral Program
%Foundation of the Ministry of Education, China. in Germany }

%%%%%%%%%%%%%%%%%%%%%%%%%%%%%%%%%%%%%%%%%%%%%%%%%%%%%%
\def\refname{ References}


\begin{thebibliography}{2}



\bibitem{Adams}R.A. Adams : \emph{Sobolev Space}. New York: Academic Press, 1975.
\bibitem {C-L-T} C. Cao, J. Li, E. Titi: \emph{Local and global well-posedness of strong solutions to the 3D primitive equations with vertical eddy diffusivity}. Arch. Rational Mech. Anal. 214 35-76 (2014).
\bibitem{C-T-2} C. Cao, E.S. Titi: \emph{Global well-posedness and finite-dimensional global attractor for a 3D planetary geostrophic viscous model.} Comm. Pure Appl. Math. 56, no. 2, 198-233 (2003).
\bibitem{C-T-1} C. Cao, E.S. Titi: \emph{Global well-posedness of the three-dimensional viscous primitive equations of large-scale ocean and atmosphere dynamics}. Ann. of Math. 166, 245-267 (2007).


\bibitem{D-G-T-Z} A. Debussche,  N. Glatt-Holtz, R. Temam, M. Ziane: \emph{Global existence and regularity for the 3D stochastic primitive equations of the ocean and atmosphere with
multiplicative white noise.}  Nonlinearity, 25, 2093-2118 (2012).
\bibitem{RR} Z. Dong, J. Zhai, R. Zhang: \emph{Exponential mixing for 3D stochastic primitive equations of the large scale ocean.} Preprint. Available at arXiv: 1506.08514.
\bibitem{B-M-R} B. Ewald, M. Petcu, P. Temam: \emph{Stochstic solutions of the two-dimensional primitive equations of the ocean and atmosphere with an additive noise}, Anal. Appl., 5, 183-198 (2007).
\bibitem{F-M} F. Flandoli, M. Romito :\emph{ Markov selection for the 3D stochastic Navier-Stokes equations} Probab. Theory Relat. Fields 140: 407-458 (2008).
\bibitem{G-S-1} H. Gao, C. Sun: \emph{Hausdorff dimension of random attractor for stochastic
Navier-Stokes-Voight equations and primitive equations}. Dyn. Partial Differ. Equ. 7 , no. 4, 307-326 (2010).
\bibitem{G-S} H. Gao, C. Sun: \emph{Well-posedness and large deviations for the stochastic primitive equations in two space dimensions}. Commun. Math. Sci. Vol.10, No.2, 575-593 (2012).

\bibitem{B-M-X} B. Goldys, M. R\"{o}ckner, X. Zhang: \emph{
Martingale solutions and Markov selections for stochastic partial differential equations}.
Stochastic Process. Appl. 119, no. 5, 1725-1764 (2009).
\bibitem{GG-MN} F. $\rm Guill\acute{e}n-Gonz\acute{a}lez$, N. Masmoudi, M.A. Rodriguez-Bellido: \emph{Anisotropic estimates and strong solutions
for the primitive equations.} Diff. Int. Equ. 14, 1381-1408(2001).
\bibitem{Guo} B. Guo, D. Huang: \emph{3D stochastic primitive equations of the large-scale ocean: global well-posedness and attractors.}  Comm. Math. Phys. 286, no. 2, 697-723 (2009).
\bibitem{H-T-Z} C. Hu, R. Temam, M. Ziane: \emph{The primitive equations on the large scale ocean under the small depth hypothesis.}
   Discrete Contin. Dynam. Systems 9, 97-131 (2003).
\bibitem{Ju} N. Ju: \emph{Existence and uniqueness of the solution to the dissipative 2D quasi-geostrophic equations in the Sobolev space}. Comm. Math. Phys. 251 365-376. MR2100059 (2004).
    \bibitem{L-T} J. Li, E.S. Titi:\emph{ Existence and uniqueness of weak solutions to viscous primitive equations for certain class of discontinuous initial data}. Arxiv:1512.00700 (2015).
\bibitem{L-T-W-1} J.L. Lions, R. Temam, S. Wang : \emph{New formulations of the primitive equations of atmosphere and applications.} Nonlinearity 5, 237-288 (1992).
\bibitem{L-T-W-2}  J.L. Lions, R. Temam, S. Wang : \emph{On the equations of the large scale ocean.} Nonlinearity 5, 1007-1053 (1992).
\bibitem{L-T-W-3}  J.L. Lions, R. Temam, S. Wang : \emph{Models of the coupled atmosphere and ocean.} Computational Mechanics Advance 1, 1-54 (1993).
\bibitem{L-T-W-4}  J.L. Lions, R. Temam, S. Wang : \emph{Mathematical theory for the coupled atmosphere-ocean models.}  J. Math. Pures Appl. 74, 105-163 (1995).
\bibitem{Resnick} S.G. Resnick: \emph{Dynamical problems in nin-linear advective partial differential equations}. PH.D. thesis, Univ. of Chicago. MR2716577 (1995).
\bibitem{ROMITO} M. Romito : \emph{Analysis of equilibrium states of markov solutions to the 3D Navier-Stokes equations driven by additive noise} J. Stat. Phys. 131, no.3, 415-444 (2008).

\bibitem{Stein} E. M. Stein: \emph{Singular Integrals and Differentiability Properties of Functions}. Princeton Univ. Press, Princeton, N.J. MR0290095 (1970).
 \bibitem{S-V} D.W. Stroock, S.R.S. Varadhan: \emph{Multidimensional Diffusion Processes}. Springer-Verlag, Berlin, 1979.
\bibitem{T-R} R. Temam : \emph{Navier-Stokes Equations: Theory and Numeriacal Analysis,} 3rd ed. Studies in Mathematical and Its Applications 2. North-Holland, Amsterdam (1984).
\bibitem{Z}	G. Zhou: \emph{Global well-posedness and random attractor
of the 3D viscous primitive equations driven
by fractional noises}. Preprint. Available at arXiv:1604.05376 (2016).

%\bibitem{C} Chavel, I. \emph{Eigenvalues in Riemannian geometry.} Pure and Applied Mathematics, 115. Academic, Orlando, Fla., 1984.
%\bibitem{D-Z-1} Da Prato, G., Zabczyk, J.: \emph{Stochastic equations in infinite dimensions.} Encyclopedia of mathematics and its Applications, Cambridge University Press, Cambridge (1992).
% \bibitem{D-Z-2} Da Prato, G., Zabczyk, J.: \emph{Ergodicity for Inifinite  Dimensional Systems.} In: London Mathematical Society Lecture Notes,n.229, Cambridge University Press, Cambridge (1996)
%\bibitem{D-G-T} Debussche, A.,  Glatt-Holtz, N.,  Temam, R.: \emph {Local martingale and pathwise solutions for an abstract fluids model.} Phys. D 240 , no. 14-15, 1123-1144 (2011)

%
%\bibitem{F-G} Flandoli, F., Gatarek, D.: \emph{Martingale and stationary solutions for stochastic Navier-Stokes equations.} Probab. Theory Related Fields 102, 367-391 (1995)
%
%
%
%\bibitem{F-H}Frankignoul, C., Hasselmann, K.: \emph{Stochastic climate models, Part II: Application to sea-surface temperature anomalies and thermocline variability.} Tellus 29, 289-305 (1977)
%
%
%\bibitem{GKVZ}Glatt-Holtz. N., Kukavica. I., Vicol. V., Ziane. M.: \emph{Existence and regularity of invariant measures for the three dimensional stochastic primitive equations.} arXiv:1311.4204v1
%

%    \bibitem{Hairer} Hairer, M.: \emph{Exponential Mixing Properties of stochastic PDEs through Asymptotic coupling.} Prob. Theory Rel. Fields 124(3), 345-380 (2002)



%
%
%
%   \bibitem{K} Kuksin, S.:  \emph{On exponential convergence to a stationary measure for nonlinear PDEs.} Partial differential equations, 161-176 (2002)
%    \bibitem{K-S} Kuksin, S., Shirikyan, A.: \emph{Ergodicity for the randomly forced 2D Navier-Stokes equations.} Math. Phys. Anal. Geom. 4, 147-195 (2001)
%        \bibitem{K-S-1} Kuksin, S., Shirikyan, A.: \emph{A coupling approach to randomly forced PDE's I.} Comm. Math. Phys. 221, no. 2, 351-366 (2001)
%             \bibitem{K-S-2} Kuksin, S., Piatnitski, A., Shirikyan, A.: \emph{A coupling approach to randomly forced PDE's II.} Comm. Math. Phys. 230(1), 81-85 (2002)

%     \bibitem{M} Mattingly, J.: \emph{Exponential convergence for the stochastically forced Navier-Stokes equations and other partially dissipative dynamics.} Comm. Math. Phys. 230, 421-462 (2002)
%
%     \bibitem{M-M}Mikolajewicz, U., Maier-Reimer, E.: \emph{Internal secular variability in an OGCM.} Climate Dyn. 4, 145-156 (1990)
%
%
%\bibitem{M-1} Tachim Medjo, T.: \emph{The exponential behavior of the stochastic three-dimensional primitive equations with multiplicative noise.} Nonlinear Anal. Real World Appl. 12, no. 2, 799-810 (2011)
%\bibitem{Nualart}Nualart, D.: \emph{ Malliavin Calculus and related topic.}  Probability and its Applications, 93B Newyork: Springer, 1995
%\bibitem{O-C} Odasso, C.: \emph{Exponential mixing for the 3D stochastic Navier-Stokes equations.}  Comm. Math. Phys. 270, no. 1, 109-139 (2007)
%\bibitem{O-C-1} Odasso, C.: \emph{Ergodicity for the stochastic Complex Ginzburg-Landau equations.} Ann. Inst. H. Poincaré Probab. Statist. 42, no. 4, 417-454 (2006)
%
%  \bibitem{P}Phillips, O.M.: \emph{On the generation of waves by turbulent winds.} J. Fluid Mech. 2, 417-445 (1957)
%
%%\bibitem{R} Romito, M. \emph{Analysis of equilibrium states of markov solutions to the 3D Navier-stokes equations driven by additive noise}.  J. Stat. Phys. 131, no. 3, 415-444 (2008).
%\bibitem{Temam} Temam, R.: \emph{Navier-Stokes equations and nonlinear functional analysis.} Second edition. CBMS-NSF Regional Conference Series in Applied Mathematics, 66. Society for Industrial and Applied Mathematics (SIAM), Philadelphia, PA, 1995


%\bibitem{Sato} Sato, K.I. \emph{L\'evy processes and infinite divisible distributions,} Cambridge University Press. London, 1999.



\end{thebibliography}
\end{document}